\documentclass[11pt,reqno]{amsart}
\makeatletter
\@namedef{subjclassname@2020}{\textup{2020} Mathematics Subject Classification}
\makeatother
\usepackage{amssymb,amscd,amsthm}
\usepackage[T1]{fontenc}
\usepackage{blkarray}

\usepackage{amsmath,amssymb,graphicx,mathrsfs} 
\usepackage[colorlinks=true,urlcolor=blue,citecolor=red,linkcolor=blue]{hyperref}

\let\frak\mathfrak
\let\Bbb\mathbb

\def\>{\relax\ifmmode\mskip.666667\thinmuskip\relax\else\kern.111111em\fi}
\def\<{\relax\ifmmode\mskip-.333333\thinmuskip\relax\else\kern-.0555556em\fi}
\def\vsk#1>{\vskip#1\baselineskip}
\def\vv#1>{\vadjust{\vsk#1>}\ignorespaces}
\def\vvn#1>{\vadjust{\nobreak\vsk#1>\nobreak}\ignorespaces}

  \let\ssize\scriptstyle
\let\sssize\scriptscriptstyle

\let\Medskip\medskip
\def\medskip{\par\Medskip}
\let\Bigskip\bigskip
\def\bigskip{\par\Bigskip}

\let\Maketitle\maketitle
\def\maketitle{\Maketitle\thispagestyle{empty}\let\maketitle\empty}

\newtheorem{thm}{Theorem}[section]
\newtheorem{cor}[thm]{Corollary}
\newtheorem{lem}[thm]{Lemma}
\newtheorem{prop}[thm]{Proposition}

\numberwithin{equation}{section}

\theoremstyle{definition}
\newtheorem{defn}[thm]{Definition}
\newtheorem{rem}[thm]{Remark}
\newtheorem{prob}[thm]{Problem}
\newtheorem*{example}{Example}

\let\mc\mathcal
\let\nc\newcommand

\let\al\alpha
\let\bt\beta
\let\dl\delta
\let\Dl\Delta
\let\eps\varepsilon
\let\gm\gamma
\let\Gm\Gamma

\let\la\lambda
\let\La\Lambda

\let\phi\varphi
\let\si\sigma

\let\ups\upsilon

\let\om\omega
\let\Om\Omega

\let\der\partial

\let\Tilde\widetilde

\let\geq\geqslant

\let\leq\leqslant

\let\on\operatorname
\let\bi\bibitem
\let\bs\boldsymbol

\def\C{{\mathbb C}}
\def\Z{{\mathbb Z}}
\def\Q{{\mathbb Q}}

\def\Pb{{\mathbb P}}

\def\F{{\mc F}}

\def\+#1{^{\{#1\}}}

\def\diag{\on{diag}}
\def\End{\on{End}}

\def\Hom{\on{Hom}}

\def\tr{\on{tr}}

\def\sln{\mathfrak{sl}_N}

\def\beq{\begin{equation}}
\def\eeq{\end{equation}}
\def\be{\begin{equation*}}
\def\ee{\end{equation*}}

\nc{\bea}{\begin{eqnarray*}}
\nc{\eea}{\end{eqnarray*}}
\nc{\bean}{\begin{eqnarray}}
\nc{\eean}{\end{eqnarray}}
\nc{\bal}{\begin{align*}}
\nc{\eal}{\end{align*}}
\nc{\baln}{\begin{align}}
\nc{\ealn}{\end{align}}

\nc{\Il}{{\mc I_{\bs\la}}}
\nc{\bla}{{\bs\la}}
\nc{\Fla}{\F_\bla}
\nc{\tfl}{{T^*\Fla}}
\nc{\GL}{{GL_n(\C)}}
\nc{\GLC}{{GL_n(\C)\times\C^*}}

\let\sd s 

\def\ddk_#1{\kk_{#1}\<\>\frac\der{\der\<\>\kk_{#1}}}

\def\bul{\mathbin{\raise.2ex\hbox{$\sssize\bullet$}}}
\def\intt{\mathchoice
{\mathop{\raise.2ex\rlap{$\,\,\ssize\backslash$}{\intop}}\nolimits}
{\mathop{\raise.3ex\rlap{$\,\sssize\backslash$}{\intop}}\nolimits}
{\mathop{\raise.1ex\rlap{$\sssize\>\backslash$}{\intop}}\nolimits}
{\mathop{\rlap{$\sssize\<\>\backslash$}{\intop}}\nolimits}}

\let\kk q 
\let\cc c

\let\Ko K

\def\GZ/{Gelfand-Zetlin}
\def\KZ/{{\slshape KZ\/}}
\def\qKZ/{{\slshape qKZ\/}}
\def\XXX/{{\slshape XXX\/}}

\nc{\slnl}{{\sln (\lambda)}}
\nc{\PCN}{{   (\C[x])^N   }}
\nc{\di}{\on{Diag}}
\nc{\dio}{\on{Diag}_0}
\nc{\Mm}{{\mc M}}
\nc{\Nn}{{\mc N}}

\nc{\A}{{\mc C}}

\nc{\PCr}{{  P  (\C[x])^n   }}

\nc{\Pk}{{(\bs{P}^1)^k}}

\nc{\N}{{\Bbb N}}

\nc{\Ll}{{\mc L}}

\nc{\ord}{{\on{ord}\,}}

\nc{\Sing}{{\on{Sing}\,}}
\nc{\sing}{{\on{Sing}\,}}

\nc{\Hess}{{\on{Hess}}}

\nc{\R}{{\Bbb R}}

\let\on\operatorname
\nc{\Kk}{{\bs K}}
\nc{\Ap}{{A_\Phi(z)}}
\nc{\ap}{{A_\Phi(z)}}

\nc{\sv}{{\sing V}}
\nc{\cd}{{\C^n-\Delta}}
\nc{\UT}{{U^0}}   
\nc{\ep}{\epsilon}

\usepackage[all,cmtip]{xy}
\usepackage{empheq}
\usepackage{bm}

\newlanguage\fakelanguage
\newcommand\cyr{\fontencoding{OT2}\fontfamily{wncyr}\selectfont
   \language\fakelanguage}
\DeclareTextFontCommand{\textcyr}{\cyr}

\numberwithin{equation}{section}

\DeclareMathOperator{\HOM}{\mathscr{H}\text{\kern -3pt {\calligra\large om}}\,}

\makeatletter
\newsavebox{\@brx}
\newcommand{\llangle}[1][]{\savebox{\@brx}{\(\m@th{#1\langle}\)}%
  \mathopen{\copy\@brx\kern-0.5\wd\@brx\usebox{\@brx}}}
\newcommand{\rrangle}[1][]{\savebox{\@brx}{\(\m@th{#1\rangle}\)}%
  \mathclose{\copy\@brx\kern-0.5\wd\@brx\usebox{\@brx}}}
\makeatother

\newcommand{\clr}{\textcolor[rgb]{1.00,0.00,0.00}}

\newcommand{\bsh}{\begin{shaded}}
\newcommand{\esh}{\end{shaded}}

\newcommand{\wPsi}{\widetilde\Psi}
\newcommand{\wV}{\widetilde V}
\newcommand{\wGm}{\widetilde\Gm}
\newcommand{\mrg}{\mathring}

\newcommand{\cfsp}{{\rm Conf}_n(\C)}
\newcommand{\cfspo}{{\rm Conf}_n^O(\C)}
\newcommand{\cfspz}{{\rm Conf}_n^Z(\C)}
\newcommand{\cfspu}{\widehat{{\rm Conf}_n(\C)}}
\newcommand{\Xn}{\mc X_n}
\newcommand{\Xno}{\mc X_n^O}
\newcommand{\Xnz}{\mc X_n^Z}
\newcommand{\Xnu}{\widehat{\mc X_n}}
\newcommand{\An}{\mc A_n}
\newcommand{\Ano}{\mc A_n^O}
\newcommand{\Anz}{\mc A_n^Z}
\newcommand{\Anu}{\widehat{\mc A_n}}
\newcommand{\Ln}{\overline{L}_n}
\newcommand{\LM}[1]{\overline{L}_{#1}}
\newcommand{\DM}[2]{\overline{\mc M}_{#1,#2}}
\newcommand{\sutp}{stable unordered 2-partition}

\begin{document}
\title[RHB inverse problem for flat $F$-manifolds]
{Riemann-Hilbert-Birkhoff inverse problem for semisimple flat $F$-manifolds, and convergence of oriented associativity potentials}
\author[Giordano Cotti ]
{Giordano Cotti$\>^{\circ,\star}$ }
\address{Faculdade de Ci\^encias da Universidade de Lisboa, Grupo de F\'isica Matem\'atica, Campo Grande Edif\'icio C6, 1749-016 Lisboa, Portugal}
\email{gcotti@fc.ul.pt, gcotti@sissa.it}
\subjclass[2020]{Primary: 53D45, 34M50; Secondary: 34M40, 34M55. }
\maketitle
\begin{center}
\textit{ $^\circ\>$Departamento de Matemática\\Faculdade de Ci\^encias da Universidade de Lisboa\\
Campo Grande Edif\'icio C6, 1749-016 Lisboa, Portugal\/}
\vskip1,5mm
\textit{ $^\star\>$Grupo de F\'isica Matem\'atica\\Universidade de Lisboa\\
Campo Grande Edif\'icio C6, 1749-016 Lisboa, Portugal\/}

\end{center}
{\let\thefootnote\relax
\footnotetext{\vskip5pt 
\noindent
$^\circ\>$\textit{ E-mail}:  gcotti@fc.ul.pt, gcotti@sissa.it}}

\vskip6mm
\begin{abstract}
In this paper, we address the problem of classification of quasi-homogeneous formal power series providing solutions of the oriented associativity equations. 
Such a classification is performed by introducing a system of monodromy local moduli on the space of formal germs of homogeneous semisimple flat $F$-manifolds. This system of local moduli is well-defined on the complement of the {\it strictly doubly resonant} locus, namely a locus of formal germs of flat $F$-manifolds manifesting both coalescences of canonical coordinates at the origin, and resonances of their {\it conformal dimensions}. It is shown how the solutions of the oriented associativity equations can be reconstructed from the knowledge of the monodromy local moduli via a Riemann-Hilbert-Birkhoff boundary value problem. Furthermore, standing on results of B.\,Malgrange and C.\,Sabbah, it is proved that any formal homogeneous semisimple flat $F$-manifold, which is not strictly doubly resonant, is actually convergent. Our semisimplicity criterion for convergence is also reformulated in terms of solutions of Losev-Manin commutativity equations, growth estimates of correlators of $F$-cohomological field theories, and solutions of open Witten-Dijkgraaf-Verlinde-Verlinde equations.
\end{abstract}

\tableofcontents

\section{Introduction}

\noindent{\bf Oriented associativity equations. }In this paper, we address both the problem of classification and the convergence issues of formal solutions, in the ring of formal power series with complex coefficients, of the {\it oriented associativity equations} \cite{LM00,LM04,Man05}. These consist of the overdetermined system of non-linear partial differential equations, in $n$ functions $F^1(\bm t),\dots, F^n(\bm t)$  depending on $n$ variables $\bm t=(t^1,\dots, t^n)$, given by
\begin{align*}
&\sum_\mu\frac{\der F^\al}{\der t^\bt\der t^\mu}\frac{\der F^\mu}{\der t^\gm\der t^\eps}=\sum_\mu\frac{\der F^\al}{\der t^\eps\der t^\mu}\frac{\der F^\mu}{\der t^\gm\der t^\bt},\quad &\al,\bt,\gm,\eps=1,\dots,n,\\
&\sum_\mu A^\mu \frac{\der F^\al}{\der t^\mu\der t^\bt}=\dl^\al_\bt,\quad A^\mu\in\C,&\al,\bt=1,\dots, n.
\end{align*}
The oriented associativity equations are a natural generalization of Witten-Dijkgraaf-Verlinde -Verlinde (WDVV) associativity equations \cite{Wit90,DVV91}. Their solutions $(F^1,\dots, F^n)$ reflect the rich geometry of {\it $F$-manifolds with a compatible flat structure}, for short {\it flat $F$-manifolds} \cite{Man05}.
\vskip2mm
\noindent{\bf Flat $F$-manifolds. }In the early 1990s, B.\,Dubrovin introduced Frobenius manifolds as geometrical materialization of solutions of WDVV equations \cite{Dub92,Dub96,Dub98,Dub99}. This notion turns up in many areas of Mathematics: for example Frobenius manifolds play a key role in mirror symmetry, singularity theory, quantum cohomology, integrable systems, and symplectic geometry.
\vskip1,5mm
It was soon realized, however, that weaker (i.e.\,\,with relaxed axioms) variants of the Frobenius structure are of interest {\it per se}. The core notion is that of $F$-manifolds, introduced by C.\,Hertling and Yu.I.\,Manin in \cite{HM99}. Such a notion not only strictly includes the Frobenius structures, but it also greatly broadens the scope of the examples and applications. Examples of $F$-manifolds, indeed, arise not only in singularity theory \cite{Her02}, but also in quantum $K$-theory \cite{Lee04}, differential-graded deformation theory \cite{Mer04,Mer06}, and even information geometry \cite{CM20}.
\vskip1,5mm
Flat $F$-manifolds -- introduced by Yu.I.\,Manin \cite{Man05}-- are an intermediate notion, weaker than Frobenius, but stronger than $F$-manifolds:
\[\text{Frobenius manifolds}\quad\subset\quad\text{flat $F$-manifolds}\quad\subset\quad\text{$F$-manifolds}.
\] Flat $F$-manifolds are equipped with the minimum amount of structures to share some of the deeper properties of Frobenius manifolds, including Dubrovin's deformed connection, Dubrovin's almost duality, and also operadic descriptions. See \cite{Man04,Man05,LPR11,AL13,AL17}. These structures are also studied in \cite{Get04}, where they are called {\it Dubrovin manifolds}.
\vskip1,5mm
A flat $F$-manifold (in the analytic category) is a complex manifold $M$ whose tangent spaces are equipped with an associative, commutative and unital algebra structure --analytically depending in the point-- whose product $\circ$ is {\it compatible} with a given flat connection $\nabla$. This means that each element of the pencil $(\nabla^z)_{z\in\C}$, defined by $\nabla^z_XY=\nabla_XY+zX\circ Y$, is required to be flat and torsionless.

The compatibility of $(\circ,\nabla)$ implies a {\it potentiality} condition. In $\nabla$-flat local coordinates $\bm t=(t^1,\dots, t^n)$ on $M$, with $n=\dim_\C M$, the product $\circ$ descends from a vector potential: there exists $\bm F=(F^1,\dots, F^n)$ such that
\beq
\label{printro}
\frac{\der}{\der t^\bt}\circ \frac{\der}{\der t^\gm}=\sum_\al \frac{\der^2 F^\al}{\der t^\bt\der t^\gm}\frac{\der}{\der t^\al},\quad \bt,\gm=1,\dots, n.
\eeq
The associativity of $\circ$ is equivalent to the oriented associativity equations for $\bm F$. Vice-versa, starting from a solution $\bm F$ of the oriented associativity equations, we can define a flat $F$-structure via equation \eqref{printro}. If the starting solution $\bm F$ is a tuple of formal power series in $k[\![\bm t]\!]$ (with $k$ a $\Q$-algebra), the resulting flat $F$-structure is said to be {\it formal} over $k$. It can be seen as a flat $F$-structure on the formal spectrum ${\rm Spf}\,k[\![\bm t]\!]$.
\vskip2mm
\noindent{\bf Homogeneity, semisimplicity, and double resonance. }In this paper we consider only {\it quasi-homogenous} solutions $\bm F$ of the oriented associativity equations, i.e.\,\,satisfying a further condition of the form
\[\sum_\al[(1-q_\al)t^\al+r^\al]\frac{\der F^\bt}{\der t^\al}=(2-q_\bt)F^\bt(\bm t)+\text{linear terms in }\bm t\text{ and constant terms},
\]for suitable complex numbers $q_\al,r^\al\in\C$. The resulting flat $F$-manifold is said to be {\it homogeneous}, or {\it conformal}. The vector field $E=\sum_\al[(1-q_\al)t^\al+r^\al]\frac{\der }{\der t^\al}$ is then an {\it Euler} vector field, i.e. it satisfies the conditions $\nabla\nabla E=0$ and $\frak L_E(\circ)=\circ$. We say that $p\in M$ is {\it tame} if the spectrum of the operator $(E\circ)_p\in\End(T_pM)$ is simple, otherwise we say that $p$ is {\it coalescing}. 
\vskip1,5mm
An analytic flat $F$-manifold is said to be {\it semisimple} if there exists an open dense subset of points $p$ whose corresponding algebra $(T_pM,\circ_p)$ is without nilpotent elements. This is equivalent to the existence of idempotent vectors $\pi_1,\dots, \pi_n\in T_pM$: $\pi_i\circ \pi_j=\pi_i\dl_{ij}$ for $i,j=1,\dots,n$.
If a manifold is both homogenous and semisimple, the eigenvalues $u^1,\dots, u^n$ of the tensor $(E\circ)$ can be chosen as local holomorphic coordinates in a neighborhood of any semisimple point $p\in M$. Tame points are necessarily semisimple, whereas coalescing points may or may not be semisimple. In the (tame/coalescing) semisimple case, the vector fields $\frac{\der}{\der u^1},\dots, \frac{\der}{\der u^n}$ are the idempotent vectors.
\vskip1,5mm
With each homogeneous semisimple (analytic/formal) flat $F$-manifold we can associate a tuple $(\dl_1,\dots, \dl_n)$ of numerical invariants called \emph{conformal dimensions}. Fix a semisimple point $p\in M$, and introduce the operator
\[\mu^0_p\in\End(T_pM),\quad \mu^0_p\left(\frac{\der}{\der t^\al}\right)=q_\al\frac{\der}{\der t^\al},\quad \al=1,\dots,n.
\]The conformal dimensions can be defined as the numbers $\dl_1,\dots, \dl_n\in\C$ satisfying
\[\mu^0_p(\pi_i)\circ \pi_i=\dl_i\pi_i,\quad i=1,\dots,n.
\]They are defined up to ordering, and they actually do not depend on the chosen semisimple point $p\in M$. The flat $F$-manifold is be said to be {\it conformally resonant} if $\dl_i-\dl_j\in\Z\setminus\{0\}$ for some $i,j$. The conformal dimensions of a Frobenius manifold are all equal, with common value $\frac{d}{2}$ (the number $d\in\C$ is the {\it charge} of the Frobenius structure). In particular, Frobenius manifolds are never conformally resonant. 
\vskip1,5mm
 In the formal case, all the conditions on points introduced above (tameness, coalescence, semisimplicity) are intended to be referred to the origin $\bm t=0$, the only geometric point of the formal spectrum ${\rm Spf}\,{ k[\![\bm t]\!]}$.
 
 A (germ of) pointed flat $F$-manifold $(M,p)$ is be said to be {\it doubly resonant} if $M$ is conformally resonant, and $p$ is coalescing. We say that $(M,p)$ is {\it strictly doubly resonant} if we have\footnote{Any ordering of the eigenvalues $u^1,\dots, u^n$ induces an ordering of the idempotent vectors $\pi_i\equiv \der_{u^i}$, and consequently of the conformal dimensions $\dl_i$.} 
 \[u^i(p)=u^j(p),\qquad \dl_i-\dl_j\in\Z\setminus\{0\},\qquad \text{for some }i,j\in\{1,\dots,n\},\,i\neq j.
 \]
\vskip2mm
\noindent{\bf Results. }
One of the main aspects of Dubrovin's analytic theory of Frobenius manifolds is their isomonodromic approach. Under the quasi-homogeneity assumption of the WDVV potential, the semisimple part of a Frobenius manifold can be locally identified with the space of isomonodromic deformation parameters of ordinary differential equations on $\Pb^1$ with rational coefficients, see \cite{Dub98,Dub99}.
\vskip1,5mm
In this paper, we extend to the case of homogenous semisimple flat $F$-manifolds both Dubrovin's analytical theory as well as its refinement developed in \cite{CG17,CDG20}. The key ingredient is a family $(\widehat\nabla^\la)_{\la\in\C}$ of flat extended deformed connections on $\pi^*TM$, with $\pi\colon M\times \C^*\to\C^*$, whose restrictions to $M\times \{z\}$ equal $\nabla^z=\nabla+z(-\circ-)$.  These families of flat connections $\widehat\nabla^\la$ on  homogeneous flat $F$-manifolds were first introduced by Yu.I.\,Manin \cite{Man05}.
\begin{rem}
In \cite[Chapter VII, \S 1]{Sab07} the notion of {\it Saito structure without metric} is introduced. These structures are defined in terms of several data on a complex manifold $M$. Among them there is a flat meromorphic connection $\widehat\nabla$ on the bundle $\pi^*TM$ on $M\times \mathbb P^1$. It turns out that the notion of Saito structure without metric is equivalent to the notion of homogeneous flat $F$-manifold, and that $\widehat\nabla$ is one of the connections $\widehat\nabla^\la$ above. See also \cite{KMS20}, in which it is shown that the space of isomonodromic deformation parameters of extended Okubo systems can be equipped with Saito structures without metrics.
\end{rem}
For any germ $(M,p)$, semisimple and not strictly doubly resonant, %
we introduce a tuple of numerical data, the {\it monodromy data} of the flat $F$-manifold. These data split into two pieces: a pair $(\mu^\la,R)$ of matrices called ``monodromy data at $z=0$'', and a 4-tuple $(S_1,S_2,\La,C)$ of matrices called ``monodromy data at $z=\infty$''. Precise definitions are given in Section \ref{sec4}. In the case of Frobenius manifolds, all these data are subjected to several constraints: the final amount of data coincides with the 4-tuple $(\mu,R,S,C)$ of monodromy data introduced by Dubrovin in \cite{Dub98,Dub99}.
If $M$ is analytic, the monodromy data define {\it local invariants} of $M$: if $p_1,p_2\in M$ are sufficiently close, the data of the germs $(M,p_1),(M,p_2)$ are equal. 
\vskip1,5mm
\begin{thm}[Cf. Theorems \ref{enF}, \ref{mthm1}]
Any homogeneous semisimple (analytic/formal) pointed germ of flat $F$-manifold, which is not strictly doubly resonant, 
is uniquely determined by its monodromy data. In particular, the vector potential $\bm F$ can be explicitly reconstructed  from the monodromy data via a Riemann-Hilbert-Birkhoff boundary value problem.
\end{thm}
We show that the totality of local isomorphism classes of germs of $n$-dimensional flat $F$-manifolds can be parametrized by points of a ``stratified'' space, whose generic stratum has dimension $n^2-n$. 
The monodromy data provide a system of local coordinates. The Frobenius structures correspond to a locus of generic dimension $\frac{1}{2}(n^2-n)$. See Theorem \ref{numbermoduli}.
\vskip1,5mm
The re-construction procedure of the flat $F$-structure is based on a crucial property of a joint system of {``generalized'' Darboux-Egoroff equations} \cite{Lor14}: solutions $\Gm(\bm u)$ of this system of non-linear partial differential equations are uniquely determined by its initial value $\Gm_o$ at one point $\bm u_o\in\C^n$ (with possibly $u_o^i=u_o^j$ for $i\neq j$), provided there is no conformal resonance. %
See Lemma \ref{miracle}. 
\vskip1,5mm
We underline that both the initial values $(\bm u_o,\Gm_o)$ and the monodromy data provide a system of local coordinates on the space of germs of flat $F$-manifolds. The reconstruction procedure of $\bm F$ in terms of the initial values, however, is generally impossible, the dependance being typically transcendental (e.g.\,\,for $n=3$, the general Darboux-Egoroff system reduces to the full-parameters family of Painlev\'e equations PVI, see \cite{Lor14}). This makes the monodromy data ``preferable'' as a system of coordinates for the classification of flat $F$-structures.
\vskip1,5mm
There is a further advantage in choosing the monodromy data as a system of local moduli. Indeed, they make possible the study of convergence issues.

\begin{thm}[Cf. Theorem \ref{thconv}]\label{th2intro}
Let $\bm F\in \C[\![\bm t]\!]^{\times n}$ be a quasi-homogeneous solution of the oriented associativity equations. If $\bm F$ defines a semisimple formal flat $F$-manifold, which is not strictly doubly resonant, %
then $\bm F$ is a tuple of convergent functions.
\end{thm}
For the proof, we invoke results of B.\,Malgrange \cite{Mal83a,Mal83b,Mal86} and C.\,Sabbah \cite{Sab18} on the solvability of families of Riemann-Hilbert-Birkhoff problems. More precisely, we use their equivalent formulations given in \cite{Cot20b}. Notice that Theorem \ref{th2intro} generalizes \cite[Theorem 5.1]{Cot20b}. 
\vskip2mm
\noindent{\bf Cohomological field theories. }%
Frobenius structures are intrinsically correlated to the notion of cohomological field theories. The geometry of Frobenius structures reflects properties of the cohomology rings $H^\bullet(\DM{0}{n},\C)$, with $n\geq 3$, and they can indeed be defined in terms of $H^\bullet(\DM{0}{n},\C)$-valued poly-linear maps. See \cite{KM94,man,Pan18}.
\vskip1,5mm
At the level of flat $F$-manifolds, such a construction has been generalized in two different (a posteriori equivalent) ways. The first one is due to A.\,Losev and Yu.I.\,Manin \cite{LM00,LM04}, the second one to A.\,Buryak and P.\,Rossi \cite{BR18}. See also \cite{ABLR20,ABLR20b,ABLR21}.
\vskip1,5mm
In \cite{LM00}, a new compactification $\LM{n}$ of $\mc M_{0,n}$ is introduced. The boundary strata represent isomorphism classes of stable $n$-pointed {\it chains} of projective lines, in which the marked points do not play a symmetric role. In \cite{LM04}, a notion of {\it genus 0 extended modular operad}, and $\mc L$-{\it algebras} over it are studied. Moreover, it is shown that the differential equations satisfied by generating functions of correlators of $\mc L$-algebras lead to two differential geometric pictures: the first one is the study of pencils of flat connections, based on the {\it commutativity equations}, the second one is the study of flat $F$-manifolds, based on the oriented associativity equations.
These two pictures are actually locally equivalent, if a further amount of data -- a {\it primitive element} -- is given. See also the constructions in \cite{Los97,Los98}.
\vskip1,5mm
Notice that the compactifications $\LM{n}$ and $\DM{0}{n}$, and their higher genus analogs, both arise in the more general construction of \cite{Has03} as compactified moduli spaces of {\it weighted} pointed stable curves, for two different choices of the weights. See also \cite{Man04}.
\vskip1,5mm
In \cite{BR18} A.\,Buryak and P.\,Rossi introduced the notion of $F$-{\it cohomological field theories} ($F$-CohFT) as a generalization of both cohomological field theories \cite{KM94,man}, and partial cohomological field theories \cite{LRZ15}. An $F$-CohFT is defined by the datum of a family of $H^\bullet(\DM{g}{n},\C)$-valued poly-linear maps on a tensor product $V^*\otimes V^{\otimes n}$ ($V$ is an arbitrary $\C$-vector space), which satisfy some natural $\frak S_n$-equivariance and gluing properties. Given an $F$-CohFT, its {\it genus zero sector} (or {\it tree level}) defines a formal flat $F$-structure on $V$. 
\vskip1,5mm
In Section \ref{sec7}, we review all these cohomological field theoretical approaches to flat $F$-manifolds, their equivalences, and we rephrase our semisimplicity criterion of convergence in terms of solutions of Losev-Manin commutativity equations, and growth estimates of correlators of $F$-CohFT's. Furthermore, in Appendix \ref{app} we prove\footnote{The author is not aware of a proof of this fact in literature.} that any formal flat $F$-manifold over $\C$ descend from a unique $F$-CohFT in the sense of P.\,Rossi and A.\,Buryak. This is in complete analogy with the Frobenius manifolds case, see \cite{man}.
\vskip2mm
\noindent{\bf Structure of the paper. }In Section \ref{sec2} we present some preliminary material and basic properties of flat $F$-manifolds, in both formal and analytic categories. We recall definitions of homogeneity, semisimplicity of flat $F$-manifolds. We also introduce the notion of local isomorphisms, pointed germs, and irreducibility of flat $F$-manifolds. 
\vskip1,5mm
In Section \ref{sec3}, we firstly describe how to reconstruct oriented associativity potentials from deformed coordinates on a flat $F$-manifold. We then introduce a family of flat extended deformed connections $\widehat\nabla^\la$, and we develop the analytical theory of its flat (co)sections. We study the differential system of $\widehat\nabla^\la$-flatness in three different frames: the flat frame, the idempotent frame, and a normalized idempotent frame (the normalizing factors are the so-called {\it Lam\'e coefficients}). We also introduce Darboux-Tsarev and Darboux-Egoroff systems of equations.
\vskip1,5mm
In Section \ref{sec4}, after introducing the notion of spectrum of a flat $F$-manifolds, we define the monodromy data of a homogeneous semisimple flat $F$-manifold. We also describe their mutual constraints.
In Section \ref{sec5}, we then clarify the dependence of the monodromy data on all the choices of normalizations involved in their definition. Different choices of normalizations affect the numerical values of the data via the action of suitable groups. We show that the analytic continuation of the flat $F$-structure is described by a braid group action on the tuple of monodromy data.
\vskip1,5mm
Section \ref{RHBsec} contains the main results of the paper. After introducing the notion of admissible data and the related Riemann-Hilbert-Birkhoff (RHB) boundary value problem, we recall some results of B.\,Malgrange and C.\,Sabbah as formulated in \cite{Cot20b}. We show that germs of flat $F$-structures can be constructed starting from solutions of RHB problems. Moreover, we show that any germ, which is not strictly doubly resonant, is of such a form: it can be reconstructed from its monodromy data. It is also proved that any formal germ of homogenous semisimple flat $F$-manifold (over $\C$), which is not strictly doubly resonant, %
is actually convergent.
\vskip1,5mm
In Section \ref{sec7}, we review equivalent approaches for defining flat $F$-manifolds. We recall the notions of Losev-Manin commutativity equations, Losev-Manin cohomological field theories (LM-CohFT), and $F$-cohomological field theories ($F$-CohFT) in the sense of A.\,Buryak and P.\,Rossi. We discuss the equivalence of these notions. Furthermore, we also discuss relations with the {\it open} WDVV (OWDVV) equations. Our semisimplicity criterion is reformulated in terms of growth estimates of correlators of LM-CohFT and $F$-CohFT, and convergence of solutions of OWDVV equations.
\vskip1,5mm
In Appendix \ref{theuelrid}, we provide a proof for the following characterization of irreducibility of flat $F$-structures in terms of Euler vector fields: a flat $F$-manifold is irreducible if and only if any two arbitrary Euler vector fields differ by a scalar multiple of the unit vector field.
\vskip1,5mm
In Appendix \ref{app}, we prove that any formal flat $F$-manifold descends from a unique tree-level $F$-CohFT.
\vskip2mm
\noindent{\bf Acknowledgements. }The author is thankful to %
D.\,Guzzetti, C.\,Hertling, A.R.\,Its, P.\,Loren\-zoni, Yu.I. Man\-in, D.\,Masoero, %
A.T.\,Ricolfi, P.\,Rossi, V.\,Roubtsov, A.\,Var\-chenko, C. Sabbah, M. Smirnov, D.\,Yang for several valuable discussions. The author is thankful to the Hausdorff Research Institute for Mathematics (HIM) in Bonn, Germany, where this project was started, for providing excellent working conditions during the JTP ``New Trends in Representation Theory''. This research was supported by HIM (Bonn, Germany), and by the FCT Projects UIDB/00208/2020, PTDC/MAT-PUR/ 30234/2017, 2021.01521.CEECIND, and 2022.03702.PTDC (GENIDE). The author is member of the COST Action CA21109 CaLISTA.
\section{Flat $F$-manifolds}\label{sec2}
\subsection{Analytic flat $F$-manifolds } Let $M$ be a complex analytic manifold with dimension $\dim_\C M=n$. Denote by $TM,\,T^*M$ the holomorphic tangent and cotangent bundles, and by $\mathscr T_M,\,\Om^1_M$ their sheaves of sections. If $E\to M$ is a holomorphic bundle with sheaf of sections $\mathscr E$, we denote by $\Gm(E)=\Gm(M,\mathscr E)$ the space of global sections. By usual abuse of notations, we will write $X\in \mathscr E$ for $X\in\Gm(U,\mathscr E)$ for some (or arbitrary) open set $U\subseteq M$.
\vskip1mm
Let $(M,\nabla,c,e)$ be the datum of
\begin{enumerate}
\item a connection $\nabla\colon \mathscr T_M\to \Om^1_M\otimes \mathscr T_M$ on $TM$;
\item a section $c\in\Gm\left(TM\otimes \bigodot^2 T^*M\right)$;
\item a vector field $e\in\Gm(TM)$ such that 
\begin{enumerate}
\item $c(-,e,-)=c(-,-,e)\in \Gm({\rm End}\, TM)$ is the identity morphism,
\item $\nabla e=0$.
\end{enumerate} 
\end{enumerate}
Denote by $X\circ Y:=c(-,X,Y)$ the commutative product defined by the tensor $c$, and introduce the one-parameter family of connections $\left(\nabla^z\right)_{z\in\C}$ defined by $\nabla^z_XY:=\nabla_XY+z X\circ Y$ for $X,Y\in\mathscr T_M$.

\begin{defn}
We say that $(M,\nabla,c,e)$ is a \emph{flat $F$-manifold} if the connection $\nabla^z$ is flat and torsionless for any $z\in\C$.
\end{defn}

Let $\bm t=(t^1,\dots, t^n)$ be a system of $\nabla$-flat local coordinates on $M$. Set $\der_\al:=\frac{\der}{\der t^\al}$, for $\al=1,\dots,n$ and define $c^\al_{\bt\gm}:=c(dt^\al,\der_\bt,\der_\gm)$. The flatness and torsionless of $\nabla^z$ is equivalent to the associativity of $\circ$, the symmetry of $\der_\al c_{\bt\gm}^\dl$ in $(\al,\bt,\gm)$, and the flatness of $\nabla$. Hence, there locally exist analytic functions $\bm F=(F^1,\dots,F^n)\in\mc O_M^n$ such that
\[c^\al_{\bt\gm}=\frac{\der^2 F^\al}{\der t^\bt\der t^\gm},\quad \al,\bt,\gm=1,\dots,n.
\]
In what follows the Einstein summation rule is used over repeated Greek indices. Let $A^\al\in\C$ be constants such that $e=A^\al\der_\al$. From the associativity of $\circ$ and the properties of $e$, we have
\begin{align}
\label{wdvv1}
A^\mu\frac{\der^2 F^\al}{\der t^\mu\der t^\bt}&=\dl^\al_\bt,\quad&\al,\bt=1,\dots,n,\\
\label{wdvv2}
\frac{\der^2 F^\al}{\der t^\mu\der t^\bt}\frac{\der^2 F^\mu}{\der t^\gm\der t^\dl}&=\frac{\der^2 F^\al}{\der t^\mu\der t^\gm}\frac{\der^2 F^\mu}{\der t^\bt\der t^\dl},\quad&\al,\bt,\gm,\dl=1,\dots,n.
\end{align}
Equations \eqref{wdvv1}, \eqref{wdvv2} are called \emph{oriented associativity equations}, and $\bm F$ is the \emph{oriented associativity potential} of the flat $F$-structure.
\vskip3mm
A flat $F$-manifold is said to be \emph{homogeneous} if there it is equipped with an {\it Euler vector field}, i.e. a vector field $E\in\Gm(TM)$ such that
\[\nabla\nabla E=0,%
\quad \frak L_Ec=c.
\]
\begin{lem}\label{commeE}
We have $[e,E]=e$.
\end{lem}
\proof
The condition $\frak L_Ec=c$ is equivalent to $[E,Y\circ Z]-[E,Y]\circ Z-[E,Z]\circ Y=Y\circ Z$, for $Y,Z\in\mathscr T_M$. If $Y=Z=e$, the identity follows.
\endproof
Homogeneous flat $F$-manifold are also called \emph{Saito structures without metric}, \cite[Ch.\,VII]{Sab07}. We assume that the $(1,1)$-tensor $\nabla E$ is diagonalizable, and in diagonal form in the $\bm t$-coordinates:
\beq
\label{evf}
E=\sum_{\al=1}^n((1-q_\al)t^\al+r^\al)\frac{\der}{\der t^\al},\quad q_\al,r^\al\in\C.
\eeq
The condition $\frak L_Ec=c$ is thus equivalent to
\beq
\label{hEF}
E^\mu\frac{\der F^\al}{\der t^\mu}=(2-q_\al)F^\al+A^\al_\bt t^\bt+B^\al,\quad A^\al_\bt,B^\al\in\C.
\eeq

\begin{defn}
A flat $F$-manifold $(M,\nabla,c,e)$ is a {\it Frobenius manifold} if there exist a symmetric non-degenerate $\mc O_M$-bilinear form $\eta\in \Gm(\bigodot^2T^*M)$, called {\it metric}, such that
\beq\label{mfrob}
\nabla\eta=0,\quad \text{and}\quad \eta(X\circ Y,Z)=\eta(X,Y\circ Z),\quad X,Y,Z\in\mathscr T_M.
\eeq
In such a case, $\nabla$ is the Levi-Civita connection of $\eta$. A vector field $E\in\Gm(TM)$ is {\it Euler} if it satisfies the conditions
\beq\label{kill}
\frak L_Ec=c,\quad \frak L_E\eta=(2-d)\eta,
\eeq where the number $d\in\C$ is the {\it conformal dimension} (or {\it charge}) of the Frobenius manifold.
\end{defn}
\begin{rem}
An Euler vector field $E$ for a Frobenius manifold is automatically an Euler vector field for the underlying flat $F$-manifold. The condition $\nabla\nabla E=0$ is indeed implied by the conformal Killing condition \eqref{kill}, and the flatness of $\nabla$.
\end{rem}

\begin{rem}
In the case a flat $F$-manifold is actually Frobenius, the oriented associativity potentials $\bm F=(F^1,\dots, F^n)$, solutions of \eqref{wdvv1} and \eqref{wdvv2}, can be shown to locally descend from a single WDVV potential $F(\bm t)$, i.e. a solution of the system of equations
\begin{align*}
&A^\mu\frac{\der^3F}{\der t^\mu\der t^\al\der t^\bt}=\eta_{\al\bt}={\rm const.,}\quad \eta=(\eta_{\al\bt})_{\al,\bt},\quad \eta^{-1}=(\eta^{\al\bt})_{\al,\bt}&\al,\bt=1,\dots, n,\\
&\frac{\der^3F}{\der t^\al\der t^\bt\der t^\mu}\eta^{\mu\nu}\frac{\der^3F}{\der t^\nu\der t^\gm\der t^\dl}=\frac{\der^3F}{\der t^\dl\der t^\bt\der t^\mu}\eta^{\mu\nu}\frac{\der^3F}{\der t^\nu\der t^\gm\der t^\al},&\quad \al,\bt,\gm,\dl=1,\dots,n.
\end{align*}
The potentials $F^\al$'s are the components of the $\eta$-gradient of $F$, that is $F^\al(\bm t)=\eta^{\al\bt}\der_\bt F(\bm t)$.
\end{rem}
\subsection{Formal flat $F$-manifolds}Let 
\begin{itemize}
\item $k$ be a commutative $\Q$-algebra, 
\item $H$ be a free $k$-module of finite rank, 
\item $K:=k[\![H^*]\!]$ be the completed symmetric algebra of $H^*:=\Hom_k(H,k)$.
\end{itemize}
Fix a basis $(\Dl_1,\dots, \Dl_n)$ of $H$, and denote by $\bm t=(t^1,\dots, t^n)$ the dual coordinates. The algebra $K$ is then identified with the algebra of formal power series $k[\![\bm t]\!]$. Denote by ${\rm Der}_k(K)$ the $K$-module of $k$-linear derivations of $K$. Put $\der_\al=\frac{\der}{\der t^\al}\colon K\to K$. The module ${\rm Der}_k(K)$ is a free $K$-module with basis $(\der_1,\dots,\der_n)$. %
\vskip 2mm
Elements of $H_K:=K\otimes_k H$ will be identified with derivations on $K$, by $\Dl_\al\mapsto \der_\al$. 
\vskip2mm
\begin{defn}
A \emph{formal flat $F$-manifold structure} on $H$ is given by an $n$-tuple $\bm \Phi=(\Phi_1,\dots,\Phi_n)\in K^n$, satisfying the oriented associativity equations \eqref{wdvv1}, \eqref{wdvv2}, where $A^\mu\in k$.
\end{defn}

Define the $K$-linear multiplication $\circ$ on $H_K$ by
\[\Dl_\al\circ \Dl_\bt:=c_{\al\bt}^{\gm}\Dl_\gm,\quad c_{\al\bt}^{\gm}:=\frac{\der^2\Phi^\gm}{\der t^\al \der t^\bt}\quad \al,\bt=1,\dots, n.
\]
The oriented associativity equations imply that $\circ$ is associative, and that $e:=A^\mu\Dl_\mu$ is the unit of the algebra $(H_K,\circ)$. A vector $E\in H_K$ is an {\it Euler} vector if it is if the form \eqref{evf}, and the pair $(E,\bm\Phi)$ satisfies equations \eqref{hEF}.
\vskip3mm
Let ${\rm Diff}_1(H_K,H_K)$ denote the set of $\mathscr D\in \Hom_k(H_K,H_K)$ such that
\[ab\mathscr D(p)-b\mathscr D(ap)-a\mathscr D(bp)+\mathscr D(abp)=0,\quad a,b\in K,\quad p\in H_K.
\]Both ${\rm Der}_k(K)$ and ${\rm Diff}_1(H_K,H_K)$ are naturally equipped with a $K$-module structure. A (formal) connection on $H_K$ is defined by a $K$-linear morphism $\nabla\colon {\rm Der}_k(K)\to {\rm Diff}_1\left(H_K,H_K\right)$, $u\mapsto\nabla_u$ satisfying the Leibniz rule
\[\nabla_u(ap)=u(a)p+a\nabla_up,\quad a\in K,\quad p\in H_K.
\]The torsion and curvature of $\nabla$ are the $K$-bilinear morphisms $T,R\colon {\rm Der}_k(K)\times {\rm Der}_k(K)\to \Hom_K(H_K,H_K)$ defined by
\bea
&T(u,v):=\nabla_uv-\nabla_vu-[u,v],\quad &u,v\in{\rm Der}_k(K)\cong H_K,\\
&R(u,v):=[\nabla_u,\nabla_v]-\nabla_{[u,v]},\quad &u,v\in{\rm Der}_k(K).
\eea
\vskip3mm
We can thus introduce the one-parameter family $(\nabla^z)_{z\in k}$ of (formal) connection given by $\nabla^z_{\der_\al}\der_\bt:=z\der_\al\circ\der_\bt$. The formal connection $\nabla^z$ is flat and torsionless for any $z\in k$.

\begin{rem}\label{resc1}
If $\bm F:=(F^1,\dots, F^n)$ is a (formal/analytic) solution of the oriented associativity equations \eqref{wdvv2}, then also $\widetilde{\bm F}:=(\la_1 F^1,\dots,\la_n F^n)$ is a solution for any $(\la_1,\dots,\la_n)\in(\C^*)^n$. If the original flat $F$-manifold has a unit $e=A^\mu\der_\mu$, then the rescaled flat $F$-manifold structure has unit $e'=\frac{1}{\la_\mu}A^\mu\der_\mu$.
\end{rem}

\begin{rem}\label{resc2}
If $E$ is an Euler vector field for a given flat $F$-manifold structure, then also $E-\la e$ is an Euler vector field, for $\la\in\C$. Under two further assumptions of {\it semisimplicity} and {\it irreducibility}, one can prove that any Euler vector field is of this form. See Theorem \ref{eulerunique}.  
\end{rem}

\subsection{Local isomorphisms and pointed germs} 
Let $(M_{i},\nabla_{i},c_{i},e_{i})$, with $i=1,2$, be two analytic flat $F$-manifolds. A biholomorphism $\phi\colon M_{1}\to M_{2}$ is an {\it isomorphism} of flat $F$-manifolds if 
\begin{enumerate}
\item $d\phi(\ker \nabla_{1})= \ker (\phi^*\nabla_{2})$, where $d\phi\colon \mathscr T_{M_{1}}\to \phi^*\mathscr T_{M_{2}}$ is the differential of $\phi$, and $\phi^*\nabla_{2}$ is the pulled-back connection on $\phi^*TM_{2}$,
\item for each $p\in M_{1}$, the map $d\phi_p\colon T_pM_{1}\to T_{\phi(p)}M_{2}$ is an isomorphism of unital algebras.
\end{enumerate}
\vskip1mm
\begin{lem}\label{lemiso}
Let $M_1$ and $M_2$ be two isomorphic flat $F$-manifolds. Given two systems of local flat coordinates, $\bm t$ on $M_1$ and $\tilde{\bm t}$ on $M_2$, the corresponding local potentials $\bm F_1(\bm t)$ and $\bm F_2(\tilde{\bm t})$ are related by
\bea 
&F_2^\al(\tilde{\bm t})=G^\al_\la F^\la_1(\bm t)+\text{linear terms in }\bm t,\quad &\al=1,\dots,n,\\
 &\tilde{\bm t}=\phi(\bm t)=\bm G\bm t+\bm c,&
\eea
where %
$\bm c\in\C^n$ is a constant vector, and the matrix $\bm G\in GL(n,\C)$ satisfies the relation $\bm GA_1=A_2$, the column vectors $A_i=(A_i^1,\dots,A_i^n)^T\in\C^n$, for $i=1,2$, being such that 
\[e_1=\sum_{\al=1}^nA_1^\al\frac{\der}{\der t^\al},\quad e_2=\sum_{\al=1}^nA_2^\al\frac{\der}{\der \tilde t^\al}.\]
\end{lem}
\proof If $\tilde{\bm t}=\phi(\bm t)$ defines the isomorphism in the chosen local flat coordinates, then the Jacobian matrix $\left(\frac{\der \phi^\al}{\der t^\bt}\right)_{\al,\bt=1}^n$ must be constant. Hence $\tilde{\bm t}=\bm G\bm t+\bm c$, for some $\bm G\in GL(n,\C)$ and $\bm c\in\C^n$. Moreover, by definition of isomorphism, we have
\[\frac{\der \phi^\gm(\bm t)}{\der t^\la}\frac{\der^2F_1^{\la}(\bm t)}{\der t^\al\der t^\bt}=\left.\frac{\der^2 F^\gm_2(\tilde{\bm t})}{\der \tilde{t}^\mu\der \tilde t^\nu}\right|_{\tilde{\bm t}=\phi(\bm t)}\frac{\der \phi^\mu(\bm t)}{\der t^\al}\frac{\der \phi^\nu(\bm t)}{\der t^\bt}=\frac{\der^2(F_2^\gm\circ\phi)(\bm t)}{\der t^\al\der t^\bt}.
\]The claim follows, by double integration.
\endproof
\begin{rem}\label{remgroupG}
The group of transformations $\bm t\mapsto \tilde{\bm t}=\bm G\bm t+\bm c$ of Lemma \ref{lemiso} has dimension $n^2$. Notice that here we are not requiring $\nabla E$ to be diagonal in both systems of coordinates. That further constraint decreases the dimension of the group of admissible transformations. If we want both $\frac{\der}{\der t^\al}$ and $\frac{\der}{\der \tilde t^\al}$, with $\al=1,\dots,n$, to be eigenvectors of $\nabla E$ we need to impose a constraint on $\bm G$: namely, the eigenspaces of $\nabla E$ must be invariant under the linear transformation $\frac{\der}{\der t^\al}\mapsto \sum_{\bt}(\bm G^{-1})^\bt_{\al}\frac{\der}{\der t^\bt}$. In particular, if $\nabla E$ has simple spectrum (i.e.\,\,the parameters $q_1,\dots,q_n$ are pairwise distinct), the matrix $\bm G$ should be diagonal. In total, we would have a group of transformations of dimension $2n$. If the spectrum of $\nabla E$ is not simple, then the dimension of possible transformations increases.
\end{rem}
\vskip1mm
A {\it pointed flat $F$-manifold} is a pair $(M,p)$, where $M$ is a flat $F$-manifold, and $p\in M$ is a fixed \emph{base point}.  Isomorphisms between pointed flat $F$-manifold will always be assumed to be base point preserving. Given $(M,p)$ we will always consider flat coordinates $\bm t$ vanishing at $p$.
\vskip1mm
Two pointed flat $F$-manifolds $(M_1,p_1)$, $(M_2,p_2)$ are {\it locally isomorphic} if there exist open neighborhoods $\Om_1\subseteq M_1$ of $p_1$, and $\Om_2\subseteq M_2$ of $p_2$ respectively, with isomorphic induced flat $F$-structures, i.e. $(\Om_1,p_1)\cong (\Om_2,p_2)$. 
\vskip1mm
A {\it pointed germ} is a (local isomorphism) equivalence class of pointed flat $F$-manifolds.
\vskip2mm
Any analytic pointed flat $F$-manifold $(M,p)$ induces a formal flat $F$-manifold $(H,\bm \Phi)$ over $k=\C$. Choose flat coordinates $\bm t$ vanishing at $p$, and set $H:=T_pM$. Let $\mc O_{M,p}$ be the local ring of germs at $p$, and $\frak m$ be its maximal ideal. The formal potential $\Phi^\al$ is given by the image of $F^\al$ in the completion $\widehat{\mc O_{M,p}}:=\varprojlim \left(\mc O_{M,p}/{\frak m}^\ell\right)$ of the local ring $\mc O_{M,p}$: this means that $\Phi^\al$ is defined by the Taylor series expansion of $F^\al$ at $p$ in coordinates $\bm t$. Moreover, the formal flat $F$-structure $(H,\bm\Phi)$ is also equipped with a flat unit $e|_p$. If $M$ has Euler vector field $E$, then $(H,\bm\Phi)$ has Euler vector field $E|_p$. We will say that the formal flat $F$-structure constructed in this way, starting from an analytic one, is \emph{convergent}.
\vskip1mm
Vice-versa, let us assume that $(H,\bm\Phi)$ is a formal flat $F$-structure over $k=\C$ (with Euler vector field). If the common domain of convergence $\Om\subseteq H$ of the power series $\Phi^{\al}\in k[\![\bm t]\!]$ is non-empty, it is easily seen that $\Om$ is equipped with an analytic flat $F$-structure (with Euler vector field).

\subsection{Semisimple flat $F$-manifolds} Let $(M,\nabla,c,e)$ be an analytic flat $F$-manifold. A point $p\in M$ is called \emph{semisimple} if the algebra $(T_pM,\circ_p)$ is without nilpotent elements. This is equivalent to 
\begin{itemize}
\item the existence of idempotent vectors $\pi_1,\dots,\pi_n\in T_pM$, i.e. such that $\pi_i\circ_p\pi_j=\pi_i\delta_{ij}$, 
\item the existence of $v\in T_pM$ such that $v\circ_p\colon T_pM\to T_pM$ has simple spectrum. %
\end{itemize}Semisimplicity is an open property: if $p\in M$ is semisimple, then there exists an open neighborhood $\mc V$ of $p$, such that any point of $\mc V$ is semisimple. Moreover, if $\mc V$ is small enough, we have well defined local idempotent holomorphic vector fields $\pi_1,\dots, \pi_n\in\Gm(\mc V,\mathscr T_M)$.  See e.g. \cite[Ch.\,II]{Her02} for a detailed discussion.
\vskip3mm
Let $(H,\bm \Phi)$ be a formal flat $F$-manifold. Denote by $\circ_0$ the product on $H$ defined by structure constants $c^\al_{\bt\gm}(0):=\der^2_{\al\bt}\Phi|_{\bm t=0}$. We will say that $(H,\bm \Phi)$ is
\begin{itemize}
\item\emph{semisimple at the origin} if we have an isomorphisms of $k$-algebras $(H,\circ_0)\cong k^n$,
\item \emph{formally semisimple} if we have an isomorphism of $K$-algebras $(H,\circ)\cong K^n$.
\end{itemize}
Formal semisimplicity is thus equivalent to the existence of vectors $\pi_i\in H_K$ such that $\pi_i\circ\pi_j=\pi_i\delta_{ij}$.

\begin{lem}
A formal flat $F$-manifold is formally semisimple iff it is semisimple at the origin.
\end{lem}
\proof The proof of \cite[Lemma 4.2]{Cot20b} works verbatim.
\endproof

\begin{rem}
In both formal and analytic case, we have $e=\sum_{i=1}^n\pi_i$.
\end{rem}

\begin{prop}\label{picomm}
For both formal and analytic semisimple flat F-manifolds, the idempotent vectors $\pi_1,\dots,\pi_n$ are pairwise commuting, i.e. $[\pi_i,\pi_j]=0$. Hence, there exist local coordinates $\bm u=(u^1,\dots, u^n)$ such that $\pi_i=\frac{\der}{\der u^i}$ for $i=1,\dots,n$. The local coordinates $\bm u$ will be called \emph{canonical}.
\qed
\end{prop}

In the formal case, the functions $u^i$ are formal functions, i.e. elements of $K$.
Canonical coordinates $\bm u$ are defined up to permutations and shifts by constants. We set $\der_i:=\frac{\der}{\der u^i}$ for $i=1,\dots,n$. If an Euler vector field is given, then the shift freedom can be actually frozen.

\begin{prop}\label{canuE}A vector field $E\in\Gm(TM)$ satisfies $\frak L_Ec=c$ iff in canonical coordinates it has the form $E=\sum_j(u^j+c^j)\der_j$. 
Up to shifts of canonical coordinates $\bm u$, we have $E=\sum_ju^j\der_j$. \qed
\end{prop}
Hence the eigenvalues of the tensor $(E\circ)\in \Gm(\End TM)$ may and will be chosen as local canonical coordinates.

\begin{defn}\label{tcp}
A point $p\in M$ will be called 
\begin{itemize}
\item {\it tame} if the operator $E\circ_p\colon T_pM\to T_pM$ has simple spectrum
\item {\it coalescing}, otherwise.
\end{itemize}If a point is tame, then it is necessarily semisimple, and with pairwise distinct canonical coordinates. 
\vskip1mm
The same definition can adapted to the formal case, relatively at the origin $\bm t=0$, by looking at the spectrum of $E\circ_0\colon H\to H$.
\end{defn}
At coalescing points $p$, we have $\bm u(p)\in\Dl$ where $\Dl\subseteq\C^n$ denotes the {\it big diagonal}
\[\Dl:=\bigcup_{i\neq j}\{u_i=u_j\}.
\]
For a given flat $F$-manifold $M$, we define an {\it algebraic symmetry of $M$} to be a diffeomorphism $\phi\colon M\to M$ whose differential $d\phi_p\colon T_pM\to T_{\phi(p)}M$ is an isomorphism of unital algebras for any $p\in M$. Denote by ${\rm AlgSym}(M)$ the group of {algebraic symmetries} of $M$. 
\begin{rem}The group ${\rm AlgSym}(M)$ is never empty: if the $\nabla$-flat coordinates $\bm t$ are normalized so that $\frac{\der}{\der t^1}=e$, any translation $t^1\mapsto t^1+c$, with $c\in\C$, defines an algebraic symmetry. Moreover, notice that an algebraic symmetry $\phi$ is an isomorphism of flat $F$-structures if and only if $\phi$ preserves the $\nabla$-affine structure, namely $d\phi(\ker\nabla)=\ker(\phi^*\nabla)$.
\end{rem}
For semisimple flat $F$-manifolds, Proposition \ref{picomm} allows to compute the connected component ${\rm AlgSym}(M)_0$ of the identity.
\begin{prop}\label{groupiso}
If $M$ is a semisimple flat $F$-manifold, then the connected component ${\rm AlgSym}(M)_0$ %
of the identity is a commutative $n$-dimensional Lie group. Moreover, it acts locally transitively on $M$.
\end{prop}
\proof
The Lie algebra of ${\rm AlgSym}(M)$ %
can be identified with the Lie algebra of vector fields $X\in\Gm(TM)$ on $M$ such that $\frak L_Xc=0$. This is equivalent to the condition
\[[X,Y\circ Z]-[X,Y]\circ Z-[X,Z]\circ Y=0,\quad Y,Z\in\mathscr T_M.
\]In local canonical coordinates $\bm u$, set $X=\sum_iX^i\der_i$ and take $Y=Z=\der_j$. We have $\der_jX^k=0$ for all $j,k$. Hence, locally $X$ is a constant linear combination of the idempotent vector fields $\der_j$. In local canonical coordinates, the flow of $X$ reads as shifts $u^i\mapsto u^i+c^i$.
\endproof
\begin{rem}
On a semisimple flat $F$-manifold, we have two different affine structures: the first one is defined by $\nabla$, while the second one is defined by the atlas of canonical coordinates. As the proof of Proposition \ref{groupiso} shows, algebraic symmetries of a semisimple flat $F$-manifold preserve the second affine structure.
In general, the two affine structures are not compatible: this only happens when the canonical coordinates $\bm u$ are $\nabla$-flat. In such a case, the oriented associativity potentials take the form $F^i(\bm u)=\frac{1}{2}u_i^2$, for $i=1,\dots,n$. The resulting flat $F$-manifold is trivial: it is isomorphic to $\C^n$ equipped with the product $\C$-algebra structure $(\C,+,\cdot,1)^{\times n}$, seen as a flat $F$-manifold. \end{rem}
\begin{rem}
In \cite{AL13} a notion of {\it bi-flat} $F$-manifold is introduced. This consists of the datum of two flat $F$-manifolds structures $(\nabla,\circ, e)$ and $(\nabla^*,*,E)$, on the same manifold $M$, satisfying the following compatibility conditions:
\begin{enumerate}
\item $E$ is $\circ$-invertible on $M$,
\item $\frak L_E(\circ)=\circ$,
\item $X*Y=E^{-1}\circ X\circ Y$, for all $X,Y\in\mathscr T_M$,
\item $(d_{\nabla}-d_{\nabla^*})(X\circ)=0$, for all $X\in\mathscr T_M$ (here $d_\nabla$ is the exterior covariant derivative).
\end{enumerate}
In \cite{AL17}, in the tame semisimple case (pairwise distinct canonical coordinates), it is proved that a bi-flat $F$ structure is actually {\it equivalent} to the datum of a homogeneous flat $F$-manifold with invertible Euler vector field $E$. We also refer the reader to \cite{KMS18}, where the identification of bi-flat $F$-structures and homogeneous flat $F$-structures (in {\it loc.\,cit.\,\,}they are called {\it Saito structures}) is established even beyond the semisimple case. Furthermore, in \cite{KMS18} Dubrovin's almost duality formalism for Frobenius manifolds is extended to homogeneous flat $F$-manifolds.
\end{rem}

\subsection{Irreducible flat $F$-manifolds}
If $M_1,M_2$ are two flat $F$-manifolds, their product $M_1\times M_2$ is naturally equipped with a flat $F$-structure, called the {\it sum} $M_1\oplus M_2$. If $M_1,M_2$ are homogenous, then also $M_1\oplus M_2$ is homogeneous.
\vskip2mm
We say that a flat $F$-manifold $M$ is {\it irreducible} if no pointed germ $(M,p)$ is locally isomorphic to a pointed sum $(M_1\oplus M_2,p')$. 
\vskip2mm
In the semisimple homogeneous case, we have the following characterization of irreducibility.

\begin{thm}\label{eulerunique}Let $M$ be a formal/analytic semisimple and homogeneous flat $F$-manifold. The following conditions are equivalent:
\begin{enumerate}
\item $M$ is irreducible;
\item if $E_1,E_2\in\Gm(TM)$ are two Euler vector fields, then $E_2=E_1-\la e$ for some $\la\in\C$.\qed
\end{enumerate}
\end{thm}
The proof of this result can be found in Appendix \ref{theuelrid}. 
\vskip1mm

\begin{rem}
Theorem \ref{eulerunique} underlines how much selective is the condition $\nabla\nabla E=0$. In the category of $F$-manifolds (not necessarily flat), Euler vector fields are defined\footnote{A more general notion of Euler vector field of {\it weight} $d\in\C$ is discussed in \cite{Her02,man}: these are vector fields $E$ such that $\frak L_Ec=d\cdot c$. If $d\neq 0$, one can always rescale $E$ in order to be of weight 1.} by the condition $\frak L_Ec=c$ only. For semisimple $F$-manifolds, given an Euler vector field $E$, all other Euler fields are of the form $E+\sum_{i=1}^n\C\pi_i$. See \cite[Ex.\,2.12(ii)]{Her02}.
\end{rem}

\section{Extended deformed connections}\label{sec3}%

\subsection{$\nabla^z$-flat coordinates and oriented associativity potentials} In both the analytic and the formal context (over $k$), we can look for $\nabla^z$-flat coordinates of the flat $F$-structure, i.e. functions $\tilde t^\al(\bm t,z)$ such that $\nabla^zd\tilde t^\al=0$. Assume they are of the form
\[\tilde t^\al(\bm t,z):=\sum_{p=0}^\infty h_{p}^\al(\bm t)z^p\in k[\![\bm t,z]\!],\quad h_{0}^\al(\bm t)=t^\al,\quad \al=1,\dots,n.
\]
\begin{thm}
The functions $h^\al_p$ satisfy the recursive equations
\[h^\al_0(\bm t)=t^\al,\quad \der_\gm\der_\bt h^\al_{p+1}=c^\la_{\gm\bt}\der_\la h^\al_p,\quad p\in\N.
\]
\end{thm}
\proof
The $\nabla^z$-flatness equations for a one-form $\xi=\xi_\al dt^\al$ are $\der_\gm\xi_\bt=z c^\la_{\gm\bt}\xi_\la$.
\endproof

\begin{cor}\label{corhf}
The functions $h^\al_1(\bm t)$ equal the oriented associativity potentials $F^\al(\bm t)$ up to linear terms.
\end{cor}
\proof
We have $\der_\gm\der_\bt h^\al_1=c^\al_{\bt\gm}$.
\endproof

\subsection{Family of extended deformed connections}Following \cite[Section 3]{Man05},  we introduce a one-parameter family $(\widehat\nabla^\la)_\la$ of flat connections, which ``rigidify'' the family $(\nabla^z)_{z}$.  See also \cite[Section 4.3]{BB19}, \cite[Section 1.4]{ABLR20}.
\vskip3mm
\noindent{\bf Analytic case. }Let $(M,\nabla,c,e,E)$ be a homogenous flat $F$-manifold. Introduce the $(1,1)$-tensors $\mc U,\mu^{(\la)}\in\Gm({\rm End}(TM))$, with $\la\in\C$, by the formulae
\[\mc U(X)=E\circ X,\quad \mu^{\la}(X):=(1-\la)X-\nabla_XE,\quad X\in\mathscr T_M.
\]By equation \eqref{evf}, in $\bm t$-coordinates we have $\mu^{\la}={\rm diag}(q_1-\la,\dots, q_n-\la)$.
\vskip2mm
Denote by $\pi\colon M\times \C^*\to M$ the canonical projection on the first factor. If $\mathscr T_M$ denotes the tangent sheaf of $M$, then $\pi^*\mathscr T_M$ is the sheaf of sections of $\pi^*TM$, and $\pi^{-1}\mathscr T_M$ is the sheaf of sections of $\pi^*TM$ constant along the fibers of $\pi$.  All the tensors $c,e,E,\mc U,\mu$ can be lifted to the pullback bundle $\pi^*TM$, and we denote these lifts with the same symbols. Consequently, also the connection $\nabla$ can be uniquely lifted on $\pi^*TM$ in such a way that $\nabla_{\frac{\der}{\der z}}Y=0$ for $Y\in\pi^{-1}\mathscr T_M$.

The extended deformed connection $\widehat\nabla^{\la}$, with $\la\in\C$, is the connection on $\pi^*TM$ defined by the formulae
\beq
\label{aedc1}
\widehat\nabla^{\la}_{\frac{\der}{\der t^\al}}Y=\nabla_{\frac{\der}{\der t^\al}}Y+z\frac{\der}{\der t^\al}\circ Y,\quad\quad
\widehat\nabla^{\la}_{\frac{\der}{\der z}}Y=\nabla_{{\frac{\der}{\der z}}}Y+\mc U(Y)-\frac{1}{z}\mu^{\la}(Y),
\eeq
where $Y\in \pi^*\mathscr T_M$. 
\vskip3mm
\noindent{\bf Formal case. }Let $k$ be a commutative $\Q$-algebra and $(H,\bm \Phi, E)$ a formal homogeneous flat $F$-manifold  over $k$. Denote by $k(\!(z)\!)$ the $k$-algebra of formal Laurent series in an auxiliary indeterminate $z$. Set $K(\!(z)\!):=k[\![\bm t]\!](\!(z)\!)$ and $H_{K(\!(z)\!)}:=H\otimes_k K(\!(z)\!)$.
\vskip2mm
In what follows we assume that the $K$-linear operator $\nabla^0E\colon {\rm Der}_k(K)\cong H_K\to H_K$ is (diagonalizable and) in diagonal form in the basis $(\Dl_1,\dots,\Dl_n)$. Define the $K$-linear operators $\mc U,\mu^\la$, with $\la\in k$, by the formulae
\bea
&\mc U\colon H_K\to H_K,\quad &X\mapsto E\circ X,\\
&\mu^\la\colon {\rm Der}_k(K)\cong H_K\to H_K,\quad &X\mapsto (1-\la)-\nabla_XE.
\eea
All the tensors $\circ, \mc U,\mu^\la$ can be $K(\!(z)\!)$-linearly extended to $H_{{K(\!(z)\!)}}$. We will denote such an extension by the same symbols.
\vskip3mm
Denote by ${\rm Diff}_1(H_{K(\!(z)\!)},H_{K(\!(z)\!)})$ the set of morphisms $\mathscr D\in \Hom_k(H_{K(\!(z)\!)},H_{K(\!(z)\!)})$ such that
\[ab\mathscr D(p)-b\mathscr D(ap)-a\mathscr D(bp)+\mathscr D(abp)=0,\quad a,b\in K(\!(z)\!),\quad p\in H_{K(\!(z)\!)}.
\]Both ${\rm Der}_k(K(\!(z)\!))$ and ${\rm Diff}_1(H_{K(\!(z)\!)},H_{K(\!(z)\!)})$ are naturally equipped with an $K(\!(z)\!)$-module structure.
\vskip3mm 
The extended deformed connection $\widehat\nabla^\la\colon{\rm Der}_k(K(\!(z)\!))\to {\rm Diff}_1(H_{K(\!(z)\!)},H_{K(\!(z)\!)})$ is the $K(\!(z)\!)$-linear operator defined by the formulae 
\[\widehat\nabla^{\la}_{\frac{\der}{\der t^\al}}X=\nabla^z_{\frac{\der}{\der t^\al}}X,\quad\quad \widehat\nabla^{\la}_{\frac{\der}{\der z}}X={\frac{\der}{\der z}}X+\mc U(Y)-\frac{1}{z}\mu^{\la}(X),
\]
where $Y\in H_{K(\!(z)\!)}$.
\vskip3mm
In both the analytic and formal pictures, the following result holds.
\begin{thm}\label{ffedc}
The connection $\widehat\nabla^\la$ is flat for any $\la\in\C$ (resp $\la\in k$).
\end{thm}
\proof The flatness of $\widehat\nabla^\la$ is equivalent to the following conditions:
$\der_{\al}c_{\bt\gm}^\dl$ is completely symmetric in $(\al,\bt,\gm)$,
the product $\circ$ is associative,
$\nabla\nabla E=0$,
and $\frak L_Ec=c$.
This can be easily checked by a straightforward computation.
\endproof

\begin{rem}\label{symderz}
For $\la=\frac{1}{2}d$, the connection $\widehat\nabla^{\la}$ equals the extended deformed connection $\widehat\nabla$ as defined by Dubrovin \cite{Dub96,Dub98,Dub99}. In that case, the tensor $\mc U$ (resp. $\mu=\mu^{(\frac{d}{2})}$) is $\eta$-self-adjoint (resp. $\eta$-anti-self-adjoint). It follows that if $\zeta_1,\zeta_2\in\pi^*\mathscr T_M$ are two $\widehat\nabla$-flat vector fields, then the pairings $\langle\zeta_1,\zeta_2\rangle_\pm:=\eta\left(\zeta_1(\bm t, e^{\pm\pi\sqrt{-1}}z),\zeta_2(\bm t,z)\right)$ do not depend on $(\bm t,z)$. See also \cite[Section 2]{CDG20}.
\end{rem}

\begin{rem}\label{remt1}
In both the analytic and formal cases, we have $\mu^\la(e)=-\la e$. This follows from the torsionless of $\nabla=\nabla^0$, the $\nabla$-flatness of $e$, and Lemma \ref{commeE}. Without loss of generality, we can then assume that the flat coordinate $t^1$ is such that $\frac{\der}{\der t^1}=e$ (analytic case), and that $\Dl_0=e$ (formal case). In such a case, the parameter $q_1$ in \eqref{evf} satisfies $q_1=0$.
\end{rem}

\subsection{$\widehat\nabla^\la$-flat covectors}
In both the analytic and formal pictures, the extended connections $\widehat\nabla^\la$ induce connections on the whole tensor algebra of $\pi^*TM$ (resp. $H_{K(\!(z)\!)}$). So, for example, let $\xi$ denote a $\widehat\nabla^\la$-flat section of the bundle $\pi^*(T^*M)$. In the co-frame $(dt^\al)_{\al=1}^n$, the equation $\widehat\nabla^\la\xi=0$ can be written, in more convenient matrix notations, as the joint system of differential equations
\beq\label{jfs}
\frac{\der \xi}{\der t^\al}=z\mc C_\al^T\xi,\quad \frac{\der \xi}{\der z}=\left(\mc U-\frac{1}{z}\mu^\la\right)^T\xi,
\eeq
where $\xi=(\xi_1,\dots,\xi_n)^T$ is a column vector of components w.r.t.\,\,$(dt^\al)_\al$, and 
\[(\mc C_\al)^\gm_\bt:=c_{\al\bt}^\gm, \quad (\mc U)^\bt_\al=E^\eps c^\bt_{\al\eps},\quad (\mu^\la)^\bt_\al=(q_\al-\la)\dl_{\al\bt}.
\]
We will refer to the second of equations \eqref{jfs} as the $\der_z$-{\it equation} of the flat $F$-manifold.

\subsection{Matrices $\wPsi,\wV_i,\wV^\la,\wGm$}\label{secwpsietc} Assume that $(M,\nabla,c,e,E)$ is a semisimple homogeneous flat $F$-manifold, 
and introduce the Jacobian matrix $\wPsi$ by \[\wPsi_{\al}^i:=\frac{\der u^i}{\der t^\al},\quad i,\al=1,\dots,n.\]

\begin{rem}\label{remtildepsi1}
If the flat coordinate $t^1$ is chosen so that $\frac{\der}{\der t^1}=e$, then we have $\Tilde\Psi_1^i=1$ for any $i=1,\dots,n$.
\end{rem} In canonical coordinates $\bm u$, under the gauge transformation $\tilde x=(\wPsi^{-1})^T\xi$, the system \eqref{jfs} becomes
\beq
\label{jfssw}
\frac{\der \tilde x}{\der u^i}=\left(zE_i-\widetilde V_i\right)^T\tilde x,\quad \frac{\der \tilde x}{\der z}=\left(U-\frac{1}{z}\widetilde V^\la\right)^T \tilde x,
\eeq
where 
\[(E_i)_{jk}=\dl_{ij}\dl_{ik},\quad U:={\rm diag}(u^1,\dots, u^n),\quad \wV_i:=\der_i\wPsi\cdot \wPsi^{-1},\quad \wV^\la:=\wPsi\cdot\mu^\la\cdot\wPsi^{-1}.\]
\vskip1,5mm
For any matrix $A\in M_n(\C)$, we denote by $A'$ and $A''$ the diagonal and off-diagonal part of $A$, respectively. So, we have the decomposition: 
\[A=A'+A'',\quad (A')_{ij}=0,\quad (A'')_{ii=0},\quad i,j=1,\dots,n,\quad i\neq j.\]
\begin{prop}\label{pids}The following facts hold true for both formal and analytic semisimple homogeneous flat $F$-manifolds.
\begin{enumerate}
\item There exists an off-diagonal matrix $\wGm$ such that 
\bean
\label{strvi}
&\wV_i=\wV_i'+[\wGm,E_i],&\quad i=1,\dots,n,\\
\label{strv}
&\wV^\la=(\wV^\la)'+[\wGm, U].&
\eean
In particular, $\wGm^i_j=-(\wV_i)^i_j$ for $i\neq j$.
\item We have 
\beq
\label{dvvi}
[E_i,\wV^\la]=[U,\wV_i],\quad \der_i\wV^\la=[\wV_i,\wV^\la].
\eeq
\item The diagonal entries of the matrix $\wV^\la$ are constant w.r.t.\,\,$\bm u$.
\item We have $\der_i\wV_j'=\der_j\wV_i'$.
\end{enumerate}
\end{prop}
\proof
The compatibility condition $\der_i\der_j=\der_j\der_i$ implies the constraints
\bean
\label{ids1}
&\der_i\wV_j-\der_j\wV_i=[\wV_i,\wV_j],\quad
[E_i,E_j]=0,\\
\label{ids2}
&[E_i,\wV_j]=[E_j,\wV_i],
\eean
Identities \eqref{ids1} are trivially satisfied, by definition of the matrices $E_i$ and $\wV_i$. From \eqref{ids2}, we deduce
\[(\dl_{jh}-\dl_{jk})(\wV_i)^h_k=(\dl_{ih}-\dl_{ik})(\wV_j)^h_k\quad\underset{h=j,\, j\neq k}{\Longrightarrow}\quad (\wV_i)^j_k=(\dl_{ij}-\dl_{ik})(\wV_j)^j_k=[\wGm, E_i]^j_k,
\]
where $\wGm=(\wGm^j_k)_{j,k}$ is defined by $\wGm^j_k:=-(\wV_j)^j_k$. This proves \eqref{strvi}.
The compatibility condition $\der_i\der_z=\der_z\der_i$ implies the constraints 
\bean
\label{ids3}
&[E_i,U]=0,\\
\label{ids4}
&[E_i,\wV^\la]=[U,\wV_i],\quad \der_i\wV^\la=[\wV_i,\wV^\la].
\eean
Identity \eqref{ids3} is trivially satisfied. Identity \eqref{strv} follows from \eqref{strvi} and the first of \eqref{ids4}.
The constancy of $(\wV^\la)'$ follows from the second identity of \eqref{ids4}, and equations \eqref{strvi}, \eqref{strv}.
Finally, from the first identity of \eqref{ids1} we have $\der_i\wV_j'-\der_j\wV_i'=[\wV_i,\wV_j]'=0$, by \eqref{strvi}, \eqref{strv}.
\endproof

\subsection{Darboux-Tsarev equations, and conformal dimensions} Let us introduce the Christoffel symbols $K^h_{ij}$ by $\nabla_{\der_i}\der_j=\sum_h K^h_{ij}\der_h$.

\begin{lem}\label{kvi}
We have $K^h_{ij}=-(\wV_i)^h_j$.
\end{lem}
\proof
The claim follows from the following computation:
\[\pushQED{\qed} 
\nabla_{\der_i}\der_j=\nabla_{\der_i}[(\wPsi^{-1})^\al_j\der_\al]=[\der_i(\wPsi^{-1})^\al_j]\der_\al=-\sum_{\ell}(\wPsi^{-1})^\al_\ell\ \der_i\wPsi^\ell_\bt\ (\Psi^{-1})^\bt_j\der_\al=-\sum_\ell (\wV_i)^\ell_j \der_\ell.\qedhere
\popQED
\]

\begin{prop}\label{DTeqs}
The following identities hold true:
\bean
\label{cs1}
&K^h_{ij}=0,\quad &i,j,h\text{ distinct,}\\
\label{cs2}
&K^i_{ij}=K^i_{ji}=-K^i_{jj}=\wGm^i_j,\quad &i\neq j,\\
\label{cs3}
&K^i_{ii}=-\sum_{h\neq i}\wGm^i_h.
\eean
Moreover, the functions $\wGm^i_j$ satisfy the Darboux-Tsarev equations 
\bean
\label{ts1}
&\der_k\wGm^i_j=-\wGm^i_j\wGm^i_k+\wGm^i_j\wGm^j_k+\wGm^i_k\wGm^k_j\quad &i,j,k\text{ distinct,}\\
\label{ts2}
&\sum_k\der_k\wGm^i_j=0,\quad &i\neq j.
\eean
\end{prop}
\proof In canonical coordinates we have $c^i_{jk}=\dl^i_j\dl^i_k$. Consequently, $(\nabla_\ell c)^i_{jk}=\sum_p\dl^p_j\dl^p_k K^i_{p\ell}-\dl^i_k K^i_{\ell j}-\dl^i_j K^i_{\ell k}$. We have $K^h_{ij}=K^h_{ji}$ because $\nabla$ is torsion free. Hence, from the symmetry $(\nabla_\ell c)^i_{jk}=(\nabla_j c)^i_{\ell k}$, one obtains \eqref{cs1}, and the first two equalities of \eqref{cs2}. The equality $K^i_{ij}=\wGm^i_j$ follows from Lemma \ref{kvi}. %
Equation \eqref{cs3} follows from the condition $\nabla e=0$. By flatness of $\nabla$, the components $R^i_{jk\ell}$ of the Riemann tensor equal zero. By definition we have $R^i_{jk\ell}=\der_k K^i_{j\ell}-\der_\ell K^i_{jk}+\sum_pK^i_{kp}K^p_{j\ell}-\sum_pK^i_{\ell p}K^p_{kj}$. Darboux-Tsarev equation \eqref{ts1} is equivalent to $R^i_{ji\ell}=0$.
From the identity $\nabla_{\der_i}\der_\al=\sum_j(\der_i \wPsi^j_\al)\der_j+\sum_j\wPsi^j_\al\nabla_{\der_i}\der_j$, summing over $i$, we obtain
\bean
&0=\nabla_e\der_\al=\sum_i\nabla_{\der_i}\der_\al=\sum_j\left(\sum_i \der_i \wPsi^j_\al\right)\der_j+\sum_j\wPsi^j_\al\underbrace{\nabla_e\der_j}_{\nabla_{\der_j}e+[e,\der_j]=0},\\
\label{scpsi}
&\Longrightarrow \sum_i \der_i \wPsi=0.
\eean
We have
\[\der_k\wGm^i_j=-\der^2_{ki}\wPsi^i_\al(\wPsi^{-1})^\al_j+\der_i\wPsi^i_\al\sum_h(\wPsi^{-1})^\al_h\ \der_k\wPsi^h_\gm\ (\wPsi^{-1})^\gm_j.
\]Summing over $k$, and using \eqref{scpsi}, we obtain \eqref{ts2}.
\endproof

\begin{cor}\label{vigm}
For $i=1,\dots, n$, the matrix $\wV_i$ has the following structure
\[\pushQED{\qed} 
\wV_i=\left(\begin{array}{cccccccc}
-\wGm^1_i&&&&\wGm^1_i&&&\\
&-\wGm^2_i&&&\wGm^2_i&&&\\
&&\ddots&&\vdots&&&\\
&&&-\wGm^{i-1}_i&\wGm^{i-1}_i&&&\\
-\wGm^i_1&-\wGm^i_2&\dots&-\wGm^{i}_{i-1}&\sum_{h\neq i}\wGm^i_h&-\wGm^i_{i+1}&\dots&-\wGm^i_n\\
&&&&\wGm^{i+1}_i&-\wGm^{i+1}_i&&\\
&&&&\vdots&&\ddots&\\
&&&&\wGm^n_i&&&-\wGm^n_i
\end{array}\right).\qedhere
\popQED
\]
\end{cor}

\begin{prop}\label{strvla}
We have $\wV^\la=-\la\cdot{\bf 1}+\sum_iu^i\wV_i$.
\end{prop}
\proof
For any $i$, we have
\[\nabla_{\der_i}E=\nabla_{\der_i}\sum_ju^j\der_j=\der_i+\sum_ju^j\nabla_{\der_i}\der_j=\der_i+\sum_{j,h}u^jK_{ij}^h\der_h=\der_i+\sum_{j,h}u^jK_{ji}^h\der_h.
\]By Lemma \ref{kvi} one concludes.
\endproof
\begin{cor}The following identities hold true:
\bean
\label{hv}
\sum_ju^j\der_j\wV_i&=&-\wV_i,\\
\label{hgm}
\sum_j u^j\der_j\wGm&=&-\wGm,\\
\label{hgm1}
(u^i-u^j)\der_i\wGm^i_j&=&\sum_{\ell\neq i,j}(u^j-u^\ell)\{-\wGm^i_j\wGm^i_\ell+\wGm^i_j\wGm^j_\ell+\wGm^i_\ell\wGm^\ell_j\}-\wGm^i_j,\\
\label{hgm2}
(u^j-u^i)\der_j\wGm^i_j&=&\sum_{\ell\neq i,j}(u^i-u^\ell)\{-\wGm^i_j\wGm^i_\ell+\wGm^i_j\wGm^j_\ell+\wGm^i_\ell\wGm^\ell_j\}-\wGm^i_j.
\eean
\end{cor}
\proof
Equation \eqref{hv} follows from the second equation of \eqref{dvvi} and the first equation of \eqref{ids1}. Equation \eqref{hgm} is easily deduced. Equations \eqref{hgm1} and \eqref{hgm2} follow from \eqref{ts1}, \eqref{ts2}, and \eqref{hgm}.
\endproof

\begin{rem}
In this section we started from a given semisimple flat $F$-manifold and we obtained a solution $\widetilde\Gm^i_j$ of the Darboux-Tsarev equations. This is also proved in \cite{AL15}. In {\it loc.\,\,cit.}, it is proved that the opposite construction works as well: %
the datum of
\begin{itemize}
\item a solution $\widetilde\Gm^i_j$ of  \eqref{ts1},\eqref{ts2},
\item the connection $\nabla$ with Christoffel symbols $K^i_{jk}$ given by \eqref{cs1},\eqref{cs2},\eqref{cs3},
\item the structure constants $c^i_{jk}=\dl^i_j\dl^i_k$,
\item the vector field $e:=\sum_i\der_i$,
\end{itemize}
locally defines a (tame) semisimple flat $F$-manifold structure on $\C^n\setminus\Dl$.
\end{rem}

\vskip2mm
\noindent{\bf Conformal dimensions. }
By Propositions \ref{pids}-(3), and \ref{strvla}, there exist complex numbers $\dl_1,\dots,\dl_n\in\C$ such that
\[(\wV^\la)'=-\la\cdot{\bf 1}+{\rm diag}(\dl_1,\dots,\dl_n).
\]

\begin{defn}The numbers $\dl_1,\dots,\dl_n$ are called {\it conformal dimensions} of the (formal or analytic) semisimple flat $F$-manifold. We will say that a (formal/analytic) semisimple flat $F$-manifold is {\it conformally resonant} if $\dl_i-\dl_j\in\Z\setminus{0}$ for some $i,j$.
\end{defn}

\begin{rem}
We have $\dl_i=\sum_ku^k(\wV_k)^i_i=\sum_{k\neq i}(u^i-u^k)\wGm^i_k$.
\end{rem}

\begin{rem}
In the case of Frobenius manifolds, all conformal dimensions equal $\frac{1}{2}d$, where $d$ is the conformal dimension of equation \eqref{kill}. In particular, a Frobenius manifold is never conformally resonant. 
\end{rem}

\subsection{Lam\'e coefficients, matrices $\Psi,V_i,V^\la,\Gm$, and Darboux-Egoroff equations} For any $j=1,\dots,n$ define the one form
\[\om_j(\bm u):=-\sum_{i=1}^n\wV_i(\bm u)^j_jdu^i.
\]
\begin{prop}
The one-forms $\om_j$ are closed. There locally exist functions $H_j(\bm u)$ such that
\[d\log H_j=\om_j,\quad j=1,\dots,n.
\]
\end{prop}
\proof
It follows from point (4) of Proposition \ref{pids}.
\endproof
The functions $H_j$ are called \emph{Lam\'e coefficients}, and they are defined up to scalar rescaling $H_j\mapsto \la_j H_j$, $\la_j\in\C^*$.\vskip2mm
Arrange the Lam\'e coefficients in the diagonal matrix $H:={\rm diag}(H_1,\dots, H_n)$, and define the matrices
\[\Psi:=H\wPsi,\quad  V_i:=H\wV_i H^{-1}+\der_i H\cdot H^{-1},\quad V^\la:=H\wV^\la H^{-1},\quad \Gm:= H\wGm H^{-1}.
\]
\begin{prop}\label{propLame}
The functions $H_1,\dots, H_n$ satisfy the following system 
\[
\der_jH_i=\Gm^i_jH_j,\quad i\neq j,\quad
\der_iH_i=-\sum_{k\neq i}\Gm^i_kH_k,\quad
\sum_ju^j\der_jH_i=-\dl_iH_i.
\]
\end{prop}
\proof
It easily follows from the definitions and identities of the previous section.
\endproof
\begin{rem}\label{magicrk2}
If the $\nabla$-flat coordinate $t^1$ is such that $\frac{\der}{\der t^1}=e$, as in Remark \ref{remt1}, then we have $\Psi^i_1=H_i$ for any $i=1,\dots,n$. In other words, the Lam\'e coefficients can be read from the first column of the $\Psi$-matrix. Moreover, we have
\beq
V^\la\bm H=-\la\bm H,\quad \bm H=(H_1,\dots,H_n)^T.
\eeq
\end{rem}
Under the gauge transformation $x=(H^{-1})^T\tilde x$, the system \eqref{jfssw} becomes
\beq
\label{jfss}
\frac{\der x}{\der u^i}=\left(zE_i-V_i\right)^Tx,\quad \frac{\der x}{\der z}=\left(U-\frac{1}{z}V^\la\right)^Tx.
\eeq

\begin{prop}\label{fundprop1}
The following identities hold true:
\bean
&V_i=[\Gm,E_i],&\quad \\
&V^\la=(V^\la)'+[\Gm,U],&\quad (V^\la)'=(\wV^\la)'={\rm diag}(\dl_1-\la,\dots,\dl_n-\la),\\
&\label{compmag}\der_i\Psi\cdot \Psi^{-1}=V_i,&\quad [E_i,V^\la]=[U,V_i],\quad \der_i V^\la=[V_i,V^\la].
\eean
\end{prop}
\proof
It easily follows from the definitions and identities of the previous section.
\endproof

\begin{prop}\label{fundprop2}
The matrix $\Gm$ satisfies the Darboux-Egoroff equations
\bean
\label{DE1}
&\der_k\Gm^i_j=\Gm^i_k\Gm^k_j,\quad &\text{$i,j,k$ distinct,}\\
\label{DE2}
&\sum_{k}\der_k\Gm^i_j=0,\quad&i\neq j,\\
\label{DE3}
&\sum_k u^k\der_k\Gm^i_j=(\dl_j-\dl_i-1)\Gm^i_j,\quad&i\neq j,\\
\label{midgm}
&(u^j-u^i)\der_i\Gm^i_j=\sum_{k\neq i,j}(u^k-u^j)\Gm^i_k\Gm^k_j-(\dl_j-\dl_i-1)\Gm^i_j,\\
\label{midgm2}
&(u^i-u^j)\der_j\Gm^i_j=\sum_{k\neq i,j}(u^k-u^i)\Gm^i_k\Gm^k_j-(\dl_j-\dl_i-1)\Gm^i_j.
\eean
\end{prop}
\proof
It easily follows from the definitions, the Darboux-Tsarev system for $\wGm$, and the homogeneity identities \eqref{hgm} of the previous section.
\endproof

\begin{rem}
For $n=3$, the Darboux-Egoroff joint system of equations \eqref{DE1}, \eqref{DE2}, \eqref{DE3} is equivalent to the full family of Painlev\'e equations PVI. See remarkable formulas of \cite[Theorem 4.1]{Lor14}.
\end{rem}

\begin{rem}\label{defrob}
In the case of Frobenius manifolds, there is a canonical choice for the Lam\'e coefficients: $H_i=\eta(\der_i,\der_i)^{\frac{1}{2}}$ for $i=1,\dots, n$. The resulting coefficients $\Gm^i_j$ are the rotation coefficients of the metric $\eta$. They satisfy the further symmetry condition $\Gm^i_j=\Gm^j_i$.
\end{rem}

\begin{rem}
In the light of Remark \ref{defrob}, given a semisimple flat $F$-manifold with a fixed choice of the Lam\'e coefficients $H_i$, we can define a metric by $\eta:=\sum_{i=1}^n H_i^2du^i$. Such a metric clearly is compatible with the product, in the sense that the second of equations \eqref{mfrob} is satisfied. The flatness of $\eta$ is the obstruction for the flat $F$-manifold to be actually Frobenius. For a more invariant description of the metric $\eta$, see \cite[Prop.\,1.8]{ABLR20}.
\vskip1mm
In \cite{ABLR20} a further notion of {\it semisimple Riemannian $F$-manifold} is introduced. In {\it loc.\,cit.} it is also proved the local equivalence of semisimple flat $F$-manifolds and semisimple Riemannian $F$-manifolds. Notice that the notion of semisimple Riemannian $F$-manifolds given in \cite{ABLR20} relaxes the axioms of analog structures introduced in \cite{DS11,LPR11}. See also the recent preprint \cite{ABLR21}.
\end{rem}

From a given homogeneous semisimple flat $F$-manifold, we obtained a joint system of equations \eqref{jfss}. In the analytic case, such a joint system defines\footnote{This will be explained in details the next section. } an {\it isomonodromic system}, because of integrability equations \eqref{compmag}. 

Vice-versa, one can start from such an isomonodromic system to construct the homogeneous semisimple flat $F$-manifold structure. This is exactly the point of view of the definition of {\it Saito structures without metric} given in \cite[Ch.\,VII, \S 1.a]{Sab07}. 

Notice that one can actually work with a companion Fuchsian system, obtained via a Fourier-Laplace transform. This is the point of view of \cite{KMS20}, in which Saito structures without metric are constructed on the space of isomonodromic deformation parameters for extended Okubo systems. 

In both \cite{Sab07,KMS20}, coalescences of the parameters of deformations (the entries of $U={\rm diag}(u^1,\dots, u^n)$ of equation \eqref{jfss}) are not taken into account. In our equation \eqref{jfss}, on the contrary, we allow coalescences, provided that the matrix $\Psi$ is not singular, i.e. provided that the geometric point of the flat $F$-structure is semisimple. In the recent paper \cite{Guz20}, D.\,Guzzetti extended the results of \cite{BJL81,Guz16} to the case of Fuchsian systems with {\it confluent  singularities}. Furthermore, in \cite{Guz20} a notion of {\it isomonodromic Laplace transform} is introduced: with such an analytic tool, the study of the correspondence  
\[\text{Monodromy data
of an irregular system}\longleftrightarrow\text{Monodromy data of a Fuchsian system,}
\]originally developed in \cite{BJL81}, has been extended to the isomonodromic case (possibly with coalescences/confluences). This also gives a new proof of the results of \cite{CG18,CDG1}.
\vskip1mm
The point of view of the current paper differs from the perspective of \cite{Sab07,KMS20}, via a Riemann-Hilbert correspondence. In Section \ref{RHBsec}, we will show a one-to-one correspondence between (local isomorphism classes of) homogeneous semisimple flat $F$-structures and solvable Riemann-Hilbert-Birkhoff problems.

\section{Monodromy moduli of admissible germs of semisimple flat $F$-manifolds}\label{sec4}
\subsection{$\mu$-nilpotent operators and $\mu$-parabolic group}\label{muop} Let $(V,\mu)$ be the datum of a $n$-dimensional complex vector space, and a diagonalizable operator $\mu\colon V\to V$. Denote by ${\rm spec}(\mu)=(\mu_1,\dots, \mu_n)$ the spectrum of $\mu$, and by $V_{\mu_\al}$ the eigenspace corresponding to the eigenvalue $\mu_\al$.
\vskip3mm
We say that $A\in{\rm End}(V)$ is $\mu$-{\it nilpotent} if 
\[AV_{\mu_\al}\subseteq\bigoplus_{m\geq 1} V_{\mu_\al+m}\quad \text{for all } \mu_\al\in{\rm spec}(\mu).
\]
In particular such an operator is nilpotent in the usual sense. Denote by $\frak c(\mu)$ the set of all $\mu$-nilpotent operators. It is easy to see that the set $\frak c(\mu)$ is a Lie algebra w.r.t.\,\,the commutator $[-,-]$ in ${\rm End}(V)$. 
We can decompose a $\mu$-nilpotent operator $A$ in components $A_k$, $k\geq 1$, such that
\[A_kV_{\mu_\alpha}\subseteq V_{\mu_\alpha+k}\quad\text{for any }\mu_\alpha\in\operatorname{spec}(\mu),
\]so that the following identities hold:
\[z^\mu Az^{-\mu}=A_1z+A_2z^2+A_3z^3+\dots,\quad [\mu,A_k]=kA_k\quad\text{for }k=1,2,3,\dots.
\]

\begin{lem}Let $(V,\mu)$ as above, and let us fix a basis $(v_i)_{i=1}^n$ of eigenvectors of $\mu$.
\begin{enumerate}
\item The operator $A\in\End(V)$ is $\mu$-nilpotent if and only if its associate matrix w.r.t. the basis  $(v_i)_{i=1}^n$ satisfies the condition
$(A)^\alpha_\beta=0\text{ unless }\mu_\alpha-\mu_\beta\in\mathbb N^*$.
\item If $A\in\End(V)$ is a $\mu$-nilpotent operator, then the matrices associated with its components $(A_k)_{k\geq 1}$ w.r.t. the basis $(v_i)_{i=1}^n$ satisfy the condition $(A_k)^\alpha_\beta=0\text{ unless }\mu_\alpha-\mu_\beta=k,$ with $k\in\mathbb N^*$.\qed
\end{enumerate}
\end{lem}
Define the $\mu$-{\it parabolic} group to be the Lie group $\mc C(\mu)$ of operator $G\in GL(V)$ such that $G={\bf 1}+A$, with $A\in\frak c(\mu)$. We have the canonical identification of Lie algebras $T_{\bf 1}\mc C(\mu)=\frak c(\mu)$, and the canonical adjoint action ${\rm Ad}\colon \mc C(\mu)\to {\rm Aut}\,\frak c(\mu)$ defined by
\[{\rm Ad}_G(A):=GA G^{-1},\quad G\in\mc C(\mu),\quad A\in\frak c(\mu).
\]
\begin{rem}\label{caniso}
Consider the space $V^*:=\Hom_{\C}(V,\C)$. Each $f\in\End_{\C}(V)$ induces a dual map $f^*\in\End_\C(V^*)$, defined by $f^*(w):=w\circ f$, where $w\in V^*$. This defines an {\it anti-isomorphism} of Lie algebras
\[(-)^*\colon \End_\C(V)\to \End_{\C}(V^*),\quad
[f_1,f_2]^*=-[f_1^*,f_2^*].
\]
The image of $\frak c(\mu)$ coincides with $\frak c(-\mu^*)$.
\end{rem}
\subsection{Spectrum of a flat $F$-manifold}\label{md2}%
Consider an analytic pointed flat $F$-manifold $(M,p)$. %
For any $\la\in\C$, we have a pair $(T_pM,\mu_p^\la)$ satisfying all the assumption of Section \ref{muop}. We can consequently introduce the Lie algebra $\frak c(\mu_p^\la)$, and the Lie group $\mc C(\mu_p^\la)$. 
\vskip2mm
Let $(V_1,\mu_1)$, $(V_2,\mu_2)$ be two pairs as in Section \ref{muop}. We say that they are {\it equivalent pairs} if there exist $\kappa\in\C$ and a $\C$-vector space isomorphism $f\colon V_1\to V_2$ such that $f\circ \mu_1=(\mu_2-\kappa\cdot{\rm id}_{V_2})\circ f$.
\vskip2mm
Given $\la_1,\la_2\in\C$, the pairs $(T_pM,\mu_p^{\la_1})$ and $(T_pM,\mu_p^{\la_2})$ are equivalent. Moreover, given $p_1,p_2\in M$,
the two pairs attached to the germs $(M,p_1)$ and $(M,p_2)$ are equivalent: using the connection $\nabla$, for any
path $\gm:[0, 1] \to M$ with $\gm(0) = p_1$ and $\gm(1) = p_2$, the parallel transport along $\gm$ provides an
isomorphism realizing the equivalence of the pairs at $p_1$ and $p_2$.
\vskip2mm
As a result, with any homogeneous flat $F$-manifold manifold $(M,\nabla,c,e,E)$ (not necessarily semisimple), we can canonically associate an equivalence %
class $[(V, \mu)]$ of pairs as above, which will be called the {\it spectrum} of $M$.
\vskip2mm
Fix a system of flat coordinates $\bm t=(t^1,\dots, t^n)$ diagonalizing $\mu^\la={\rm diag}(q_1-\la,\dots,q_n-\la)$. We can thus introduce the $\la$-independent matrix Lie algebras 
\bea
\frak c(\mu):=\left\{R\in\frak{gl}(n,\C)\colon R^\al_\bt=0\text{ unless }q_\al-q_\bt\in\Z_{>0}\right\},\\
\frak c(-\mu^*):=\left\{R\in\frak{gl}(n,\C)\colon R^\al_\bt=0\text{ unless }q_\al-q_\bt\in\Z_{<0}\right\},
\eea
which are canonically anti-isomorphic via transposition, see Remark \ref{caniso}. We also denote by $\mc C(\mu)$ and $\mc C(-\mu^*)$ the corresponding parabolic Lie groups.
\subsection{Solutions in Levelt normal forms and monodromy data at $z=0$}\label{md3} %
We now introduce some formal invariant of the given analytic flat $F$-manifold, by studying {\it Levelt normal forms} of solutions of the joint system of differential equations \eqref{jfs}.
\begin{thm}\label{solor}$\quad$\newline
\noindent $\rm (1)$ There exist $n\times n$-matrix valued functions $(G_p(\bm t))_{p\geq 1}$, analytic in $\bm t$, and a $\bm t$-independent matrix $R\in \frak c(-\mu^*)$, such that the matrix
\[\Xi(\bm t,z)=G(\bm t,z)z^{-\mu^\la} z^R,\quad G(\bm t,z)={\bf 1}+\sum_{p=1}^\infty G_p(\bm t)z^p,
\]is a (formal) solution of the joint system \eqref{jfs}.\newline 
\noindent $\rm (2)$ The series $G(\bm t, z)$ converges to an analytic function in $(\bm t,z)$. The matrix $\Xi(\bm t_o,z)$ is a fundamental system of solutions of the $\der_z$-equation of \eqref{jfs} for any fixed $\bm t_o$.  
\end{thm}
\proof
Consider $n$ functions $\tilde t^\al(\bm t,z)=\sum_{p=0}^\infty h^\al_p(\bm t)z^p$ such that $\nabla^z d\tilde t^\al=0$, and $\tilde t^\al(\bm t,0)=t^\al$. This translates in the following recursive equations for the coefficients $h^\al_p$:
\[h^{\al}_0(\bm t)=t^\al, \quad \der_\gm\der_\bt h^\al_{p+1}=c^\eps_{\gm\bt}\der_\eps h^\al_p,\quad p\geq 0.
\]Introduce the Jacobian matrix $J(\bm t, z):=(J^\al_\bt)_{\al,\bt}$, with $J^\al_\bt(\bm t,z):=\frac{\der \tilde t^\al}{\der t^\bt}$. Under the gauge transformation $\xi=J^T\tilde \xi$, the joint system \eqref{jfs} becomes
\begin{align}
\nonumber
\frac{\der\tilde \xi}{\der t^\al}&=\left(zJ\mc C_\al J^{-1}-\frac{\der J}{\der t^\al}J^{-1}\right)^T\tilde\xi=0,\\
\label{difeq1}
\frac{\der\tilde \xi}{\der z}&=\left[J\left(\mc U-\frac{1}{z}\mu^\la\right)J^{-1}-\frac{\der J}{\der z}J^{-1}\right]^T\tilde\xi\\
\nonumber
&=\left(-\frac{1}{z}(\mu^\la)^T+U_1^T+zU_2^T+z^2 U_3^T+\dots\right)\tilde\xi,
\end{align}
for suitable matrices $U_k$. From the compatibility $\der_\al\der_z=\der_z\der_\al$, it follows that the matrices $U_k$ are $\bm t$-independent. Up to a further gauge transformation $\tilde\xi\mapsto G(z)\tilde\xi$, of the form $G(z)=1+\sum_{k=1}^\infty G_kz^k$, the differential equation \eqref{difeq1} can be put in a \emph{normal form} 
\bean
\label{difeq2}
\frac{\der\tilde \xi}{\der z}=\left(-\frac{1}{z}\mu^\la+R_1+zR_2+z^2R_3+\dots\right),\\
\nonumber
 (R_k)^\al_\bt\neq 0\text{ only if }\mu^\la_\al-\mu^\la_\bt=-k,\quad k\geq 1.
\eean
Indeed, from the recursion relations
\[R_n=U_n^T+nG_n-[G_n,\mu^\la]+\sum_{k=1}^{n-1}(G_{n-k}U_k^T-R_kG_{n-k})
\]we determine the entries $(R_n)^\al_\bt$ for $(\mu^\la)_\al-(\mu^\la)_\bt=-n$, and $(G_n)^\al_\bt$ for $(\mu^\la)_\al-(\mu^\la)_\bt\neq -n$. We set $(G_n)_{\bt}^\al=0$ for $(\mu^\la)_\al-(\mu^\la)_\bt=-n$. See also \cite{Gan59}. A fundamental system of solutions of \eqref{difeq2} is given by $\tilde\xi(z)=z^{-\mu^\la} z^R$, where $R:=\sum_kR_k$. The proof of the convergence of the series $G(\bm t,z)$ is standard, the reader can consult e.g. \cite{wasow,Sib90}.
\endproof
\begin{cor}
The monodromy matrix $M_0$, defined by $\Xi(\bm t,e^{2\pi\sqrt{-1}}z)=\Xi(\bm t,z)M_0$, is independent of $\bm t$. We have  $M_0:=\exp(-2\pi\sqrt{-1}\mu^\la)\exp(2\pi\sqrt{-1}R)$.\qed
\end{cor}

Solutions $\Xi(\bm t,z)$ of the form above will be said to be in \emph{Levelt normal form}.

\begin{thm}
Assume that 
\bea
\Xi(\bm t,z)=G(\bm t,z)z^{-\mu^\la} z^R,\quad G(\bm t,z)={\bf 1}+\sum_{p=1}^\infty G_p(\bm t)z^p,\quad R\in\frak c(-\mu^*)\\
\widetilde\Xi(\bm t,z)=\widetilde G(\bm t,z)z^{-\mu^\la} z^{\widetilde R},\quad \widetilde G(\bm t,z)={\bf 1}+\sum_{p=1}^\infty \widetilde G_p(\bm t)z^p,\quad \widetilde R\in\frak c(-\mu^*),
\eea
are two solutions of the joint system \eqref{jfs}, in Levelt normal form. Then there exists a unique $C\in \mc C(-\mu^*)$ such that $\widetilde R=C^{-1} R C$, the function $p_C(z):=z^{-\mu^\la}z^RCz^{-\widetilde R}z^{\mu^\la}$ is polynomial in $z$, and $\widetilde G(\bm t,z)=p_C(z)G(\bm t,z)$. 
\end{thm}
\proof
By assumption there exist a unique invertible matrix $C\in M_n(\C)$ such that $\widetilde\Xi =\Xi C$. This implies that
\[G^{-1}\widetilde G=z^{-\mu^\la}z^RCz^{-\widetilde R}z^{\mu^\la}.
\]We deduce that the r.h.s.\,\,is a series in $z$ of the form $z^{-\mu^\la}z^RCz^{\widetilde R}z^{\mu^\la}={\bf 1}+H_1z+H_2z^2+\dots$. Actually, this sum is finite (i.e. a polynomial in $z$) because both $R$ and $\widetilde R$ are nilpotent. We can also re-write this identity as follows
\beq\label{mid1}
z^RCz^{-\widetilde R}=z^{\mu^\la}({\bf 1}+H_1z+H_2z^2+\dots)z^{-\mu^\la}.
\eeq
The l.h.s.\,\,is a polynomial in $\log z$, the r.h.s.\,\,contains only powers of $z$. Hence, both sides should actually be independent of $z$. The $(\al,\bt)$-entry of the r.h.s.\,\,equals 
\[\dl^\al_\bt+(H_1)^\al_\bt\ z^{(\mu^\la)_\al-(\mu^\la)_\bt+1}+(H_2)^\al_\bt\ z^{(\mu^\la)_\al-(\mu^\la)_\bt+2}+\dots,\]
which is $z$-independent iff $(H_k)^\al_\bt=0$ for $(\mu^\la)_\al-(\mu^\la)_\bt\neq -k$. Set $z=1$ in \eqref{mid1}: we have $C={\bf 1}+\sum_kH_k$. This shows that $C\in\mc C(-\mu^*)$.

We have just shown that both sides of \eqref{mid1} are $z$-independent and they equal $C$. The l.h.s.\,\,of equation \eqref{mid1} can also be written as $Cz^{C^{-1}RC}z^{-\widetilde R}$. Thus, $Cz^{C^{-1}RC}z^{-\widetilde R}=C$. It follows that $\widetilde R=C^{-1}RC$.
\endproof

\begin{defn}
We call \emph{monodromy data at $z=0$} of the flat $F$-manifold the datum  $(\la,\mu^\la,[R])$, where $[R]$ is the adjoint orbit, in the Lie algebra $\frak c(-\mu^*)$,  of the exponents of solutions of \eqref{jfs} in Levelt normal form. 
\end{defn}

\begin{rem}
The notion of spectrum of a flat $F$-manifold generalizes the corresponding notion  for Frobenius manifolds given in \cite{Dub99,CDG20}. In the Frobenius manifolds case, the group $\mc C(\mu)$ is replaced by its subgroup $\mc G(\eta,\mu)$ called {\it $(\eta,\mu)$-parabolic orthogonal group}: this is due to the fact that solutions in Levelt normal form satisfy a further $\eta$-orthogonality requirement described in Remark \ref{symderz}. Notice that in the Frobenius case we have $-\mu^*=\eta\mu\eta^{-1}$. See \cite[Section 2.1]{CDG20}.
\end{rem}

\subsection{Admissible germs and monodromy data at $z=\infty$}\label{md4}Let $(M,\nabla,c,e,E)$ be an analytic \emph{semisimple} homogenous flat $F$-manifold. Under semisimplicity assumption, the joint system of differential equations \eqref{jfs} is gauge equivalent to the joint system \eqref{jfss}. By studying this system, we are going to introduce another set of invariants of pointed germs of the flat $F$-manifold.
\vskip2mm
\noindent{\bf Semisimple, doubly resonant, and admissible germs. }An analytic pointed germ $(M,p)$ will be called
\begin{itemize}
\item {\it (tame/coalescing) semisimple} if the base point $p$ is (tame/coalescing) semisimple,
\item {\it doubly resonant} if $p$ is coalescing and $M$ is conformally resonant,
\item {\it strictly doubly resonant} if, for any arbitrarily fixed ordering $\bm u_o=(u_o^1,\dots,u_o^n)\in\C^n$ of the eigenvalues of $\mc U(p)$, we have
\[u^i_o=u^j_o,\qquad \dl_i-\dl_j\in\Z\setminus\{0\},\qquad \text{for some }i,j\in\{1,\dots,n\},\, i\neq j,
\]
\item {\it admissible} if it is semisimple but not strictly doubly resonant. %
\end{itemize}In this section, we will consider admissible pointed germs $(M,p)$. 
\begin{rem}
According to Definition \ref{tcp}, the specification ``tame/coalescing'' depends on the choice of the Euler vector field. In case $M$ is irreducible, it does not depend on such a  choice. This follows from Theorem \ref{eulerunique}.
\end{rem}
\vskip2mm
\noindent{\bf Formal solutions. }Let $(M,p)$ be an admissible germ. Fix an ordering $\bm u_o=(u_o^1,\dots,u_o^n)\in\C^n$ of the spectrum of the operator $\mc U(p)\colon T_pM\to T_pM$. Consider the $\der_z$-equation of the joint system \eqref{jfss} specialized at $\bm u=\bm u_o$. 
\begin{thm}\label{mrgAk}
There exist unique $n\times n$-matrices $(\mathring{A}_k)_{k\geq 1}$ such that the matrix 
\[\mrg X_{\rm for}(z)=\left({\bf 1}+\sum_{k=1}^\infty\frac{\mrg A_k}{z^k}\right)z^{\La}e^{z U_o},\quad {\La}:=\la\cdot{\bf 1}-\diag(\dl_1,\dots,\dl_n),\quad U_o={\rm diag}(u^1_o,\dots,u_o^n),
\]is a formal solution of the $\der_z$-equation of \eqref{jfss}, specialized at $\bm u=\bm u_o$. Moreover\footnote{Recall the notations introduced in Section \ref{secwpsietc}: for any matrix $A\in M_n(\C)$, we have the diagonal/off-diagonal decomposition $A=A'+A''$.}, we have $\mrg A_1''=\Gm(\bm u_o)^T$.
\end{thm}
\proof
The matrix $\mrg X_{\rm for}(z)$ is a solution of the $\der_z$-equation of \eqref{jfss} if and only if we have
\beq\label{recak}
(1-k)\mrg A_{k-1}+\mrg A_{k-1}\La=[U_o,\mrg A_k]-V(\bm u_o)^T\mrg A_{k-1},\quad k\geq 1,\quad \mrg A_0:={\bf 1}.
\eeq
We can compute recursively the matrices $\mrg A_k$. Let us start with $\mrg A_1$. \newline
\noindent$\bullet$ For $(i,j)$, with $i\neq j$, and so that $u_o^i\neq u_o^j$: from \eqref{recak} specialized at $k=1$, we deduce
\[(u_o^i-u_o^j)(\mrg A_1)^i_j=V^j_i=(u_o^i-u_o^j)\Gm(\bm u_o)^j_i\quad\Longrightarrow\quad (\mrg A_1)^i_j=\Gm(\bm u_o)^j_i.
\]
\noindent$\bullet$ For $(i,j)$, with $i\neq j$, and so that $u_o^i=u_o^j$: from \eqref{recak} specialized at $k=2$, we deduce
\begin{align*}
\nonumber
(\mrg A_1)^i_j&=\frac{1}{1-\dl_i+\dl_j}\sum_{\ell\neq i}V(\bm u_o)^\ell_i (\mrg A_1)^\ell_j
=\frac{1}{1-\dl_i+\dl_j}\sum_{\ell\neq i,j}(u^i_o-u^\ell_o)\Gm(\bm u_o)^j_\ell \Gm(\bm u_o)^\ell_i
=\Gm(\bm u_o)^j_i.
\end{align*}In the last equality, we used identity \eqref{midgm} specialized at $\bm u=\bm u_o$.\newline
\noindent$\bullet$ For the diagonal entries: from \eqref{recak} specialized at $k=2$, we deduce
\[(\mrg A_1)^i_i=\sum_{\ell\neq i}V(\bm u_o)^\ell_i (\mrg A_1)^\ell_i=\sum_{\ell\neq i}(u^i_o-u^\ell_o)\Gm(\bm u_o)^i_\ell\Gm(\bm u_o)^\ell_i.
\]
This completes the computation of $\mrg A_1$, and also proves that $\mrg A_1''=\Gm(\bm u_o)^T$. 

Assume now to have computed $\mrg A_1,\mrg A_2,\dots, \mrg A_{h-1}$. The matrix $\mrg A_h$ can be computed by repeating the same procedure. Namely, from \eqref{recak}, with $k=h$, one can compute the entries $(\mrg A_h)^i_j$ for $i\neq j$ such that $u_o^i\neq u_o^j$. From \eqref{recak}, with $k=h+1$, one can compute the remaining entries of $\mrg A_h$. %
\endproof

\begin{thm}\label{forsoldef}
Let $\Om\subseteq\C^n$ be a simply connected open neighborhood of $\bm u_o$. %
If $\Om$ is sufficiently small, then:\newline
\noindent{$\rm (1)$} For any $\bm u\in\Om$ there exist unique $n\times n$-matrices $({A}_k(\bm u))_{k\geq 1}$ such that the matrix 
\beq\label{exfsol}
X_{\rm for}(\bm u,z)=\left({\bf 1}+\sum_{k=1}^\infty\frac{A_k(\bm u)}{z^k}\right)z^{\La}e^{z U}\eeq
is a formal solution of the $\der_z$-equation of \eqref{jfss}. Moreover, we have $A_1(\bm u)''=\Gm(\bm u)^T$.\newline
\noindent{$\rm (2)$} We have $A_k(\bm u_o)=\mrg A_k$ for $k\geq 1$, and $X_{\rm for}(\bm u_o,z)=\mrg X(z)$.%
\end{thm}

\proof
Point (1) can be proved following the same computations as for Theorem \ref{mrgAk}. Point (2) follows by uniqueness. %
\endproof

\begin{rem}
From the computations above, it is clear that the coefficients $A_k$ are holomorphic at point $\bm u\in\Om$ such that $u^i\neq u^j$ for $i\neq j$. Below we will prove that the coefficients $A_k$ are actually holomorphic on the whole $\Om$.
\end{rem}

The identity $X_{\rm for}(\bm u,z e^{2\pi\sqrt{-1}})=X_{\rm for}(\bm u,z) e^{2\pi\sqrt{-1}\La}$ justifies the following terminology.
\begin{defn}
The matrix $\La={\rm diag}(\la-\dl_1,\dots, \la-\dl_n)$ is called {\it formal monodromy} matrix.
\end{defn}

\vskip2mm
\noindent{\bf Admissible directions $\tau$. }Let $q\in M$ be an arbitrary point, and fix an ordering $\bm u(q):=(u^1(q),\dots,u^n(q))\in\C^n$ of the eigenvalues of the operator $\mc U(q)\colon T_qM\to T_qM$. Denote by ${\rm Arg}(z)\in ]-\pi,\pi]$ the principal branch of the argument of the complex number $z$.  Set
\[
\mathscr S(q):=\left\{{\rm Arg}\left(-\sqrt{-1}(\overline{u^i(q)}-\overline{u^j(q)}\right)+2\pi k\colon k\in\Z,\ i,j\text{ are s.t. }u^i(q)\neq u^j(q)\right\}.
\]
Any element $\tau\in\R\setminus \mathscr S(q)$ will be called an \emph{admissible direction at $q$}.
\begin{rem}
The notion of admissibility only depends on the {\it set} $\{u^1(q),\dots, u^n(q)\}$.
\end{rem}
\noindent{\bf Asymptotic solutions. }Though the formal series defining $X_{\rm for}$ are typically divergent, $X_{\rm for}$ contains asymptotical information about genuine analytic solutions of the $\der_z$-equation of \eqref{jfss}.
\vskip2mm
Let $(M,p)$ be an admissible germ, $\Om$ as in Theorem \ref{forsoldef}, and $\tau$ an admissible direction at $p$. Consider a sufficiently small simply connected open neighborhood $\widetilde\Om\subseteq M$ of $p$ such that:
\begin{enumerate}
\item $\widetilde\Om\subseteq M_{ss}$, i.e. any point $q\in\widetilde\Om$ is semisimple,
\item a coherent choice of ordering $\bm u$ of eigenvalues of $\mc U$ is fixed on $\widetilde\Om$, so that $\bm u\colon \widetilde\Om\to \C^n$ defines a local system of canonical coordinates, with $\bm u(p)=\bm u_o$,
\item $\bm u(\widetilde\Om)\subseteq\Om$,
\item $\tau$ is admissible at any $q\in\widetilde\Om$.
\end{enumerate}

\begin{thm}\label{asymsol}
 If $\widetilde\Om$ is as above, then the following facts hold.
 \vskip1mm
 \noindent $\rm (1)$ For any $q\in\widetilde\Om$ there exist three fundamental systems of solutions $X_1,X_2,X_3$, of the $\der_z$-equation of \eqref{jfss} specialized at $\bm u=\bm u(q)$, uniquely determined by the asymptotics
 \beq\label{asym}
 X_i(\bm u,z)\sim X_{\rm for}(\bm u,z),\quad |z|\to+\infty,\quad \tau-(3-h)\pi<\arg z<\tau+(h-2)\pi,\quad h=1,2,3.
 \eeq
 \noindent$\rm (2)$ The functions $X_i$ are holomorphic w.r.t. $\bm u\in\bm u(\widetilde\Om)$, and the asymptotics \eqref{asym} holds true uniformly in $\bm u$.\newline
 \noindent $\rm (3)$ The functions $X_i(\bm u,z)$ are solutions of the joint system of differential equations \eqref{jfss}.\newline
  \noindent$\rm (4)$ The solutions $X_1$ and $X_3$ satisfy the identity $X_3(\bm u, ze^{2\pi\sqrt{-1}})=X_1(\bm u, z)e^{2\pi\sqrt{-1}\La}$, for $z\in\widehat{\C^*}$.
\end{thm}
\proof
Let us temporarily assume that $q$ is tame, i.e. $u^i(q)\neq u^j(q)$ for $i\neq j$.
For the proof of points (1) and (2), %
see e.g. \cite{BJL79,wasow}. %
Fix $h\in\{1,2,3\}$, and set $W_{i}(\bm u,z):=\der_iX_h(\bm u,z)-(zE_i-V_i)^TX_h(\bm u,z)$ for $i=1,\dots,n$ and $h=1,2,3$. A simple computation, invoking the identities \eqref{compmag} shows that $W_{i}(\bm u,z)$ is a solution of the $\der_z$-equation of \eqref{jfss}. Hence, there exist a matrix $C(\bm u)$ such that $W_{i}(\bm u,z)=X_h(\bm u,z)C(\bm u)$. Denote by $F(\bm u,z)={\bf 1}+z^{-1}A_1(\bm u)+O(z^{-2})$ the formal power series in \eqref{exfsol}. For $|z|\to+\infty$ in the sector
\[\mc V_{\tau,h}:=\left\{z\in\widehat{\mathbb C^*}\colon \tau-(3-h)\pi<\arg z<\tau+(h-2)\pi\right\},
\]the function $W_i(\bm u,z)$ has asymptotics
\begin{align*}W_i(\bm u,z)&\sim\der_i F(\bm u,z) z^{\La}e^{zU}+F(\bm u,z)zz^\La \overbrace{e^{zE_iu_i}}^{=E_ie^{zU}}-zE_iF(\bm u,z)z^{\La}e^{zU}+V_i^TF(\bm u,z)z^{\La}e^{zU}\\
&=\left(\der_i F(\bm u,z)+zF(\bm u,z)E_i-zE_iF(\bm u,z)+V_i^TF(\bm u,z)\right)z^{\La}e^{zU}.
\end{align*}
But we also have $W_i(\bm u,z)\sim F(\bm u,z)z^\La e^{zU}C(\bm u)$. As a consequence, we deduce
\beq
\label{czero}z^\La e^{zU}C(\bm u)e^{-zU}z^{-\La}=\text{formal power series in }\frac{1}{z}.
\eeq
For $j\neq k$, the sector $\mc V_{\tau,h}$ contains rays of points $z$ along which ${\rm Re}(z(u^j-u^k))>0$. Hence, necessarily, we deduce that the $(j,k)$-entry of $C(\bm u)$ vanishes, otherwise we would have a divergence on the l.h.s.\,\,of \eqref{czero}. So the matrix $C(\bm u)$ is diagonal, and 
\bea
C(\bm u)&=&z^\La e^{zU}C(\bm u)e^{-zU}z^{-\La}\\
&=&F(\bm u,z)^{-1}\left(\der_i F(\bm u,z)+zF(\bm u,z)E_i-zE_iF(\bm u,z)+V_i^TF(\bm u,z)\right)\\
&=&z(E_i-E_i)+(A_1E_i-E_iA_1+V_i^T)+O\left(\frac{1}{z}\right)=O\left(\frac{1}{z}\right),
\eea
where we used the identity $V_i^T=[E_i,\Gm^T]=[E_i,A_1]$. Hence $C(\bm u)=0$. This proves point (3) in the case $q$ is tame. 

The coefficients $A_k$ of Theorem \ref{forsoldef} are holomorphic at $\bm u$ such that $u^i\neq u^j$. Moreover, from the computations above, we deduce that
\beq\label{mrecak1}
[A_{k+1},E_i]=[A_1,E_i]A_k-\der_iA_k,\quad k\geq 1.
\eeq
This formula recursively determines the off-diagonal matrix $A_{k+1}''$ in terms of $A_1,\dots, A_k$. On the other hand, the diagonal entries of $A_{k+1}$ can be computed as in Theorem \ref{forsoldef}, so that
\beq\label{mrecak2}
(k+1)(A_{k+1})^i_i=\sum_{\ell\neq i}V^\ell_i (A_{k+1})^\ell_i=\sum_{\ell\neq i}(u^i-u^\ell)\Gm^i_\ell (A_{k+1})^\ell_i.
\eeq
Since $A_1''=\Gm^T$ is holomorphic also at coalescing points $\bm u$, an inductive argument shows that all the matrices $A_k$ are holomorphic at coalescing point, by using formulae \eqref{mrecak1} and \eqref{mrecak2}. See also \cite[Prop. 19.3]{CDG1}.
The system \eqref{jfss} is a completely integrable Pfaffian system with holomorphic coefficients on $\bm u(\widetilde\Om)$: the solutions $X_i(z,\bm u)$ can be $\bm u$-analytically continued as single-valued holomorphic functions on $\bm u(\widetilde\Om)$, see \cite[Cor. 19.1]{CDG1}. The assumptions of \cite[Th. 14.1]{CDG1} are thus satisfied, and (1),(2),(3) hold true also at coalescing points.\newline
Finally, notice that the two functions $X_1(\bm u,z)e^{2\pi\sqrt{-1}\La}$ and $X_3(\bm u,ze^{2\pi\sqrt{-1}})$ have the same asymptotics on the sector $\tau-2\pi< \arg z<\tau-\pi$. By uniqueness, it follows point (4).
\endproof

\begin{rem}\label{remvktau}For any $h=1,2,3$, the precise meaning of the uniform asymptotic relation \eqref{asym} is the following: for any compact $K\subseteq\bm u(\widetilde\Om)$, for any $\ell\in\N$, and for any unbounded closed subsector $\overline{\mathcal V}$ of $\mc V_{\tau,h}:=\left\{z\in\widehat{\mathbb C^*}\colon \tau-(3-h)\pi<\arg z<\tau+(h-2)\pi\right\}$, there exists a constant $C_{h,K,\ell,\overline{\mc V}}\in\R_{>0}$ such that
\[z\in\overline{\mc V}\setminus\{0\}\quad\Longrightarrow\quad
\sup_{u\in K}\left\| X_h(\bm u,z)e^{-z U}z^{-\La}-\left(\bm 1+\sum_{m=1}^{\ell-1}\frac{A_m(\bm u)}{z^m}\right)\right\|<\frac{C_{h,K,\ell,\overline{\mathcal V}}}{|z|^\ell}.
\]
\end{rem}
\vskip2mm
\noindent{\bf Stokes and central connection matrices. }Let $(M,p)$ be an admissible germ, and $\tau$ an admissible direction at $p$. Let $\Xi(\bm t,z)$ be a solution in Levelt form of the joint system \eqref{jfs}, and $X_h(\bm u,z)$, with $h=1,2,3$, be the solutions of the joint system \eqref{jfss} as in Theorem \ref{asymsol}. Let $\bm t_o=\bm t(p)$ and $\bm u_o=\bm u(p)$ the values of the flat and canonical coordinates at $p$, respectively. 
\vskip2mm
We define the {\it Stokes matrices} $\mrg S_1,\mrg S_2$ at $p$ to be the matrices defined by
\beq\label{defsto}
X_2(\bm u_o,z)=X_1(\bm u_o,z)\mrg S_1,\quad X_3(\bm u_o,z)=X_2(\bm u_o,z)\mrg S_2.
\eeq
We define the {\it central connection matrix} $\mrg C$ at $p$ to be the matrix defined by
\beq\label{defcc}
X_2(\bm u_o,z)=(\Psi(\bm u_o)^{-1})^T\cdot \Xi(\bm t_o,z)\cdot \mrg C.
\eeq

\begin{prop}\label{propsc1}We have
\begin{enumerate}
\item the matrices $\mrg S_1,\mrg S_2,\mrg C$ are invertible, with $\det \mrg S_1=\det \mrg S_2=1$,
\item $(\mrg S_1)_{ii}=(\mrg S_2)_{ii}=1$, 
\item if $i\neq j$, then $(\mrg S_1^{-1})_{ij}=0$ if ${\rm Re}\left(e^{\sqrt{-1}(\tau-\pi)}(u_o^i-u_o^j)\right)>0$, 
\item if $i\neq j$, then $(\mrg S_2)_{ij}=0$ if ${\rm Re}\left(e^{\sqrt{-1}\tau}(u_o^i-u_o^j)\right)>0$,
\item we have
\beq\label{c2} \mrg S_1^{-1}e^{2\pi\sqrt{-1}\La}\mrg S_2^{-1}=\mrg C^{-1}e^{-2\pi\sqrt{-1}\mu^\la}e^{2\pi\sqrt{-1}R}\mrg C.
\eeq
\end{enumerate}
\end{prop}
\proof
The proof of points (1)-(4) is standard, see \cite{wasow}. Point (5) follows from point (4) of Theorem \ref{asymsol}.
\endproof
\begin{prop}[\cite{CDG1,CG18}]\label{propsc2}
If $p$ is a coalescing point, define the partition $\{1,\dots,n\}$ $=\coprod_{r\in J}I_r$ such that for any $r\in J$ we have
$\{i,j\}\subseteq I_r\quad\text{if and only if}\quad u_o^i=u_o^j$. We then have the further vanishing condition 
\[\pushQED{\qed} 
(S_1)_{ij}=(S_1)_{ji}=(S_2)_{ij}=(S_2)_{ji}=0\quad\text{ if $i,j\in I_r$ for some $r\in J$.}\qedhere
\popQED\]
\end{prop}
For a $\mc D$-modules theoretical proof of Proposition \ref{propsc2} see the recent preprint \cite{Sab21}.
\vskip1mm
\begin{defn}
We call {\it monodromy data at $z=\infty$} of the admissible germ $(M,p)$, computed w.r.t.\,\,the admissible direction $\tau$, the 4-tuple of matrices $(\mrg S_1,\mrg S_2,\La, \mrg C)$. %
\end{defn}

\begin{rem}\label{mondatfrob}
In the case of Frobenius manifolds, with the standard choice $\la=\frac{d}{2}$, we have $\La=0$ and $\mrg S_1^{-1}=\mrg S_2^T$. This follow from the (anti-)self-adjointness properties of $\mc U$ and $\mu$. For detailed proofs see \cite[Th.\,2.42]{CDG20}. In the notations of {\it loc.\,cit.} we have $\mrg S_1=S$ and $\mrg S_2=S_-^{-1}$. Moreover, the Stokes matrices are uniquely determined by the metric, the central connection matrix, and the monodromy data at $z=0$:
\beq\label{StoFrob}
S=C^{-1}e^{-\pi\sqrt{-1}R}e^{-\pi\sqrt{-1}\mu}\eta^{-1}(C^{-1})^T,\qquad S_-=S^T=C^{-1}e^{\pi\sqrt{-1}R}e^{\pi\sqrt{-1}\mu}\eta^{-1}(C^{-1})^T.
\eeq
This is a direct consequence of the symmetries of the joint system \eqref{jfs}, see Remark \ref{symderz}.
\end{rem}
In the next paragraph we show that the monodromy data at $z=\infty$ define local invariants of the germ, i.e. that they are locally constant w.r.t.\,small perturbations of both the point $p$ and the admissible direction $\tau$.
\vskip2mm
\noindent{\bf Isomonodromicity Property. } Let $\widetilde\Om$ be an open neighborhood of $p$ as above. By Theorems \ref{solor} and \ref{asymsol}, if we let vary the point $q$ in $\widetilde\Om$, we have well defined solutions $\Xi(\bm t(q),z),X_i(\bm u(q),z)$, $i=1,2,3$, of the joint systems \eqref{jfs} and \eqref{jfss} respectively. We can thus introduce the Stokes and central connections matrices $(S_1,S_2,C)$ as functions of $q\in\widetilde\Om$ by the formulae
\bean
X_2(\bm u(q),z)&=&X_1(\bm u(q),z)\,S_1(\bm u(q)),\\ 
X_3(\bm u(q),z)&=&X_2(\bm u(q),z)\,S_2(\bm u(q)),\\ 
X_2(\bm u(q),z)&=&(\Psi(\bm u(q))^{-1})^T\cdot \Xi(\bm t(q),z)\cdot C(\bm u(q)).
\eean

\begin{thm}\label{isoth1}
The functions $S_1,S_2,C$ are constant on $\widetilde\Om$. In particular, we have $S_1(q)=\mrg S_1,S_2(q)=\mrg S_2, C(q)=\mrg C$, for all $q\in\widetilde\Om$.
\end{thm}

\proof
Let us prove the statement for $S_1$. We have
\bea\der_i S_1(\bm u)&=&\der_i\left[X_1(\bm u(q),z)^{-1}X_2(\bm u(q),z)\right]\\
&=&-X_1(\bm u(q),z)^{-1}\cdot\der_iX_1(\bm u(q),z)\cdot X_1(\bm u(q),z)^{-1}\cdot X_2(\bm u(q),z)\\
&&+X_1(\bm u(q),z)^{-1}\der_iX_2(\bm u(q),z)\\
&=&-X_1(\bm u(q),z)^{-1}\cdot(zE_i-V_i)^T\cdot X_2(\bm u(q),z)\\
&&+X_1(\bm u(q),z)^{-1}\cdot(zE_i-V_i)^T\cdot X_2(\bm u(q),z)=0.
\eea
The proof for $S_2,C$ is similar.
\endproof

\noindent{\bf Small perturbations of the admissible direction. }Let $(M,p)$ be an admissible germ, and $\tau$ be an admissible direction at $p$. 

\begin{thm}\label{isoth2}If $\tau'\in\R$ is such that $|\tau-\tau'|<\inf_{\phi\in\mathscr S(p)}|\tau-\phi|$, then $\tau'$ is admissible at $p$. Moreover, the monodromy data at $p$ computed w.r.t.\,\,$\tau$ and $\tau'$ are equal.
\end{thm}

\proof
The first claim is straightforward. Let us prove the second claim. There exist fundamental systems of solutions $X_i(\bm u(p),z)$ and $X_i'(\bm u(p),z)$, for $i=1,2,3$, such that
\bea
X_h(\bm u(p),z)\sim X_{\rm for}(\bm u(p),z),\quad |z|\to+\infty,\quad z\in\mc V_{\tau,h},\quad h=1,2,3,\\
X'_h(\bm u(p),z)\sim X_{\rm for}(\bm u(p),z),\quad |z|\to+\infty,\quad z\in\mc V_{\tau',h},\quad h=1,2,3,
\eea
by Theorem \ref{asymsol}, see also Remark \ref{remvktau}.
We prove that $X_h=X_h'$ for all $h=1,2,3$. 

Let $K_h$ be the connection matrix s.t.\,\,$X_h'=X_h\,K_h$. We have
\[z^\La e^{zU}K_he^{-zU}z^{-\La}\sim {\bm 1},\quad |z|\to +\infty,\quad z\in\mc V_{\tau,h}\cap \mc V_{\tau',h}.
\]By taking the $(j,k)$-entry, for any $\ell\in\N$ we have
\[(K_h)_{jk} e^{z(u^j(p)-u^k(p))}z^{\dl_k-\dl_j}= \dl_{jk}+O\left(|z|^{-\ell}\right),\quad |z|\to +\infty,\quad z\in\mc V_{\tau,h}\cap \mc V_{\tau',h}.
\]
Assume $j\neq k$. If $u^j(p)=u^k(p)$, then necessarily $(K_h)_{jk}=0$. If $u^j(p)\neq u^k(p)$, notice that in $\mc V_{\tau,h}\cap \mc V_{\tau',h}$ there are rays along which ${\rm Re}(z(u^j(p)-u^k(p)))$ is negative, and also rays along which it is positive. So, we necessarily have $(K_h)_{jk}=0$. This proves that $K_h$ is diagonal. It follows that $K_h=\bm 1$ for $h=1,2,3$.
\endproof

\subsection{Monodromy data for a formal admissible germ} In Sections \ref{md2}, \ref{md3}, \ref{md4}, the flat $F$-manifold structure on $M$ is assumed to be analytic. The notion of admissible germs and of their monodromy data can however be extended to the formal case.
\vskip1mm
Let $(H,\bm \Phi)$ be a semisimple formal $F$-manifold over $\C$, with Euler field $E$. Associated with it we have two joint systems of differential equations \eqref{jfs} and \eqref{jfss} whose coefficients are matrix-valued formal power series in the coordinates $\bm t$ and $\bm u$, respectively.
\vskip1mm
We will say that $(H,\bm \Phi)$ is
\begin{itemize} 
\item {\it doubly resonant}, if the origin is coalescing and the formal flat $F$-manifold is conformally resonant, 
\item {\it strictly doubly resonant}, if we have 
\[u^i_o=u^j_o,\qquad \dl_i-\dl_j\in\Z\setminus\{0\},\qquad\text{for some }i,j\in\{1,\dots,n\},\, i\neq j,
\]
\item {\it admissible} if it is semisimple but not strictly doubly resonant. %
\end{itemize}
The $\der_z$-equations of the joint systems \eqref{jfs} and \eqref{jfss} can be specialized at $\bm t=0$ and $\bm u=\bm u_o$, respectively. For these specialized systems of equations we can introduce a triple $(\la,\mu^{\la},[R])$ of monodromy data at $z=0$, and a 4-tuple $(\mrg S_1,\mrg S_2,\La,\mrg C)$ of monodromy data at $z=\infty$, exactly as in the case of an analytic germ $(M,p)$. 

The system $(\la,\mu^\la,[R],\mrg S_1,\mrg S_2,\La,\mrg C)$ will be referred to as the {\it monodromy data} of the formal structure $(H,\bm \Phi)$. {\it A priori}, Theorem \ref{isoth1} cannot be adapted to this formal picture, but Theorem \ref{isoth2} still holds true, and its proof works verbatim. 

In Section \ref{seconv}, we will prove that an admissible formal germ is actually convergent: it defines an analytic flat $F$-manifold, so that all the results of Sections \ref{md2}, \ref{md3}, \ref{md4} apply.

\section{Normalizations, and analytic continuation} \label{sec5}

\subsection{Choices of normalizations}\label{chnor}The monodromy data of an admissible germ $(M,p)$ are defined up to several non-canonical choices: 
\begin{enumerate}
\item the choice of $\la\in\C$,
\item the choice of a base point in the universal cover $\widehat{\C^*}$,
\item the choice of the solution $\Xi$ in Levelt normal form,
\item the choice of Lam\'e coefficients $(H_1,\dots, H_n)$,
\item the choice of ordering of canonical coordinates $(u^1(p),\dots, u^n(p))$,
\item the choice of an admissible direction $\tau\in\R\setminus\mathscr S(p)$.
\end{enumerate}
Different choices of normalizations affect the numerical values of the monodromy data. These transformations of the data can be described by actions of corresponding suitable groups:
\begin{enumerate}
\item the group $\C$,
\item the deck transformation group ${\rm Deck}(\widehat{\C^*})\cong \Z$,
\item the group $\mc C(-\mu^*)$,
\item the torus $(\C^*)^n$,
\item the symmetric group $\frak S_n$,
\item the braid group $\mc B_n$.
\end{enumerate}

We first describe actions (1)-(5), and postpone the description of action (6) to the next sections.
\vskip2mm
\noindent $\bullet$ {\it Action of $\C$}: the transformation $\lambda\mapsto \lambda'$ implies the following transformations of the monodromy data by translations
\[\mu^\la\mapsto \mu^{\la'}=\mu^\la+(\la-\la'){\bf 1},\quad\quad \La\mapsto \La-(\la-\la'){\bf 1}.
\]For irreducible flat $F$-manifolds, the choice of $\la$ is equivalent to the choice of an Euler vector field, see Theorem \ref{eulerunique}.
\vskip2mm
\noindent $\bullet$ {\it Action of ${\rm Deck}(\widehat{\C^*})\cong \Z$}: a different choice of the base point in $\widehat{\C^*}$ is equivalent to the choice of a different determination of the logarithm (i.e.\,\,of the argument $\arg z$). In particular, by changing $\log z\mapsto \log z+2\pi k\sqrt{-1}$ with $k\in\Z$, we have the transformations
\[\quad S_1\mapsto e^{-2\pi k \sqrt{-1}\La}\,S_1\,e^{2\pi k \sqrt{-1}\La}, \quad S_2\mapsto e^{-2\pi k \sqrt{-1}\La}\,S_2\,e^{2\pi k \sqrt{-1}\La},
\]
\[C\mapsto M_0^{-k}\,C\,e^{2\pi k \sqrt{-1}\La},\quad M_0:=e^{-2\pi\sqrt{-1}\mu^\la}e^{2\pi\sqrt{-1}R},\qquad k\in\Z.
\]
\vskip2mm
\noindent $\bullet$ {\it Action of $\mc C(-\mu^*)$}: for $A\in\mc C(-\mu^*)$, the change of solutions $\Xi\mapsto \Xi A$ implies the transformation of the central connection matrix
\[C\mapsto A^{-1}C.
\]
\noindent $\bullet$ {\it Action of $(\C^*)^n$}: for $(h_1,\dots,h_n)\in (\C^*)^n$, consider the transformation $(H_1,\dots, H_n)\mapsto (H_1h_1,\dots, H_nh_n)$. The monodromy data transform as follows
\[S_1\mapsto h^{-1}S_1h,\quad S_2\mapsto h^{-1}S_2h,\quad C\mapsto Ch,
\]
\[\text{where }h:={\rm diag}(h_1,\dots, h_n).
\]
\noindent $\bullet$ {\it Action of $\frak S_n$}: for $\sigma\in\frak S_n$, consider the permutation of canonical coordinates $(u^1,\dots, u^n)\mapsto (u^{\sigma(1)},\dots, u^{\sigma(n)})$. The monodromy data transform as follows
\[S_1\mapsto PS_1P^{-1},\quad S_2\mapsto PS_2P^{-1},\quad C\mapsto CP^{-1},\quad \La\mapsto P\La P^{-1},\]
\[\text{where }\quad P=(P_{ij})_{i,j},\quad P_{ij}:=\dl_{\sigma(i)j}.
\]
\begin{rem}
In the above discussion, we have fixed once for all a system of local flat coordinates $\bm t$ in which $\mu$ is diagonal. Different choices of systems of local coordinates as described in Remark \ref{remgroupG} affect the monodromy data. For example, if we have the transformation %
\[\bm t\mapsto \bm G\bm t+\bm c,\quad \bm G\in GL(n,\C),\,\,\bm G\text{ diagonal},\quad \bm c\in\C^n,
\]then the monodromy data change as follows:
\[%
R\mapsto \bm G^{-1} R\bm G,\quad C\mapsto \bm G^{-1}C.
\]This is not the most general form of the admissible transformations: if $\mu$ has not simple spectrum, there are also transformations with $\bm G$ not diagonal.
\end{rem}
\subsection{Triangular and lexicographical orders} If an admissible direction $\tau$ at $p$ is fixed, we will say that the canonical coordinates $(u^i(p))_{i=1}^n$ at $p$ are in {\it triangular order w.r.t.\,\,the admissible direction $\tau$}  if the Stokes matrix $S_1$ is upper triangular, and $S_2$ is lower triangular.
\vskip1mm
On the one hand, in general, triangular orders at $p$ are not unique. This happens for example if $p$ is a semisimple coalescing point. In such a case, we have $(S_1)_{ij}=(S_1)_{ji}=0$ if $u_i=u_j$ with $i\neq j$, by Proposition \ref{propsc2}. If $S_1$ is upper triangular, then so is $PS_1P^{-1}$ for $P$ corresponding to the transposition $i\leftrightarrow j$. Similarly, the lower triangular structure of $S_2$ is preserved.
\vskip1mm
On the other hand, we always have a distinguished triangular order, called {\it lexicographical w.r.t.\,\,$\tau$}. Introduce the following rays in the complex plane 
\[L_j:=\{u^j(p)+\rho e^{\sqrt{-1}(\frac{\pi}{2}-\tau)}\colon \rho\in\R_+\},\quad j=1,\dots,n.
\]The ray $L_j$ originates from the point $u^j(p)$, and it is oriented from $u^j(p)$ to $\infty$.
\vskip1mm
The canonical coordinates $(u^1(p),\dots, u^n(p))$ are in lexicographical order if $L_j$ is to the left of $L_k$ (w.r.t.\,\,the orientation above), for any $1\leq j<k\leq n$ such that $u^j(p)\neq u^k(p)$.
\vskip1mm
The lexicographical order is the unique triangular order at $p$ if the number of nonzero entries of $S_1$ or $S_2$ is maximal, i.e. $\frac{n(n-1)}{2}$.

\subsection{Braid group action on matrices }\label{secbrmon}Denote by $U_n$ and $L_n$ the groups of unipotent upper and lower triangular $n\times n$-matrices, and by $\frak t$ the Lie algebra of diagonal $n\times n$-matrices. 
\vskip2mm
The (abstract) {\it Artin braid group $\mc B_n$ with $n$-strings} is the group with $n-1$ generators $\bt_1,\dots, \bt_{n-1}$ satisfying the relations
\beq\label{cbr}
\bt_i\bt_j=\bt_j\bt_i,\quad \text{if }|i-j|>1,\quad \bt_i\bt_{i+1}\bt_i=\bt_{i+1}\bt_i\bt_{i+1}.
\eeq
Given $\bm g=(g_1,g_2,g_3)\in U_n\times L_n\times \frak t$, define $3(n-1)$ block-diagonal matrices $B^{(i)}_1(\bm g),B^{(i)}_2(\bm g),$ $B^{(i)}_3(\bm g)$, with $i=1,\dots, n-1$, as follows:
\[B^{(i)}_1(\bm g):={\bf 1}_{i-1}\oplus\left[\begin{array}{cc}
(g_1)_{i,i+1}&1\\
1&0
\end{array}\right]\oplus {\bf 1}_{n-i-1},
\]
\beq
\label{brb2}
B^{(i)}_2(\bm g):={\bf 1}_{i-1}\oplus\left[\begin{array}{cc}
0&1\\
1&(g_2)_{i+1,i}
\end{array}\right]\oplus {\bf 1}_{n-i-1},
\eeq
\[B^{(i)}_3(\bm g):={\bf 1}_{i-1}\oplus\left[\begin{array}{cc}
\hbar\cdot (g_1)_{i,i+1}&1\\
1&0
\end{array}\right]\oplus {\bf 1}_{n-i-1},
\]
where
\[\hbar:=e^{2\pi\sqrt{-1}[(g_3)_{i+1}-(g_3)_{i}]}.
\]
For any $\bt_i\in\mc B_n$, define the triple $\bm g^{\bt_i}\in U_n\times L_n\times\frak t$ by
\beq\label{actbr}
\bm g^{\bt_i}:=\left(B_1^{(i)}(\bm g)^{-1}g_1B_2^{(i)}(\bm g),\quad B_2^{(i)}(\bm g)^{-1}g_2B_3^{(i)}(\bm g),\quad P_i\,g_3\,P_i\right),
\eeq
where $P_i$ is the permutation matrix $i\leftrightarrow i+1$.
\begin{lem}
The braid group $\mc B_n$ acts on $U_n\times L_n\times \frak t$ by mapping $(\bt_i,\bm g)\mapsto \bm g^{\bt_i}$ for $i=1,\dots, n-1$.
\end{lem}
\proof
By a direct computation, one checks that $g^{\bt}=\rm id$ for any relator $\bt$ in \eqref{cbr}.
\endproof

\begin{example}
Let $n=3$, and
\[\bm g=\left(\left(
\begin{array}{ccc}
 1 & a & b \\
 0 & 1 & c \\
 0 & 0 & 1 \\
\end{array}
\right),\left(
\begin{array}{ccc}
 1 & 0 & 0 \\
 \alpha  & 1 & 0 \\
 \beta  & \gamma  & 1 \\
\end{array}
\right),\left(
\begin{array}{ccc}
 d_1 & 0 & 0 \\
 0 & d_2 & 0 \\
 0 & 0 & d_3 \\
\end{array}
\right)\right).
\]

We have
\[\bm g^{\bt_1}=\left(\left(
\begin{array}{ccc}
 1 & \alpha  & c \\
 0 & 1 & b-a c \\
 0 & 0 & 1 \\
\end{array}
\right),\left(
\begin{array}{ccc}
 1 & 0 & 0 \\
 a e^{2 \sqrt{-1} (d_2-d_1) \pi } & 1 & 0 \\
 a\bt e^{2  \sqrt{-1} (d_2-d_1) \pi } +\gamma  & \beta  & 1 \\
\end{array}
\right),\left(
\begin{array}{ccc}
 d_2 & 0 & 0 \\
 0 & d_1 & 0 \\
 0 & 0 & d_3 \\
\end{array}
\right)\right),
\]
\[ 
\bm g^{\bt_2}=\left(
\left(
\begin{array}{ccc}
 1 & b & a+b \gamma  \\
 0 & 1 & \gamma  \\
 0 & 0 & 1 \\
\end{array}
\right),\left(
\begin{array}{ccc}
 1 & 0 & 0 \\
 \beta -\alpha  \gamma  & 1 & 0 \\
 \alpha  & c e^{2 \sqrt{-1} (d_3-d_2) \pi } & 1 \\
\end{array}
\right),\left(
\begin{array}{ccc}
 d_1 & 0 & 0 \\
 0 & d_3 & 0 \\
 0 & 0 & d_2 \\
\end{array}
\right)\right).
\]
If $\bt=(\bt_1\bt_2)^3$, the triple ${\bm g}^\bt=(g'_1,g'_2,g'_3)$ equals
\[g_1'=\left(
\begin{array}{ccc}
 1 & a e^{2 \sqrt{-1} (d_2-d_1) \pi } & b e^{2 \sqrt{-1} (d_3-d_1) \pi } \\
 0 & 1 & c e^{2 \sqrt{-1} (d_3-d_2) \pi } \\
 0 & 0 & 1 \\
\end{array}
\right)=e^{-2\pi\sqrt{-1}g_3}\,g_1\,e^{2\pi\sqrt{-1}g_3},\]
\[\quad g_2'=\left(
\begin{array}{ccc}
 1 & 0 & 0 \\
 e^{2 \sqrt{-1} (d_1-d_2) \pi } \alpha  & 1 & 0 \\
 e^{2 \sqrt{-1} (d_1-d_3) \pi } \beta  & e^{2 \sqrt{-1} (d_2-d_3) \pi } \gamma  & 1 \\
\end{array}
\right)=e^{-2\pi\sqrt{-1}g_3}\,g_2\,e^{2\pi\sqrt{-1}g_3},\]
\[
g'_3=
\left(
\begin{array}{ccc}
 d_1 & 0 & 0 \\
 0 & d_2 & 0 \\
 0 & 0 & d_3 \\
\end{array}
\right)=g_3.
\]
\end{example}

\subsection{Braid mutations of monodromy data}\label{brmut}Let $M$ be an analytic homogeneous semisimple flat $F$-manifold, and denote by $M'$ the open set of tame semisimple points $p\in M$, i.e.\,\,at which the spectrum of the operator $\mc U(p)\colon T_pM\to T_pM$ is simple.
\vskip1mm
Consider the following two different settings:
\begin{enumerate}
\item[(I)] Assume $g \colon[0,1]\to M'$ to be a continuous path such that $g([0,1])$ is contained in a simply connected open set, on which a coherent choice of normalizations (1)-(5) can be done. Assume also that
\begin{itemize}
\item $\tau$ is admissible at both $g(0)$ and $g(1)$,
\item there exists $\bar t\in[0,1]$ such that $\tau$ is not admissible at $g(\bar t)$.
\end{itemize}

\item[(II)] Assume $p\in M'$ is a semisimple point, and fix some choice of normalizations (1)-(5). Let $\tau_0,\tau_1\in\R\setminus\mathscr S(p)$, and $\tau\colon[0,1]\to M$ to be a continuous map such that 
\begin{itemize}
\item $\tau(0)=\tau_0$ and $\tau(1)=\tau_1$,
\item there exists $\bar t\in[0,1]$ such that $\tau(\bar t)$ is not admissible at $p$.
\end{itemize}
\end{enumerate}
In both cases (I) and (II), for each $t\in\{0,1\}$, we can introduce a set $\mc M_t$ of monodromy data. \vskip2mm
{\bf Problem: }In both settings (I) and (II), how to describe the transformation $\mc M_0\mapsto \mc M_1$?
\vskip2mm
The matrices $\mu^\la,R$ will not depend on $t$, due to the results of Section \ref{md3}. Hence, we need to describe how the matrices $(S_1,S_2,C,\La)$ will transform. In this section, we prove that this is described by an action of the braid group, which on the triple $(S_1,S_2,\La)$ reduces to \eqref{actbr}. %

\begin{rem}
Pictures (I) and (II) are ``dual'' to each other. 
In (I), we have a fixed $\tau\in\R$ and a variable set $\mathscr S(g(t))$ of non-admissible directions such that $\tau\in \mathscr S(g(0))\cap\mathscr S(g(1))$.

In (II), we have a fixed set $\mathscr S(p)\subseteq\R$ of non-admissible directions and a continuous map $\tau\colon [0,1]\to \R$ with $\tau(0),\tau(1)\in\R\setminus\mathscr S(p)$.  In both cases, we have to face a {\it wall-crossing phenomenon}: the fixed (resp. variable) point $\tau$ is not admissible for some values of the time parameter.
\end{rem}

Given $\bm u\in\C^n$ introduce a family of {\it Stokes rays} in the universal cover $\widehat{\C^*}$: for any pair $(i,j)$ such that $u^i\neq u^j$ set 
\[\tau_{ij}(\bm u):=\frac{3\pi}{2}-{\rm Arg}(u^i-u^j),\qquad R_{ij}^{(k)}(\bm u):=\{z\in\widehat{\C^*}\colon \arg z=\tau_{ij}(\bm u)+2\pi k\},\quad k\in\Z.
\]Also, for any $\tau\in \R$ introduce the {\it admissible ray}
\[\ell_\tau:=\{z\in\widehat{\C^*}\colon \arg z=\tau\}.
\]
Both Stokes and admissible rays are equipped with the natural orientation, from 0 to $\infty$. Any continuous transformations of $\bm u$ and $\tau$ induce continuous rotations of the Stokes and admissible rays. In the case of settings (I) and (II), the oriented ray crosses some of the Stokes rays during the transformation. We will call {\it elementary} any such transformation of rays, along which $\ell_\tau$ crosses one Stokes ray $R_{ij}^{(k)}$ only.
\vskip1mm
Let us focus on the picture (I). 
Fix $\bm u_o\in\C^n\setminus\Dl$ with components in $\tau$-lexicographical order. Consider a continuous map $b_i\colon [0,1]\to (\C^n\setminus\Dl)$, with $i=1,\dots,n-1$, such that: 
\begin{enumerate}
\item $b_i(0)=\bm u_o$,
\item $b_i(t)^h=b_i(0)^h$ for all $h\neq i,i+1$, 
\item $b_i(t)^i$ counter-clockwise rotates by one half turn w.r.t.\,\,$b_i(t)^{i+1}$ in the plane $\C$, 
\item $b_i(0)^i=b_i(1)^{i+1}$ and $b_i(1)^i=b_i(0)^{i+1}$. 
\end{enumerate}The map $b_i$ can be seen as a loop on ${\rm Conf}_n(\C):=(\C^n\setminus\Dl)/\frak S_n$,  the configuration space of $n$ pairwise distinct points in $\C$. The space $\cfsp$ is aspherical (i.e. $\pi_k(\cfsp)=0$ for $k\geq 2$), and its fundamental group is isomorphic to the braid group $\mc B_n$, see \cite{KT08}. Consider the homotopy classes $[b_i]$ in $\pi_1({\rm Conf}(\C^n),\{u^i(0)\})\cong \mc B_n$. It is easily seen that
\[[b_i]*[b_j]=[b_j]*[b_i],\quad |i-j|>1,\qquad [b_i]*[b_{i+1}]*[b_i]=[b_{i+1}]*[b_i]*[b_{i+1}],
\]where $*$ denotes the concatenation of loops. We identify $[b_i]$ with the elementary braid $\bt_i$.
\vskip1mm
In the case of picture (II), any of the maps $b_i$'s can be seen as a map with target $M$, by working in a local chart with canonical coordinates in $\tau$-lexicographical order. It is an elementary transformation: one of the Stokes rays $R_{i,i+1}^{(k)}$ clockwise crosses the ray $\ell_\tau$. 
\vskip2mm
In summary, elementary transformations of type (I) can be identified with elements of $\mc B_n$. Dually, by exchanging orientations (counter-clockwise $\leftrightarrow$ clockwise), we can identify $\bt_i$ with the type (II) transformation defined by a counter-clockwise rotation of $\ell_\tau$ across one of the Stokes rays $R_{i,i+1}^{(k)}$.
\vskip2mm
Let $(S_1,S_2,\La,C)$ be the 4-tuple of Stokes, formal monodromy, and central connection matrices computed 
\begin{itemize}
\item w.r.t.\,\,the point $g(0)$, in case (I);
\item w.r.t.\,\,the line $\tau(0)$, in case (II);
\end{itemize}
In both cases (I) and (II), the monodromy data are always computed w.r.t.\,\,the lexicographical order of canonical coordinates, so that $(S_1,S_2,\La)\in U_n\times L_n\times \frak t$.

\begin{thm}\label{thmbmmd}
Along the elementary transformation $\bt_i$, with $i=1,\dots, n-1$, the monodromy data transform as follows:
\[(S_1,S_2,\La)\mapsto (S_1,S_2,\La)^{\bt_i},\quad C\mapsto C B^{-1},
\]where 
\[B=B_2^{(i)}(S_1,S_2,\La)={\bf 1}_{i-1}\oplus\left[\begin{array}{cc}
0&1\\
1&(S_2)_{i+1,i}
\end{array}\right]\oplus {\bf 1}_{n-i-1}.
\]
\end{thm}
Cf. equations \eqref{brb2} and \eqref{actbr}. 
\proof
Whatever is the case under consideration, (I) or (II), let us consider the initial ``frozen'' configuration of Stokes and admissible rays, $R_{ij}^{(h)}$ and $\ell_\tau$. 
\vskip1mm
Label the Stokes rays as follows: let $R^{(1)}$ be the first Stokes ray on the left of $\ell_\tau$, $R^{(0)}$ the first Stokes ray on the right of $\ell_\tau$, and extend the numeration $R^{(k)}$, with $k\in\Z$, so that the label $k$ increases in counter-clockwise order.
\vskip1mm
Let $m$ be the number of Stokes rays in any sector of $\widehat{\C^*}$ defined by $$(2j-1)\pi<|\arg z-\tau|<2\pi j,\qquad j\in\Z.$$
The number $m$ also equals the number of Stokes rays in the sectors 
$$2\pi j<|\arg z-\tau|<(2j+1)\pi,\qquad j\in\Z.$$
For generic initial points $\bm u_o\in\C^n\setminus\Dl$, we have $m=\frac{n(n-1)}{2}$, but some Stokes rays may coincide\footnote{This happens for example if there are three indices $(i,j,k)$ such that $u_o^i,u_o^j,u_o^k\in\C$ are collinear, or there are four indices $(i,j,k,l)$ such that $u_o^i,u_o^j$ and $u_o^k,u_o^l$ define two parallel lines.}.
\vskip1mm
Define $\Pi^{(k)}$, with $k\in\Z$, to be the sector in $\widehat{\C^*}$ from the ray $R^{(k-1)}$ to the ray $R^{(m+k)}$. For each $k\in\Z$, there exists a unique solution $X^{(k)}(\bm u_o,z)$ of the $\der_z$-equation of \eqref{jfss}, specialized at $\bm u=\bm u_o$, such that 
\[X^{(k)}(\bm u_o,z)\sim X_{\rm for}(\bm u_o,z),\qquad |z|\to+\infty,\quad z\in\Pi^{(k)}.
\]Introduce invertible matrices $K_k$, called {\it Stokes factors}, such that
\[X^{(k+1)}(\bm u_o,z)=X^{(k)}(\bm u_o,z)\,K_k,\quad k\in\Z.
\]The matrices $K_k$ have the following structure: all diagonal entries are 1, and the entry $(K_k)_{ab}$ is non-zero only if $R^{(m+k)}$ is one of the rays $R^{(h)}_{ab}$ with $h\in\Z$, see \cite{BJL79}.

Recall that the Stokes matrices $S_1,S_2$ are defined in terms of solutions $X_1,X_2,X_3$, see equation \eqref{defsto}. We have
\[X_1(\bm u_o, z)\equiv X^{(1-m)}(\bm u_o,z),\quad X_2(\bm u_o, z)\equiv X^{(1)}(\bm u_o,z),\quad X_3(\bm u_o,z)\equiv X^{(1+m)}(\bm u_o,z).
\]Hence, we deduce
\begin{align*}X_{2}(\bm u_o,z)&=X^{(0)}(\bm u_o,z)\,K_0=X^{(-1)}(\bm u_o,z)\,K_{-1}K_0=\dots=\underbrace{X^{(1-m)}(\bm u_o,z)}_{X_1(\bm u_o,z)}\,{K_{1-m}\dots K_{-1}K_0}.
\end{align*}
From equation \eqref{defsto}, we deduce
\[S_1=K_{1-m}\dots K_{-1}K_0.
\]
Analogously, we have 
\[S_2=K_1\dots K_m.
\]Up to now we have considered a ``static'' picture, at the initial time $t=0$ of the transformation $\bt_i$. By letting the time parameter $t$ vary, the Stokes rays and/or the ray $\ell_\tau$ rotate. In particular, immediately before the collision of Stokes and oriented rays, we have $R_{i,i+1}^{(h)}=R^{(1)}$ for a suitable $h\in\Z$. After the collision we have $R_{i,i+1}^{(h)}\equiv R^{(0)}$. Hence, after the transformation $\bt_i$, we have the following transformation of Stokes matrices
\bea
S_1&\mapsto& S_1'=K_{-m}\dots K_{0}K_{1}= K_{1-m}^{-1}\,S_1\,K_1,\\
S_2&\mapsto& S_2'=K_2\dots K_mK_{m+1}=K_1^{-1}\,S_2\,K_{1+m}.
\eea
Similarly, the central connection matrix transforms as follows
\[C\mapsto C'= C\,K_1^{-1}
\]
The only non-zero off-diagonal entries of $K_1,K_{1-m},K_{1+m}$ are
\[(K_{1-m})_{i,i+1}=(S_1)_{i,i+1},\qquad (K_1)_{i+1,i}=(S_2)_{i+1,i},\]
\[(K_{1+m})_{i,i+1}=[e^{-2\pi\sqrt{-1}\La}S_1e^{2\pi\sqrt{-1}\La}]_{i,i+1}.
\]The last identity follows from point (4) of Theorem \ref{asymsol}. 
Finally, we also need to recover the lexicographical order, which is lost after the transformation $\bt_i$. By applying the permutation $i\leftrightarrow i+1$, we complete the proof.
\endproof

\begin{rem}
Theorem \ref{thmbmmd} generalizes the braid group action in Dubrovin's analytic theory of Frobenius manifolds, see \cite[Th.\,4.8]{Dub99}. The action of the braid group on $U_n\times L_n\times \frak t$ is just the simplest case of a more general picture described in \cite{Boa01,Boa02}. The starting point is the observation that the ``monodromy manifold''\footnote{In \cite{Boa01,Boa02}, meromorphic connections on a closed disk $D\subseteq \C$, with an irregular singularity only, are studied. We have a connections on $\C^*$ with two singularities, and we need a further piece of information for a global description of the monodromy: the central connection matrix.} $U_n\times L_n\times\frak t$ is isomorphic to the dual Poisson-Lie group $G^*$ of $G=GL(n,\C)$. In \cite[Section 2]{Boa02} P.\,Boalch generalized the notion of Stokes multipliers and isomonodromic deformations for general connected complex reductive groups $G$. It was also proved that $G^*$ can be identified with the space of meromorphic connections on principal $G$-bundles over the disc. In such a case one has an action of $\pi_1(\frak t_{\rm reg})$ on $G^*$, where $\frak t_{\rm reg}$ is the regular subset of a Cartan subalgebra $\frak t\subseteq\frak g$. This induces an action on $G^*$ of the full braid group $\pi_1(\frak t_{\rm reg}/W)$, with $W$ the Weyl group.
Such an action coincides with the De Concini-Kac-Procesi action of $\pi_1(\frak t_{\rm reg}/W)$ on $G^*$, obtained in \cite{DKP92} as classical limit of the {\it quantum Weyl group action} on the corresponding quantum group, due to Lusztig \cite{Lus90}, Kirillov and Reshetikhin \cite{KR90}, Soibelman \cite{Soi90}.
\end{rem}

\noindent{\bf Action of the center $Z(\mc B_n)$. }Consider the shift of the admissible direction $\tau\mapsto \tau+2\pi$. We have the following facts:
\begin{enumerate}
\item In the generic case (i.e. for canonical coordinates in general position), the number of Stokes rays in any sector of $\C^*$ of width $2\pi$ equals $n(n-1)$. An elementary braid acts whenever the line $\ell_\tau$ crosses a Stokes ray. So, in total, we expect that a complete rotation of $\ell_\tau$ correspond to the product of $n(n-1)$ elementary braids.
\item The effect of the shift $\tau\mapsto \tau+2\pi$ on the monodromy data can be identified with a transformation of different nature, namely a different choice of normalization (2). This consists in a different choice of the branch of the logarithm. From this it follows that the braid corresponding to $\tau\mapsto \tau+2\pi$ must commute with any other braids.
\end{enumerate}
From point (2), we deduce that the braid corresponding to $\tau\mapsto \tau+2\pi$ is an element of the center $$Z(\mc B_n)\cong \Z\cong {\rm Deck}(\widehat{\C^*}).$$
The center $Z(\mc B_n)$ is the cyclic group generated by the braid $\bt=(\bt_1\dots\bt_{n-1})^n$. From point (1), and the action of ${\rm Deck}(\widehat{\C^*})$, we deduce the following result.
\begin{prop}\label{propcen}
The braid corresponding to the shift $\tau\mapsto \tau+2\pi$ of the admissible direction is the generator $\bt$ of the center $Z(\mc B_n)$. It acts as follows:
\[(S_1,S_2,\La)^\bt=\left(e^{-2\pi\sqrt{-1}\La} S_1 e^{2\pi\sqrt{-1}\La},\quad e^{-2\pi\sqrt{-1}\La} S_2 e^{2\pi\sqrt{-1}\La},\quad\La\right)
\]
\[\pushQED{\qed} 
C\mapsto M_0^{-1} C e^{2\pi\sqrt{-1}\La}.\qedhere
\popQED
\]
\end{prop}
This proposition extends a computation (for $n=3$) of the Example of Section \ref{secbrmon}.
\begin{cor}
The generator $\bt=(\bt_1\dots\bt_{n-1})^n$ of the center $Z(\mc B_n)$ acts on $U_n\times L_n\times\frak t$ as follows:
\beq\label{cenact}
\bm g^\bt=\left(e^{-2\pi\sqrt{-1}g_3}\,g_1\,e^{2\pi\sqrt{-1}g_3},\quad e^{-2\pi\sqrt{-1}g_3}\,g_2\,e^{2\pi\sqrt{-1}g_3},\quad g_3\right).
\eeq
\end{cor}
\proof
Let $\bm g'$ be the r.h.s.\,of \eqref{cenact}. Assume there exists $\tilde\bt\in\mc B_n$ such that $\bm g^{\tilde\bt}=\bm g'$. A simple computation shows that $B^{(i)}_j(\bm g)e^{2\pi\sqrt{-1}g_3}=B^{(i)}_j(\bm g')$ for $i=1,\dots, n-1$ and $j=1,2,3$. Thus, we have $(\bm g^{\tilde\bt})^{\bt_i}=(\bm g^{\bt_i})^{\tilde\bt}$ for any $i=1,\dots,n-1$ and any $\bm g\in U_n\times L_n\times\frak t$. Denote by $\frak S_{U_n\times L_n\times\frak t}$ the group of bijections of $U_n\times L_n\times\frak t$, and by $\iota\colon\mc B_n\to\frak S_{U_n\times L_n\times\frak t}$ the group morphism defining the action. We have $\iota(\tilde\bt)\in\iota\left(Z(\mc B_n)\right)$. Hence there exists $k\in\Z$ such that $\bm g^{\tilde\bt}=\bm g^{\bt^k}$ for any $\bm g\in U_n\times L_n\times\frak t$. %
We deduce $k=1$, by Proposition \ref{propcen}.
\endproof

\subsection{Analytic continuation of the flat $F$-structure} There is a more global point of view from which one can reinterpret the results of the previous sections. It both makes transparent the appearance of a braid group action on the monodromy data, and clarifies the ``duality'' of settings (I) and (II) of the previous section. Moreover, it also describes the analytic continuation of the flat $F$-manifold structure.
\vskip2mm
\noindent{\bf Admissibility sets $\An,\Ano,\Anz,\Anu$. }Consider the configuration space $\cfsp:=(\C^n\setminus\Dl)/\frak S_n$ of $n$ points in the plane, together with the following covering space:
\begin{itemize}
\item the ordered $\frak S_n$-covering $\cfspo:=\C^n\setminus\Dl$,
\item the {covering} $\cfspz$ associated with the center $Z(\pi_1(\cfsp))\cong\Z$,
\item and the universal cover  $\cfspu$. %
\end{itemize}
The fundamental group $\pi_1(\cfspo)$ is isomorphic to the group $\mc P_n$ of pure braids. The deck transformation group of the covering $\cfspo\to\cfsp$ is isomorphic to $\frak S_n\cong \mc B_n/\mc P_n$. 
\vskip1mm
The fundamental group of the space $\cfspz$ is isomorphic to $Z(\pi_1(\cfsp))\cong\Z$. The deck transformation group of the covering $\cfspz\to\cfsp$ is isomorphic to $\mc B_n/Z(\mc B_n)\cong M_n(\R^2)$, the mapping class group of the $n$-punctured plane, see \cite{Bir74}.
\vskip1mm
Notice that $Z(\mc B_n)= Z(\mc P_n)\cong \Z$. The coverings spaces above fit into the chain 
\[\cfspu\to \cfspz\to\cfspo\to\cfsp.
\]

\begin{rem}
The space $\cfspu$ has been described for the first time  by S.\,Kaliman \cite{Kal75,Kal77,Kal93}: in {\it loc.\,\,cit.\,\,}it is proved that it is isomorphic to $\C^2\times \mc T(0,n+1)$, where $\mc T(0,n+1)$ denotes the Teichm\"uller space of the Riemann sphere with $n+1$ punctures. The space $\mc T(0,n+1)$ is homeomorphic to $\R^{2n-4}$, and it is biholomorphic to a holomorphically convex Bergmann domain in $\C^{n-2}$. For further details see the interesting paper \cite{Lin04}.
\end{rem}
\vskip1mm
Given $\bm p=\{p_1,\dots, p_n\}\in\cfsp$, define the set $\mathscr S(p)\subseteq \R$ by
\[\mathscr S(\bm p):=\{{\rm Arg}[-\sqrt{-1}(\overline{p_i}-\overline{p_j})]+2k\pi\colon\quad k\in\Z\},\quad {\rm Arg}(z)\in]-\pi,\pi].
\]Any number $\tau\in\R\setminus \mathscr S(p)$ is said to be an {\it admissible direction} at $\bm p$.
\vskip1mm
Introduce the smooth $(2n+1)$-dimensional real manifolds
\[\mc X_n:=\cfsp\times \R,\quad \mc X_n^O:=\cfspo\times \R,\quad \Xnz:=\cfspz\times \R,\quad \widehat{\mc X_n}:=\cfspu\times \R,
\]Define the {\it admissibility open subset} $\An\subseteq \Xn$ by
\[\An:=\{(\bm p,\tau)\in\Xn\colon \tau\in \R\setminus \mathscr S(\bm p)\}.
\]Analogously, define the open subsets $\Ano\subseteq \Xno$, $\Anz\subseteq \Xnz$, and $\Anu\subseteq \Xnu$ as the pre-images of $\An$ along the projections $\Xnu\to\Xnz\to\Xno\to\Xn$. We have the following commutative diagram
\[
\xymatrix @R=1pc {
\Anu\ar@{^{(}->}[d]\ar[r]&\Anz\ar@{^{(}->}[d]\ar[r]&\Ano\ar@{^{(}->}[d]\ar[r]&\An\ar@{^{(}->}[d]&\\
\Xnu\ar[d]\ar[r]%
&\Xnz\ar[d]\ar[r]%
&\Xno\ar[d]\ar[r]%
&\Xn\ar[d]\ar[r]%
&\R\\
\cfspu\ar[r]&\cfspz\ar[r]&\cfspo\ar[r]&\cfsp&\\
}
\]
\vskip2mm
\noindent{\bf Homotopy groups of $\An,\Ano,\Anz,\Anu$. }Consider the subspace $\An'$ of admissibility set $\An$ defined by$$\An':=\{(\bm p,0)\in\Xn\colon 0 \text{ is admissible at $\bm p$}\}.$$
\begin{lem}\label{lem1.5}
The subspace $\An'$ is a strong deformation retract of $\An$.
\end{lem}
\proof
Let $(\bm p,\tau)\in\An$. If $\bm p=\{p_1,\dots, p_n\}$, denote by $p^{b}:=\frac{1}{n}\sum_{j=1}^np_j$ the barycenter of the configuration.
Let $F\colon[0,1]\times \An\to\An$ be defined by a rotation w.r.t.\,the barycenter
\[F(t,\bm p,\tau):=\left(\left\{e^{\sqrt{-1}t\tau}(p_j-p^b)+p^b\colon\quad j=1,\dots,n\right\},\quad \tau(1-t)\right).
\]For all $(\bm p,\tau)\in \An$, $(\bm a,0)\in\An'$ and $t\in[0,1]$, we have $F(0,\bm p,\tau)=(\bm p,\tau)$,  $F(1,\bm p,\tau)\in\An'$, and  $F(t,\bm a,0)=(\bm a,0)$.
\endproof

\begin{lem}\label{lem2}
The space $\An'$ is contractible.
\end{lem}
\proof
We show that the point $(\{1,\dots,n\},0)$ is a strong deformation retract of $\An'$. Given $(\bm p,0)\in\An'$, with $\bm p=\{p_1,\dots, p_n\}$, without loss of generality we may assume that the $p_j$'s are labelled in $0$-lexicographical order. Consider the continuous map $F\colon [0,1]\times \An'\to \An'$ defined by
\begin{empheq}[left={F(t,\bm p)=\empheqlbrace\,}]{align*}
      & \left(\left\{(1-2t)p_j+2t{\rm Re}(p_j)\colon j=1,\dots, n\right\},\quad 0\right),\qquad 0\leq t\leq \frac{1}{2},\\
      & \left(\left\{2t{\rm Re}(p_j)+(2t-1)j\colon j=1,\dots, n\right\},\quad 0\right),\qquad \frac{1}{2}\leq t\leq 1.
\end{empheq}
The map $F$ defines a strong deformation retraction of $\An'$ onto $(\{1,\dots, n\}, 0)$.
\endproof

\begin{thm}\label{hgan}
We have 
\bea
&\pi_i(\An)=0,\quad &i=0,1,2,\dots,\\
&\pi_0(\Ano)=\frak S_n,\quad \pi_i(\Ano)=0,\quad &i=1,2,3,\dots,\\
&\pi_0(\Anz)=M_n(\R^2),\quad \pi_i(\Anz)=0,\quad &i=1,2,3,\dots,\\
&\pi_0(\Anu)=\mc B_n,\quad \pi_i(\Anu)=0,\quad &i=1,2,3,\dots,
\eea
where the homotopy groups are based at an arbitrary point.
\end{thm}
\proof
All homotopy groups of $\An$ vanish, since $\An$ is contractible by Lemmata \ref{lem1.5}, \ref{lem2}. Since $\An$ is simply connected, we have the homeomorphisms $\Ano\cong \An\times \frak S_n$, $\Anz\cong\An\times M_n(\R^2)$, and $\Anu\cong \An\times \mc B_n$. Here $\frak S_n$, $\Z$, and $\mc B_n$ are equipped with the discrete topology. The claim follows.
\endproof
\noindent{\bf Relative homotopy groups. }
Given a triple $(A,B,c)$ of pointed topological spaces, with $c\in B\subseteq A$, the relative homotopy group $\pi_k(A,B,c)$, with $k\geq 1$, is the set\footnote{The group structure is well defined only for $k\geq 2$.} of homotopy classes of continuous maps $f\colon (\mathbb D^k,\mathbb S^{k-1},s_0)\to (A,B,c)$, with $s_0\in\mathbb S^{k-1}$. In particular, the set $\pi_1(A,B,c)$ is the set of homotopy classes of paths $f\colon [0,1]\to A$ such that $f(0)\in B$, $f(1)=c$.
\vskip1mm
Fix a point $\bm x\in\An$, and three points $\bm x^o\in\Ano,\bm x^z\in\Anz, \hat{\bm x}\in\Anu$ over it.
\begin{thm}
We have
\bea
\label{rhg1}
&\pi_1(\Xn,\An,{\bm x})\cong\pi_1(\Xno,\Ano,{\bm x}^o) \cong\pi_1(\Xnz,\Anz,{\bm x}^z)\cong\pi_1(\Xnu,\Anu,\hat{\bm x}) \cong\mc B_n,&\\
\label{rhg2}
&\pi_k(\Xn,\An,{\bm x})\cong\pi_k(\Xno,\Ano,{\bm x}^o) \cong\pi_k(\Xnz,\Anz,{\bm x}^z)\cong\pi_k(\Xnu,\Anu,\hat{\bm x}) \cong 0,\quad& k\geq 2.
\eea
\end{thm}

\proof
With the morphisms of triples of topological spaces 
\beq
\label{mxa}
(\Xnu,\Anu,\hat{\bm x})\to (\Xnz,\Anz,{\bm x}^z)\to (\Xno,\Ano,\bm x^o)\to (\Xn,\An,\bm x),
\eeq
we can associate the following commutative diagram of relative homotopy groups
\[
\xymatrix @C=1pc{
\dots\ar[r]&\pi_1(\Anu,\hat{\bm x})\ar[d]\ar[r]&\pi_1(\Xnu,\hat{\bm x})\ar[d]\ar[r]&\pi_1(\Xnu,\Anu,\hat{\bm x})\ar[d]^{\al_1}\ar[r]&\pi_0(\Anu,\hat{\bm x})\ar[d]\ar[r]&\pi_0(\Xnu,\hat{\bm x})\ar[d]\\
\dots\ar[r]&\pi_1(\Anz,{\bm x}^z)\ar[d]\ar[r]&\pi_1(\Xnz,{\bm x}^z)\ar[d]\ar[r]&\pi_1(\Xnz,\Anz,{\bm x}^z)\ar[d]^{\al_2}\ar[r]&\pi_0(\Anz,{\bm x}^z)\ar[d]\ar[r]&\pi_0(\Xnz,{\bm x}^z)\ar[d]\\
\dots\ar[r]&\pi_1(\Ano,{\bm x}^o)\ar[d]\ar[r]&\pi_1(\Xno,{\bm x}^o)\ar[d]\ar[r]&\pi_1(\Xno,\Ano,{\bm x}^o)\ar[d]^{\al_3}\ar[r]&\pi_0(\Ano,{\bm x}^o)\ar[d]\ar[r]&\pi_0(\Xno,{\bm x}^o)\ar[d]\\
\dots\ar[r]&\pi_1(\An,{\bm x})\ar[r]&\pi_1(\Xn,{\bm x})\ar[r]&\pi_1(\Xn,\An,{\bm x})\ar[r]&\pi_0(\An,{\bm x})\ar[r]&\pi_0(\Xn,{\bm x})
}
\]
The rows are the long exact relative homotopy sequences for each triples, and columns are the maps induced by \eqref{mxa}. The maps $\al_1,\al_2,\al_3$ are bijections: this follows from the unique lifting property of paths for coverings.
The claim then follows from Theorem \ref{hgan}.
\endproof

\vskip2mm
\noindent{\bf Monodromy data as functions on $\mc A_M^Z$. }Let $M$ be a flat $F$-manifold with Euler field $E$, and denote by $M'$ the set of points $p\in M$ at which the spectrum ${\rm spec}(E\circ_p)$ is simple. We have the local biholomorphism
\beq\label{upsmap}
\ups\colon M'\to\cfsp,\quad p\mapsto {\rm spec}(E\circ_p).
\eeq
Consider the pulled-back fiber bundles on $M'$
\[\mc X_M:=\ups^*\Xn,\quad \mc X_M^O:=\ups^*\Xno,\quad \mc X_M^Z:=\ups^*\Xnz,\quad \widehat{\mc X_M}:=\ups^*\Xnu,
\]together with their open subsets 
\[\mc A_M:=(\ups^*)^{-1}\An,\quad \mc A_M^O:=(\ups^*)^{-1}\Ano, \quad\mc A_M^Z:=(\ups^*)^{-1}\Anz, \quad \widehat{\mc A_M}:=(\ups^*)^{-1}\Anu.
\]
Given $p_o\in M'$, the monodromy data of $(M',p_o)$ are well-defined after fixing the choice of normalizations (1)-(6) of Section \ref{chnor}.
\vskip1mm
The choice of (5) only, i.e. an ordering of canonical coordinates at $p_o$, is equivalent to the choice of a point of $\cfspo$ %
over $v(p_o)={\rm spec}(E\circ_{p_o})$. 
\vskip1mm
The choice of (6) only, i.e. an admissible direction at $p_o$, is equivalent to the choice of a point of $\mc A_M$ over $p_o$.
\vskip1mm
The choice of both (5) and (6) is equivalent to the choice of a point of $\mc A_M^O$ over $p_o$.
\vskip1mm
If choices of (1),(2),(3),(4) are fixed, however, the 4-tuple $(S_1,S_2,\La,C)$ is not well-defined as a single-valued function on $\mc A_M^O$. Indeed, if $(p_o,\bm u_o,\tau)$ is a fixed point of $\mc A_M^O$, for any $k\in\Z$ there exist paths $\gm_k\colon [0,1]\to \mc A_M^O$ such that $\gm(0)=(p_o,\bm u_o,\tau)$ and $\gm_k(1)=(p_o,\bm u_o,\tau+2\pi k)$. Namely $\gm_k$ are lifts of loops in the center $Z(\pi_1(\cfspo))$. 
\vskip1mm
Thus, the joint choice of (2),(5),(6) is equivalent to the choice of a point of $\mc A_M^Z$ over $p_o$. 
\vskip1mm
Theorems \ref{isoth1} and \ref{isoth2} can be reformulated as follows.
\begin{thm}
Fix a choice of normalizations (1),(3),(4), and a point in $\mc A_M^Z$. The monodromy data $(S_1,S_2,\La,C)$ are locally constant functions on $\mc A_M^Z$.\qed
\end{thm}
In total we have ${\rm card}\,\pi_0(\mc A_M^Z)$ possible values of the monodromy data at $z=\infty$. Different values at different connected components of $\mc A_M^Z$, are labelled by paths in $\pi_1(\mc X_M^Z, \mc A_M^Z)$. The map \eqref{upsmap} induces a morphism in homotopy
\[\ups_*\colon \pi_1(\mc X_M^Z, \mc A_M^Z)\to \pi_1(\Xnz,\Anz)\cong \mc B_n.
\]
The paths of settings (I) and (II) of Section \ref{brmut} are representatives of homotopy classes in $\pi_1(\mc X_M^O, \mc A_M^O)\cong \pi_1(\mc X_M^Z, \mc A_M^Z)$. More precisely, consider the double fibrations
\[
\xymatrix {
&\mc A_M^O\ar[d]\ar@/_2pc/[ddl]_{p_1}\ar@/^2pc/[ddr]^{p_2}&\\
&\mc X_M^O\ar[dl]^{\rho_1}\ar[dr]_{\rho_2}&\\
M'&&\R
}
\]Paths of setting (I) represent classes in $\pi_1(\rho_2^{-1}(\tau),p_2^{-1}(\tau))$ for fixed $\tau\in\R$. 
\vskip1mm
\noindent Paths of setting (II) represent classes in $\pi_1(\rho_1^{-1}(p),p_1^{-1}(p))$ for fixed $p\in M'$.
\vskip1mm
\noindent Thus the ``duality'' mentioned in Section \ref{brmut} reflects the underlying double fibrations above. In both cases (I) and (II), we have induced paths in $\pi_1(\mc X_M^O, \mc A_M^O)$.
\vskip1mm
\noindent Following the terminology of \cite{CDG20}, given $\tau\in\R$ we call $\tau$-\emph{chamber} of $M$ any connected component of the open set $p_1(p_2^{-1}(\tau))$. 
\vskip2mm
\noindent{\bf Analytic continuation. }Let $(M,p_o)$ be the germ of a semisimple analytic flat $F$-manifold with Euler vector field. Assume that $p_o$ is tame semisimple. The whole flat $F$-structure can be analytically continued. The picture described in this section gives an insight on this continuation procedure. 
\vskip1mm
By shrinking $M$, we can assume that the germ is defined on a simply connected set, sufficiently small so that \eqref{upsmap} is an embedding. We can thus identify $M'$ with $\ups(M')\subseteq \cfsp$. By fixing a point $\hat{\bm u}_o\in\cfspu$ over $\ups(p_o)$, we have an open embedding of $(M',p_o)\cong (\ups(M'),\ups(p_o))\subseteq (\cfspu,\hat{\bm u}_o)$. In this way, one finds a maximal tame analytic continuation of the initial germ. Notice that the coefficients of the joint system of differential equations \eqref{jfss} continue to meromorphic functions of $\bm u\in\cfspu$: this is the {\it Painlev\'e property} of the solution $V^\la(\bm u)$ of the isomonodromic differential equations \eqref{compmag}. By fixing choices of normalizations (1),(3),(4), the monodromy data $(S_1,S_2,\La,C)$ of the system \eqref{jfss} can be seen as locally constant functions on the space $\Anu$. This space has countably many connected components in bijection with the braid group $\mc B_n$. All possible values of $(S_1,S_2,\La,C)$ are given by the action of the braid group $\mc B_n$ of Theorem \ref{thmbmmd}.

\section{Riemann-Hilbert-Birkhoff inverse problem for semisimple flat $F$-manifolds}\label{RHBsec}
\subsection{RHB problem $\mc P[${\boldmath $u$}$,\tau,\frak M]$ and the Malgrange-Sabbah Theorem}$\quad$\newline 

\noindent{\bf Admissible data. }Denote by ${\rm Arg}(z)\in ]-\pi,\pi]$ the principal branch of the argument of the complex number $z$.  Let $\bm u\in\C^n$, and set
\[
\mathscr S(\bm u):=\left\{{\rm Arg}\left(-\sqrt{-1}(\overline{u^i}-\overline{u^j}\right)+2\pi k\colon k\in\Z,\ i,j\text{ are s.t. }u^i\neq u^j\right\}.
\]
Any element $\tau\in\R\setminus \mathscr S(\bm u)$ will be said to be \emph{admissible at $\bm u$}.

\begin{defn}\label{datum}
Let $\bm u\in\C^n$ and $\tau$ admissible at $\bm u$. A $(\bm u,\tau)$-{\it admissible datum} is a $6$-tuple $\frak M:=(B,D,L,S_1,S_2,C)$ of matrices in $M_n(\C)$ such that:
\begin{enumerate} 
\item the matrix $B$ is diagonal, i.e. $B=B'$,
\item $D={\rm diag}(D_1,\dots,D_n)$ is a diagonal matrix of integers,
\item we have $L_{ij}=0$ if $D_i-D_j\in\Z_{<0}$,
\item we have
\beq
\label{vad1}
\tr B=\tr D+\tr L.
\eeq
\item the matrices $S_1,S_2,C$ are invertible, with $\det S_1=\det S_2=1$,
\item $(S_1)_{ii}=(S_2)_{ii}=1$, 
\item if $i\neq j$, then $(S_1^{-1})_{ij}=0$ if ${\rm Re}\left(e^{\sqrt{-1}(\tau-\pi)}(u^i-u^j)\right)>0$, 
\item if $i\neq j$, then $(S_2)_{ij}=0$ if ${\rm Re}\left(e^{\sqrt{-1}\tau}(u^i-u^j)\right)>0$,
\item we have
\beq\label{c2bis} S_1^{-1}e^{2\pi\sqrt{-1}B}S_2^{-1}=C^{-1}e^{2\pi\sqrt{-1}L}C.
\eeq
\end{enumerate}
If $\bm u\in\Dl$, define the partition $\{1,\dots,n\}=\coprod_{r\in J}I_r$ such that for any $r\in J$ we have
$\{i,j\}\subseteq I_r$ if and only if $u^i=u^j$. We then require the further vanishing condition 
\begin{enumerate}
\item[(10)] $(S_1^{-1})_{ij}=(S_2)_{ij}=0$ if $i,j\in I_r$ for some $r\in J$.
\end{enumerate}
\end{defn}

\begin{lem}
Let $\bm u_o\in \C^n$ and $\tau$ admissible at $\bm u_o$. If $\frak M$ is $(\bm u_o,\tau)$-admissible, then there exists a sufficiently small neighborhood $\mc V$ of $\bm u_o$ such 
\begin{enumerate}
\item $\tau$ is admissible at $\bm u$, for all $\bm u\in\mc V$,
\item $\frak M$ is $(\bm u,\tau)$-admissible for all $\bm u\in\mc V$. \qed
\end{enumerate}
\end{lem}

\begin{figure}
\centering
\def\svgscale{.6}
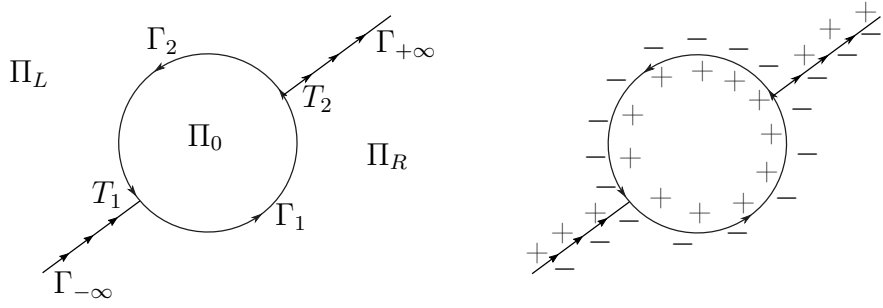
\caption{Contour $\Gm$, paths $\Gm_{\pm\infty},\Gm_1,\Gm_2$, domains $\Pi_0,\Pi_L,\Pi_R$, and $\pm$ sides of $\Gm$.}
\label{disegno1}
\end{figure}

Let $\bm u\in\C^n$ and $\tau$ admissible at $\bm u$. Consider the complex $z$-plane with a branch cut from $0$ to $\infty$:
\[\tau-\pi<\arg z<\tau+\pi.
\]
Let $r>0$ and denote by $\Gm=\Gm(\tau,r)$ the union of the following oriented paths, see Figure \ref{disegno1}:
\begin{enumerate}
\item the half-line $\Gm_{-\infty}$ defined by $\arg z=\tau\pm\pi$, $|z|>r$, originating from $\infty$;
\item the half-line $\Gm_{+\infty}$ defined by $\arg z=\tau$, $|z|>r$, ending to $\infty$;
\item the half-circle $\Gm_1$ defined by $\tau-\pi<\arg z<\tau$, $|z|=r$, counterclockwise oriented;
\item the half-circle $\Gm_2$ defined by $\tau<\arg z<\tau+\pi$, $|z|=r$, counterclockwise oriented.
\end{enumerate}
The orientations uniquely define the + and - side for each path $\Gm_{\pm\infty},\Gm_1,\Gm_2$. For $z\in\Gm_{-\infty}$ we use the symbol $z_\pm$ if $\arg z=\tau\pm\pi$. Set $\Pi_0,\Pi_L,\Pi_R$ to be the components of complement $\C\setminus\Gm$, and $T_1, T_2$ to be the two nodes of $\Gm$, as in Figure \ref{disegno1}.
\vskip2mm
Let $\frak M:=(B,D,L,S_1,S_2,C)$ be a $(\bm u,\tau)$-admissible datum. Define two functions
\[Q(-;\bm u), \mc H(-;\bm u)\colon \Gm\to GL(n,\C),
\] by
\[Q(z;\bm u):=U(\bm u)z+B\log z,\quad U(\bm u):={\rm diag}(u^1,\dots, u^n),
\]
 \begin{empheq}[left=\mc H(z;\bm u){:=}\empheqlbrace]{align*} 
       &e^{Q(z_-;\bm u)}S_1^{-1} e^{-Q(z_-;\bm u)},\quad\text{ along }\Gm_{-\infty},\\
       &e^{Q(z;\bm u)}S_2 e^{-Q(z;\bm u)},\quad\text{ along }\Gm_{+\infty},\\
       &e^{Q(z;\bm u)} C^{-1}z^{-L}z^{-D},\quad\text{ along }\Gm_{1},\\
       &e^{Q(z;\bm u)} S_2^{-1}C^{-1}z^{-L}z^{-D},\quad\text{ along }\Gm_{2}.
    \end{empheq}

\begin{prob}[{Problem $\mc P[\bm u,\tau,\frak M]$}] Find an analytic function $G\colon \C\setminus\Gm\to M_n(\C)$ such that 
\begin{enumerate}
\item $G|_{\Pi_\nu}$ extends continuously to $\overline{\Pi_\nu}$ for $\nu=0,L,R$;
\item the non-tangential limits $G_\pm\colon \Gm\to M_n(\C)$ of $G$ from the - and + sides of $\Gm$ exist, and are continuous;
\item they are related by
\[G_+(z)=G_-(z)\mc H(z;\bm u);
\]
\item $G(z)$ tends to the identity matrix $I$ as $z\to\infty$.
\end{enumerate}
\end{prob}

\begin{thm}[{\cite[Section 3]{Cot20b}}]\label{teoimp}
Let $\bm u_o\in\C^n$. Assume that the pair $(\tau,\frak M)$ is admissible at each point of a sufficiently small open neighborhood $\mc V$ of $\bm u_o$. If $\mc P[\bm u_o,\tau,\frak M]$ is solvable, there exists an analytic set $\Theta\subseteq \mc V\setminus \{\bm u_o\}$ such that $\mc P[\bm u,\tau,\frak M]$ is solvable for all $\bm u\in\mc V\setminus\Theta$. Moreover, the solution $G(z;\bm u)$ is unique and holomorphic w.r.t. $\bm u\in \mc V\setminus\Theta$.\qed
\end{thm}

\begin{rem}
In \cite{Cot20b}, we showed that Theorem \ref{teoimp} is essentially equivalent to a (weaker) extension, %
due to C.\,Sabbah \cite[Th.\,4.9]{Sab18}, of a previous result of B.\,Malgrange \cite{Mal83b}. For this reason, we refer to Theorem \ref{teoimp} as Malgrange-Sabbah Theorem. The original result of Malgrange concerns the case $\bm u_o\in\C^n\setminus\Dl$. The result of Sabbah concerns the case $\bm u_o\in\Dl$.
\end{rem}

\subsection{Construction of semisimple flat $F$-manifolds via a RHB inverse problem}\label{RHBcon}Let $\bm u_o\in\C^n$, $\tau$ be an admissible direction at $\bm u_o$, and $\frak M=(B,D,L,S_1,S_2,C)$ be a $(\bm u_o,\tau)$-admissible datum.  Assume that the RHB boundary value problem $\mc P[\bm u_o,\tau,\frak M]$ is solvable. Let $\mc V$ and $\Theta$ as in Malgrange-Sabbah Theorem \ref{teoimp}: the problem $\mc P[\bm u,\tau,\frak M]$ is well-defined, solvable and with unique solution $G(z,\bm u)$, holomorphic w.r.t.\,\,$\bm u\in\mc V\setminus\Theta$. Consider the asymptotic expansions of $G(z,\bm u)$ for $z\to 0$ and $z\to\infty$:
\bea
G(z,\bm u)={\bf 1}+{z^{-1}}F_1(\bm u)+O\left(z^{-2}\right),\quad z\to\infty,\quad z\in\Pi_{L/R},\\
G(z,\bm u)=G_0(\bm u)+z G_1(\bm u)+z^2 G_2(\bm u)+O(z^3),\quad z\to 0,
\eea
with coefficients $F_1,G_i$'s holomorphic w.r.t.\,\,$\bm u$. Define the functions
\bea
&X_{L/R}(z,\bm u):=G(z,\bm u)z^Bz^{zU},\quad &z\in\Pi_{L/R},\\
&X_0(z,\bm u):=G(z,\bm u)z^Dz^L,\quad &z\in \Pi_0.
\eea
\begin{lem}
The functions $X_0(z,\bm u),X_{L/R}(z,\bm u)$ are solutions of the joint system of differential equations
\bean
\label{rhjs1}
&\frac{\der}{\der u^i}X=\left(zE_i-V_i(\bm u)^T\right)X,\quad &V_i(\bm u):=[F_1(\bm u)^T,E_i]\equiv -\left(\frac{\der G_0}{\der u^i}\cdot G_0^{-1}\right)^T,\\
\label{rhjs2}
&\frac{\der}{\der z}X=\left(U-\frac{1}{z}V(\bm u)^T\right)X,\quad &V(\bm u):=[F_1(\bm u)^T,U]-B\equiv -(G_0^T)^{-1}(D+L')G_0^T.
\eean
\end{lem}
\proof
We have $X_L(z,\bm u)=X_R(z,\bm u)S_2$ and $X_R(z,\bm u)=X_0(z,\bm u)C$. It follows that $\der_zX_0\cdot X_0^{-1}=\der_zX_L\cdot X_L^{-1}=\der_zX_R\cdot X_R^{-1}$, and the resulting function $f(z,\bm u):=\der_zX_i(z,\bm u)X_i(z,\bm u)^{-1}$, with $i=0,R,L$, is analytic with respect to $z\in\C^*$. We have
\bea
\frac{\der X_{L/R} }{\der z}\cdot X_{L/R}^{-1}&=&\der_zG\cdot G^{-1}+\frac{1}{z}GBG^{-1}+GU G^{-1}\\
&=&U+\frac{1}{z}\left([F_1(\bm u),\La(\bm u)]+B\right)+O\left(\frac{1}{z^2}\right),\quad z\to\infty,\\
\frac{\der X_{0} }{\der z}\cdot X_{0}^{-1}&=&\der_zG\cdot G^{-1}+\frac{1}{z}\left(GDG^{-1}+Gz^DLz^{-D}G^{-1}\right)\\
&=& \frac{1}{z}G_0(\bm u)\left(D+L'\right)G_0(\bm u)^{-1}+O(1),\quad z\to0.
\eea
The last equality follows from the identity $z^DLz^{-D}=L'+O(z)$, which is deduced from the admissibility condition (3) on the matrix $L$. The matrices $S_1,S_2,C$ are constant with respect to both $\bm u$ and $z$: we deduce that the r.h.s.\,\,of the two equalities above are equal. This implies that $X_{L/R}$ and $X_0$ are solutions of the differential equation \eqref{rhjs2}.
Similarly, we have
\bea
\frac{\der X_{L/R} }{\der u^i}\cdot X_{L/R}^{-1}&=&\frac{\der G}{\der u^i}\cdot G^{-1}+zGE_iG^{-1}=zE_i+[F_1,E_i]+O\left(\frac{1}{z}\right),\\
\frac{\der X_{0} }{\der u^i}\cdot X_{0}^{-1}&=&\frac{\der G}{\der u^i}\cdot G^{-1}=\frac{\der G_0}{\der u^i}\cdot G^{-1}_0+O(z),
\eea
where $(E_i)_{ab}=\dl_{ai}\dl_{bi}$. The matrices $S_1,S_2,C$ being constant, we deduce that the r.h.s.\,\,of the two equalities above are equal. Hence $X_{L/R}$ and $X_0$ are solutions of the differential systems \eqref{rhjs1}.
\endproof
\begin{rem}
Notice that $V(\bm u)$ is diagonalizable with diagonal Jordan form $-D-L'$. In particular, the eigenvalues do not depend on $\bm u$.
\end{rem}

\begin{lem}\label{def1}
The off-diagonal entries $(F_1'')^T$ satisfy the Darboux-Egoroff equations \eqref{DE1}, \eqref{DE2}, and the homogeneity conditions
\beq
\label{homof1}
\sum_{k=1}^nu^k\der_k F_1(\bm u)^i_j=(b_i-b_j-1)F_1(\bm u)^i_j,\quad B={\rm diag}(b_1,\dots,b_n).
\eeq
\end{lem}
\proof
The compatibility condition $\der_i\der_j=\der_j\der_i$ for the joint system \eqref{rhjs1},\eqref{rhjs2} reads
\[[E_j,\der_iF_1^T]-[E_i,\der_jF_1^T]+[[E_i,F_1^T],[E_j,F_1^T]]=0.
\]This coincides with equations \eqref{DE1} and \eqref{DE2}. Let $\kappa\in\C^*$, and set $$h(z):=z^Dz^LC\kappa^{-B}C^{-1}\kappa^{-L}z^{-L}\kappa^{-D}z^{-D}.$$ The piecewise analytic function $\widetilde G\colon (\Pi_0\cup\Pi_L\cup\Pi_R)\times(\kappa\mc V\setminus\kappa\Theta)\to\C$ defined by
\bea
&\widetilde G(z;\bm u):=\kappa^{-B}G(\kappa z;\kappa^{-1} \bm u)h(z)^{-1},\quad &z\in\Pi_0,\\
& \quad \widetilde G(z;\bm u):= \kappa^{-B}G(\kappa z;\kappa^{-1} \bm u)\kappa^B,\quad &z\in\Pi_{L/R},
\eea %
solves the same RHB problem $\mc P[\bm u,\tau,\frak M]$ as $G$. By uniqueness of solution we have $\widetilde G=G$. This implies that $F_1(\kappa^{-1}\bm u)=\kappa\cdot\kappa^{-B}F_1(\bm u)\kappa^B$, and \eqref{homof1} follows.
\endproof

Define the off-diagonal matrix $\Gm(\bm u)$ by $\Gm(\bm u)^i_j:=F_1(\bm u)^j_i$. 
\begin{cor}
For any fixed $\bm H_o\in(\C^*)^n$, there exists a unique $\bm H(\bm u)=(H_1(\bm u),\dots, H_n(\bm u))^T$, analytic in $\mc V\setminus\Theta$, satisfying
\beq\label{pfaff}
\der_j H_i=\Gm^i_jH_j,\quad i\neq j,\quad\quad \der_iH_i=-\sum_{k\neq i}\Gm^i_kH_k,
\eeq and such that $\bm H(\bm u_o)=\bm H_o$. Moreover, the functions $H_i$ are never vanishing.
\end{cor}
\proof The linear Pfaffian system \eqref{pfaff} is completely integrable, by Lemma \ref{def1}. This ensures uniqueness and existence of solutions $H_i$. The non-vanishing of solutions is a standard result, see e.g. \cite[Ch.\,11]{Har20}.  \endproof
If we compare equations \eqref{pfaff} with the equations of Proposition \ref{propLame}, we notice that an homogeneity condition of the functions $H_i(\bm u)$ is missing. In general, for arbitrary choices of the initial datum $\bm H_o$, there are no constants $\dl_1,\dots,\dl_n$ such that
\beq\label{furtherconstraint}
\sum_{j=1}^nu^j\der_jH_i=-\dl_i H_i,\quad i=1,\dots,n.
\eeq
The following lemma will clarify how to recover \eqref{furtherconstraint}.
\begin{lem}\label{magicV}
Let $\mathbb V$ be the linear space of (column) vector solutions $\bm H(\bm u)=(H_1(\bm u),\dots, H_n(\bm u))^T$, analytic in $\bm u\in\mc V\setminus\Theta$, of the Pfaffian system \eqref{pfaff}.
The matrix $V(\bm u)$ acts on $\mathbb V$. Moreover, if $V(\bm u)\bm H(\bm u)=\ell \bm H(\bm u)$, with $\ell\in\C$, we have
\[\sum_{j=1}^nu^j\der_jH_i=(\ell+b_i) H_i,\quad i=1,\dots,n.
\]
\end{lem}
\proof
The matrices $V_i(\bm u)$ defined in \eqref{rhjs1} have entries $(V_i)^j_h=\Gm^j_{i}\dl_{ih}-\dl_{ji}\Gm^i_h$. Thus, the joint system \eqref{pfaff} can be written in matrix form $\der_k\bm H=V_k\bm H$, for any $k=1,\dots,n$. We have
\[\der_k(V\bm H)=(\der_kV)\bm H+V\der_k\bm H=[V_k,V]\bm H+VV_k\bm H=V_kV\bm H,
\]where we have used the equation $\der_kV=[V_k,V]$ (compatibility condition of the joint system \eqref{rhjs1},\eqref{rhjs2}). If $V\bm H=\ell\bm H$, we deduce
\[
\pushQED{\qed}
\sum_{j=1}^nu^j\der_j\bm H=\sum_{j=1}^nu^jV_j\bm H=\sum_{j=1}^nu^j[\Gm,E_j]\bm H=[\Gm,U]\bm H=(V+B)\bm H=(\ell\cdot{\bf 1}+B)\bm H.\qedhere
\popQED
\]
\begin{rem}
The incompatibility of the joint system of equations \eqref{pfaff} and \eqref{furtherconstraint}, for arbitrary choices of $\bm H_o$ and $\dl_1,\dots,\dl_n$ was already noticed in \cite[Section 2]{Lor14}, in the construction procedure of  bi-flat $F$-manifold starting from a solution of the Darboux--Egoroff system \eqref{DE1},\eqref{DE2},\eqref{DE3}. Lemma \ref{magicV} is essentially formulated in \cite[Remark 2.3]{Lor14}.
\end{rem}
\begin{cor}\label{magichoice}
Assume that 
\beq\label{assume}
\prod_{j=1}^n\left[G_0(\bm u_o)^{-1}\right]^1_j\neq 0.
\eeq
Then the functions $H_i(\bm u):=\left[G_0(\bm u)^{-1}\right]^1_i$, with $i=1,\dots,n$, are non-vanishing in a sufficiently small neighborhood of $\bm u_o$, and satisfy the joint system \eqref{pfaff},\eqref{furtherconstraint} with $\dl_i=D_1+L_{11}-b_i$, for $i=1,\dots, n$.
\end{cor}
\proof Let $\bm H(\bm u)$ be the first column of $(G_0(\bm u)^{-1})^T$. From \eqref{rhjs1}, we deduce $\der_k\bm H=V_k\bm H$ for any $k=1,\dots, n$. This is equivalent to \eqref{pfaff}. Moreover, from \eqref{rhjs2}, we deduce $V(\bm u)(G_0(\bm u)^{-1})^T=(G_0(\bm u)^{-1})^T(-D-L')$. By taking the first column of both sides, we obtain $V\bm H=\ell\bm H$ with $\ell=-D_1-L_{11}$. We conclude by Lemma \ref{magicV}.
\endproof
Thus, if the non-vanishing assumption \eqref{assume} holds, we have an automatic choice for a solution $\bm H(\bm u)$ of the joint system of equations \eqref{pfaff},\eqref{furtherconstraint}, which is entry-by entry non-vanishing in a sufficiently small neighborhood of $\bm u_o$. Such a choice is given by the first column of $[G_0(\bm u)^{-1}]^T$. In what follows, we will assume \eqref{assume} and we will fix such a choice for $\bm H$.
\begin{lem}
Let $G_0(\bm u),H_i(\bm u)$ as above. For any $\al=1,\dots,n$, the one form 
\[\varpi_\al(\bm u):=\sum_{k=1}^n G_0(\bm u)^k_\al H_k(\bm u)du^k\] is closed.
\end{lem}
\proof
Set $H:={\rm diag}(H_1,\dots, H_n)$. We have $\der_j[G_0^TH]^\al_i=(G_0^T)^\al_j H_i\Gm^j_i+(G_0^T)^\al_i H_j\Gm^i_j=\der_i[G_0^TH]^\al_j$. This can be easily seen by invoking equations \eqref{rhjs1}, and \eqref{pfaff}. 
\endproof

\begin{lem}
Let $G_1(\bm u),H_i(\bm u)$ as above.
For any $\al=1,\dots, n$, the one form 
\[\phi_\al(\bm u):=\sum_{k=1}^n G_1(\bm u)^k_\al H_k(\bm u)du^k\] is closed.  
\end{lem}
\proof
Firstly, notice that we have $\der_iG_1=E_iG_0+\der_iG_0\cdot G_0^{-1}\cdot G_1$. This follows from the fact that $X_0$ is a solution of the joint system \eqref{rhjs1}, \eqref{rhjs2}. Furthermore, we have the identity $\phi_\nu(\bm u):=\sum_\al [G_0(\bm u)^{-1}G_1(\bm u)]^\al_\nu\cdot\varpi_\al(\bm u)$. Hence, we deduce
\[
\pushQED{\qed}
d\phi_\nu=\sum_{\al,i} \der_i[G_0(\bm u)^{-1}G_1(\bm u)]^\al_\nu\ du^i\wedge\varpi_\al=\sum_{\al,i,k}(G_0^{-1})^\al_i\ (G_0)^i_\nu\ (G_0)^k_\al\ H_k\ du^i\wedge du^k=0.\qedhere
\popQED
\]
Consider a polydisc $D(\bm u_o)\subseteq\mc V\setminus\Theta$ centered in $\bm u_o$, and sufficiently small so that $H_i(\bm u)\neq 0$ for $\bm u\in D(\bm u_o)$, and $i=1,\dots,n$. Define the functions
\beq\label{recprocedure}
t^\al(\bm u):=\int_{\bm u_o}^{\bm u}\varpi_\al,\quad F^\nu(\bm u):=\int_{\bm u_o}^{\bm u}\phi_\nu,\quad\al,\nu=1,\dots, n,
\eeq
where $\bm u\in D(\bm u_o)$. By definition, the functions $\bm t=(t^\al)_\al$ define a system of coordinates on $D(\bm u_o)$, with Jacobian matrix
\beq\label{recpsitilde}
\left(\frac{\der t^\al}{\der u^i}\right)_{\al,i}=G_0(\bm u)^T H(\bm u).
\eeq
Introduce the connection $\nabla$ on the holomorphic tangent bundle of $D(\bm u_o)$ by 
\beq\label{30.10.23-1}
\nabla_{\frac{\der}{\der u^i}}\frac{\der}{\der u^j}=\sum_{h=1}^nK^h_{ij}\frac{\der}{\der u^h},\quad i,j=1,\dots,n,
\eeq
where the Christoffel symbols $K^h_{ij}$, with $i,j,h=1,\dots,n$, are defined by
\begin{align}
\label{30.10.23-2}
&K^h_{ij}=0,&i,j,k\text{ distinct,}\\
\label{30.10.23-3}
&K^i_{ij}=K^i_{ji}=-K^i_{jj}=\frac{H_j}{H_i}\Gm^i_j,&i\neq j,\\
\label{30.10.23-4}
&K^i_{ii}=-\sum_{h\neq i}K^i_{ih}.
\end{align}Moreover, define the commutative and associative product $\circ$ of holomorphic vector fields on $D(\bm u_o)$ by
\beq
\label{30.10.23-5}
\frac{\der}{\der u^i}\circ \frac{\der}{\der u^j}=\sum_{h=1}^nc_{ij}^h\frac{\der}{\der u^h},\quad c_{ij}^h=\dl^h_i\dl^h_j,\quad i,j,h,=1,\dots,n,
\eeq
with unit $e=\sum_{j=1}^n\frac{\der}{\der u^j}$.
\begin{thm}\label{enF}
The datum $(D(\bm u_o),\nabla,c,e)$ defines an analytic semisimple homogeneous flat $F$-manifold. More precisely:
\begin{enumerate}
\item The coordinates $\bm t$ are flat coordinates, the coordinates $\bm u$ are canonical coordinates.
\item The coordinate $t^1$ is such that $\frac{\der}{\der t^1}=e$.
\item The functions $F^\al(\bm u(\bm t))$ are solutions of the oriented associativity equations \eqref{wdvv2}, and equations \eqref{wdvv1} hold for $A^1=1,A^\mu=0$ for $\mu\neq 1$.
\item The $\Tilde\Psi$-matrix is $\Tilde\Psi(\bm u)=\left(G_0(\bm u)^TH(\bm u)\right)^{-1}$.
\item The vector field $E=\sum_iu^i\frac{\der}{\der u^i}$ defines an Euler vector field.
\item The tensor $\mu^\la$, with $\la=D_{1}+L_{11}$, is given by $\mu^\la=-D-L'$.
\item The homogeneous degrees $q_\al$ are given by $q_\al=-D_\al-L_{\al\al}+D_{1}+L_{11}$ for $\al=1,\dots,n.$
\item The conformal dimensions are $\dl_j=D_1+L_{11}-b_j$, for $j=1,\dots,n$.
\item The monodromy data $(\mu^\la,[R],\mathring{S}_1,\mathring{S}_2,\La,\mathring{C})$ at $\bm u_o$, computed with respect to $\la=D_1+L_{11}$, $\tau$, the Lam\'e coefficients $(H_i)_{i=1}^n$, can be reconstructed from the admissible datum $\frak M$:
\[%
\mu^\la=-D-L',\quad [R]=[L''],\quad\mathring{S}_1=S_1,\quad \mathring{S}_2=S_2,\quad \La=B,\quad \mathring C=C.
\]
\end{enumerate}
\end{thm}
\proof
The proof that $(D(\bm u_o),\nabla,c,e)$ is a flat $F$-manifold is the same as in \cite[Th.\,2.2, Lem.\,2.3]{AL13}. More precisely, the flatness of $\nabla$ is checked via explicit computations of the Riemann curvature tensor. The identity $\nabla e=0$ is equivalent to the identities $K^i_{jj}=-K^j_{ij}$, for $i\neq j$, and $K^i_{ii}=-\sum_{h\neq i}K^i_{ih}$ for the Christoffel symbols. The symmetry of $\nabla c$ in its covariant (i.e.\,lower) indices is equivalent to the identities $K^i_{jj}=-K^j_{ij}$, for $i\neq j$, and $K^h_{ij}=0$, for $i,j,k$ distinct. Furthermore, the flat $F$-structure is semisimple, and $\bm u$ are canonical coordinates, as \eqref{30.10.23-5} manifestly shows.
\vskip1mm
The matrix $M(\bm u):=G_0^T(\bm u)H(\bm u)$ satisfies
\begin{align*}
\der_i M&=\der_i G_0^T\cdot H+G_0^T\der_iH\\
&=-G_0^TV_iH+G_0^T\der_iH\qquad\text{by equation \eqref{rhjs1}}\\
&=M\left(-H^{-1}V_iH+H^{-1}\der_iH\right).
\end{align*}
Arrange the Christoffel symbols $K^h_{ij}$ into $n$ matrices $K_1,\dots, K_n$ by setting $(K_i)^h_j=K^h_{ij}$, for $i,j,h=1,\dots,n$. 
We claim that 
\beq
\label{30.10.23-6}
-H^{-1}V_iH+H^{-1}\der_iH=K_i,\quad\text{for any $i=1,\dots, n$.}
\eeq
Indeed, the entry $(h,j)$ of the l.h.s.\,\,of \eqref{30.10.23-6} equals
\begin{empheq}[left={-\frac{H_j}{H_h}(V_i)^h_j+\frac{1}{H_h}\der_iH_j\dl^h_j=-\frac{H_j}{H_h}\left(\Gm^h_i\dl_{ij}-\Gm^i_j\dl_{ih}\right)+\frac{1}{H_h}\der_iH_j\dl^h_j}{=}\empheqlbrace]{align*} 
       &0,\quad\text{ $i,j,h$ distinct,}\\
       &\frac{H_j}{H_i}\Gm^i_j,\quad h=i\neq j,\\
       &-\frac{H_i}{H_h}\Gm^h_i,\quad h\neq i=j,\\
       &-\sum_{\ell\neq i}\frac{H_{\ell}}{H_i}\Gm^i_\ell,\quad i=j=h.
 \end{empheq}
    This proves \eqref{30.10.23-6}. Since $M^\al_j=\frac{\der t^\al}{\der u^j}$, for $\al,j=1,\dots,n$, the identity $M^{-1}\der_iM=K_i$ reads
    \[K^h_{ij}=\sum_{\al}\frac{\der u^h}{\der t^\al}\frac{\der}{\der u^i}\left(\frac{\der t^\al}{\der u^j}\right).
    \]By recalling the transformation rule for Christoffel symbols, this proves that in the coordinates $\bm t$ all the Christoffel symbols of $\nabla$ vanish. So $\bm t$ are $\nabla$-flat coordinates. Since the coefficients $H_i$ are chosen as in Corollary \ref{magichoice}, that is $H_i(\bm u)=\left[G_0(\bm u)^{-1}\right]^1_i$ for $i=1,\dots,n$, we have
    \[\frac{\der u^i}{\der t^1}=\left[H^{-1}(G_0^{-1})^T\right]^i_1=\frac{1}{H_i}(G_0^{-1})^1_i=1,\quad i=1,\dots,n,\qquad\text{so that }\frac{\der}{\der t^1}=\sum_{j=1}^n\frac{\der}{\der u^j}=e.
    \]The functions $F^\al(\bm u(\bm t))$ are the potentials, by Corollary \ref{corhf}. The validity of equations \eqref{wdvv1}, with $A^1=1$ and $A^\mu\neq 0$ for $\mu\neq 1$, is clear. Point (4) is clear from \eqref{recpsitilde}. All the remaining claims directly follows by inspection of the joint system \eqref{rhjs1},\eqref{rhjs2} written in the coordinate frame $(\frac{\der}{\der t^\al})_{\al=1}^n$, Remarks \ref{remt1} and \ref{magicrk2}, and by the definition of the monodromy data and of the Riemann-Hilbert-Birkhoff problem $\mc P[\bm u_o,\tau,\frak M]$.
\endproof

\begin{defn}
We denote by $\mc F[\bm u_o,\tau,\frak M]$ the germ (pointed at $\bm u_o$) of the analytic flat $F$-manifold $(D(\bm u_o),\nabla,c,e)$ described in Theorem \ref{enF}.
\end{defn}

\begin{prop}\label{propisorec1}Let $h={\rm diag}(h_1,\dots,h_n)\in GL(n,\C)$. %
If $\frak M=(B,D,L,S_1,S_2,C)$ is $(\bm u_o,\tau)$-admissible, then also the datum $$h\,\frak M:=(B,\quad D,\quad L,\quad h^{-1}S_1h,\quad h^{-1}S_2h,\quad Ch)$$
is $(\bm u_o,\tau)$-admissible. If $\mc P[\bm u_o,\tau,\frak M]$ is solvable, then so is $\mc P[\bm u_o,\tau,h\,\frak M]$. Moreover, if the solution of $\mc P[\bm u_o,\tau,\frak M]$ satisfies the assumption \eqref{assume}, then also the solution of $\mc P[\bm u_o,\tau,h\,\frak M]$ does. The resulting pointed flat $F$-manifolds $\mc F[\bm u_o,\tau,\frak M]$ and $\mc F[\bm u_o,\tau,h\,\frak M]$ coincides, being defined by the same potentials $(F^1(\bm t),\dots, F^n(\bm t))$, up to linear terms.
\end{prop}
\proof
The admissibility of $h\,\frak M$ is easily checked. Denote by $\mc H(z,\bm u)$ and $\mc H'(z,\bm u)$ the coefficients of the problems $\mc P[\bm u,\tau,\frak M]$ and $\mc P[\bm u,\tau,h\,\frak M]$, respectively. We have
 \begin{empheq}[left=\mc H'{:=}\empheqlbrace]{align*} 
       &h^{-1}\mc H h,\quad\text{ along }\Gm_{-\infty},\\
       &h^{-1}\mc H h,\quad\text{ along }\Gm_{+\infty},\\
       &h^{-1}\mc H,\quad\text{ along }\Gm_{1},\\
       &h^{-1}\mc H,\quad\text{ along }\Gm_{2}.
    \end{empheq}
    
   If $G(z,\bm u)$ denotes the solution of the problem of $\mc P[\bm u,\tau,\frak M]$, with $\bm u$ varying in a sufficiently small neighborhood of $\bm u_o$, then the function 
    \begin{empheq}[left={ G'(z,\bm u)}{:=}\empheqlbrace]{align*} 
       &h^{-1} G(z,\bm u) h,\quad\text{ on }\Pi_{L/R},\\
       &h^{-1} G(z,\bm u),\quad\text{ on }\Pi_0
    \end{empheq}
    is the solution of $\mc P[\bm u,\tau,h\,\frak M]$. We then have the following rescaling:
    \begin{align*}
    G_0(\bm u)&\mapsto G_0'(\bm u)=h^{-1}G_0(\bm u),\\
     G_1(\bm u)&\mapsto G_1'(\bm u)=h^{-1}G_1(\bm u),\\
      H_j(\bm u)&\mapsto H_j'(\bm u)= h_jH_j(\bm u),\quad j=1,\dots,n.
    \end{align*}
    The forms $\varpi_\al(\bm u)$ and $\phi_\al(\bm u)$ are invariant. 
\endproof

\begin{prop}\label{propisorec2}
Let $\bm G={\rm diag}(g_1,\dots, g_n)\in GL(n,\C)$. If $\frak M=(B,D,L,S_1,S_2,C)$ is $(\bm u_o,\tau)$-admissible, then also the datum $$\bm G\,\frak M:=(B,\quad D,\quad \bm G^{-1}L\bm G,\quad S_1,\quad S_2,\quad \bm G^{-1}C)$$
is $(\bm u_o,\tau)$-admissible. If $\mc P[\bm u_o,\tau,\frak M]$ is solvable, then so is $\mc P[\bm u_o,\tau,\bm G\,\frak M]$. Moreover, if the solution of $\mc P[\bm u_o,\tau,\frak M]$ satisfies the assumption \eqref{assume}, then also the solution of $\mc P[\bm u_o,\tau,\bm G\,\frak M]$ does. The resulting pointed flat $F$-manifolds $\mc F[\bm u_o,\tau,\frak M]$ and $\mc F[\bm u_o,\tau,\bm G\,\frak M]$ are isomorphic: if $(F^\al(\bm t))_{\al=1}^n$ (resp.\,\,$(\Tilde F^\al(\tilde{\bm t}))_{\al=1}^n$) are the potentials defining the structure $\mc F[\bm u_o,\tau,\frak M]$ (resp.\,\,$\mc F[\bm u_o,\tau,\bm G\,\frak M]$), we have
\[\tilde t^\al=\frac{g_\al}{g_1} t^\al+c^\al,\quad \Tilde F^\al(\tilde{\bm t})=\frac{g_\al}{g_1} F^\al(\bm t)+\text{linear terms in $\bm t$},\quad\al=1,\dots,n.
\]In particular, $\frac{\der}{\der t^1}=\frac{\der}{\der \tilde t^1}=e$.
\end{prop}

\proof
The admissibility of $\bm G\,\frak M$ is easily checked. Denote by $\mc H(z,\bm u)$ and $\mc H'(z,\bm u)$ the coefficients of the problems $\mc P[\bm u,\tau,\frak M]$ and $\mc P[\bm u,\tau,\bm G\,\frak M]$, respectively. We have
\[\mc H'=\mc H\quad \text{along }\Gm_{\pm\infty},\qquad \mc H'=\mc H\bm G\quad\text{along }\Gm_1,\Gm_2.
\]If $G(z,\bm u)$ denotes the solution of the problem of $\mc P[\bm u,\tau,\frak M]$, with $\bm u$ varying in a sufficiently small neighborhood of $\bm u_o$, then the function 
\[G'(z,\bm u):=G(z,\bm u)\quad\text{ on }\Pi_{L/R},\qquad G'(z,\bm u)=G(z,\bm u)\bm G\quad\text{on }\Pi_0,
\]is the solution of $\mc P[\bm u_o,\tau,\bm G\,\frak M]$. We then have the following rescaling
 \begin{align*}
    G_0(\bm u)&\mapsto G_0'(\bm u)=G_0(\bm u)\bm G,\\
     G_1(\bm u)&\mapsto G_1'(\bm u)=G_1(\bm u)\bm G,\\
      H_j(\bm u)&\mapsto H_j'(\bm u)= g_1^{-1}H_j(\bm u),\quad j=1,\dots,n.
    \end{align*}
    The forms $\varpi_\al(\bm u)$ and $\phi_\al(\bm u)$ transform as follows: if $\bm G={\rm  diag}(g_1,\dots,g_n)$, then 
    \[\varpi_\al(\bm u)\mapsto \varpi_\al'(\bm u)=\frac{g_\al}{g_1}\varpi_\al(\bm u),\quad %
    \phi_\al(\bm u)\mapsto\phi_\al'(\bm u)=\frac{g_\al}{g_1}\phi_\al(\bm u),\quad  %
    \al=1,\dots,n.
    \]The claim follows.
\endproof

\subsection{Reconstruction of admissible germs of semisimple flat $F$-manifold } In this section we prove that all admissible germs of analytic flat $F$-manifolds are of the form $\mc F[\bm u_o,\tau,\frak M]$. %
\vskip1mm
Let $(M,p)$ be an admissible germ of an analytic flat $F$-manifold. Without loss of generality, we fix flat local coordinates $\bm t$ so that $\frac{\der}{\der t^1}=e$. Fix choices of normalizations (1)-(6) of Section \ref{chnor}. We then have a well defined system $(\la,\mu^\la,R,S_1,S_2,\La,C)$ of monodromy data at $p$, computed w.r.t.\,the chosen ordering of $\bm u_o=\bm u(p)$, an admissible direction $\tau$, and the normalization $(H_{o,1},\dots, H_{o,n})\in(\C^*)^n$ of Lam\'e coefficients at $p$.

\begin{lem}
The matrix $\mu^\la={\rm diag}(q_1,\dots, q_n)-\la\cdot{\bf 1}$ has a unique decomposition $\mu^\la=D^\la+S^\la$ with
\bea
&D^\la={\rm diag}(d_1,\dots, d_n),& d_i\in\Z,\quad i=1,\dots,n\\
&S^\la={\rm diag}(\rho_1,\dots, \rho_n),& {\rm Re}(\rho_i)\in [0,1[,\quad i=1,\dots,n.
\eea
We have 
\beq\label{commrel}
[D^\la,S^\la]=0,\qquad [R,S^\la]=0.
\eeq
\end{lem}
\proof
The uniqueness of the decomposition, and the first commutation relation of \eqref{commrel} are clear. For any pair $(i,j)$, we have $[S^\la, R]_{ij}=(\rho_i-\rho_j)R_{ij}=0$. Indeed, for $i\neq j$ we have two possibilities: if $q_i-q_j\notin\Z_{<0}$, then $R_{ij}=0$; if $q_i-q_j\in\Z_{<0}$, then $\rho_i-\rho_j=0$.
\endproof

\begin{rem}
It follows that $z^{-\mu^\la}z^R=z^{-D^\la}z^{R-S^\la}$.
\end{rem}

\begin{prop}\label{propmonodatum}
The $6$-tuple $\frak M=(\La, -D^\la, R-S^\la, S_1,S_2,C)$ is $(\bm u_o,\tau)$-admissible. The problem $\mc P[\bm u_o,\tau,\frak M]$ is solvable. Moreover, the solution $G(z,\bm u_o)$ satisfies the non-vanishing assumption \eqref{assume}.
\end{prop}
\proof
Point (1),(2),(3),(4) of Definition \ref{datum} directly follow from the definitions and properties of $R,\La,D^\la,S^\la$. Points (5)-(10) of Definition \ref{datum} follow from Propositions \ref{propsc1} and \ref{propsc2}. Let 
\begin{itemize}
\item $\Xi_0(\bm t,z)$ be the fixed solution of the joint system \eqref{jfs} in Levelt normal form,
\item $X_i(\bm u,z)$, with $i=1,2,3$, be the solutions of the joint system \eqref{jfss} uniquely defined by the asymptotics \eqref{asym}.
\end{itemize}
The (unique) solution of the problem $\mc P[\bm u_o,\tau,\frak M]$ is
\begin{empheq}[left={G(z,\bm u_o)=\empheqlbrace\,}]{align*}
      & \left(\widetilde{\Psi}(\bm t_o) H_o^{-1}\right)^T\Xi_0(\bm t_o,z)z^{-R}z^{\mu^\la}, &z\in\Pi_0,&\qquad\qquad\\
      & X_2(\bm u_o,z)e^{-zU_o}z^{-\La}, &z\in\Pi_R,&\qquad\qquad\\
      & X_3(\bm u_o,z)e^{-zU_o}z^{-\La}, &z\in\Pi_L.&\qquad\qquad
\end{empheq}
Here we set $\bm t_o:=\bm t(p)$, $U_o:={\rm diag}(u_o^1,\dots, u_o^n)$, and $ H_o:={\rm diag}(H_{o,1},\dots, H_{o,n})$. Notice that assumption \eqref{assume} holds, since $\Tilde\Psi(\bm t_o)_1^i=1$ for any $i=1,\dots,n$, see Remark \ref{remtildepsi1}.
\endproof

\begin{rem}
In \cite{Cot20b} it is proved that solutions of $\mc P[\bm u,\tau,\frak M]$ can be factorized via two auxiliary RHB problems $\mc P_1[\bm u,\tau,\frak M]$ and $\mc P_2[\bm u,\tau,\frak M]$. The problem $\mc P_1[\bm u,\tau,\frak M]$ is shown to admit unique solution $\Psi(z,\bm u)$ holomorphically depending on $\bm u$ varying in a neighborhood of $\bm u_o$, see \cite[Th.\,3.7]{Cot20b}. Given $\Psi(z,\bm u)$, the problem $\mc P_2[\bm u,\tau,\frak M]$ is formulated, and it is shown to be locally uniquely solvable, see \cite[proof of Th.\,3.13]{Cot20b}.

If $\bm u_o\in\Dl$, assumption (10) of Definition \ref{datum} is crucial for the proof of the unique solvability of $\mc P_1[\bm u,\tau,\frak M]$, see \cite[proof of Lem.\,3.6]{Cot20b}. If $(M,p)$ is not an admissible germ, then the monodromy data are not well defined. Indeed Theorem \ref{asymsol} does not hold, solutions $X_i(\bm u_o,z)$, with $i=1,2,3$, are not unique. With each such a triple of solutions we can associated a pair $(S_1,S_2)$ of Stokes matrices. In general these Stokes matrices do not satisfy Proposition \ref{propsc2}.
\end{rem}

\begin{thm}\label{mthm1}
The analytic germ $\mc F[\bm u_o,\tau,\frak M]$ %
of flat $F$-manifold is isomorphic to the original admissible germ $(M,p)$. They are defined by the same oriented associativity potentials (modulo linear terms).  \qed %
\end{thm}

In the light of the construction of Section \ref{RHBcon}, Theorem \ref{mthm1} follows from the following crucial result.

\begin{lem}\label{miracle}
Let one of the following assumptions hold:
\begin{enumerate}
\item $\bm u_o\in\C^n\setminus\Dl$,
\item $\bm u_o\in\Dl$, and if $u^i_o=u^j_o$, with $i,j\in\{1,\dots,n\}$, $i\neq j$, then $\dl_i-\dl_j\notin\Z\setminus\{0\}$. %
\end{enumerate}
Let
\[ F(\bm u)= F_o+\sum_{k=1}^\infty\sum_{\ell_1,\dots,\ell_k=1}^n \frac{1}{k!}F^{(\bm \ell)}\prod_{j=1}^k\overline{u}_{\ell_j},\quad \overline{u}_i:=u_i-u_{o,i},
\]be a matrix-valued formal power series whose off-diagonal entries $F^i_{j}$ are formal solutions of the Darboux-Egoroff system \eqref{DE1}, \eqref{DE2}, \eqref{DE3}, \eqref{midgm}, \eqref{midgm2}. The off-diagonal entries of the coefficients $F^{(\bm \ell)}$ can be uniquely reconstructed from the off-diagonal entries of $F_o$.
\end{lem}

\proof We have to show that the derivatives $\der_{i_1}\dots\der_{i_N}F^i_{j}(\bm u_o)$ can be computed from the only knowledge of the numbers $F^i_{j}(\bm u_o)$. We proceed by induction on $N$. Let us start with the case $N=1$.\newline
{\bf Step 1.} For $i,j,k$ distinct, by expanding both sides of $\der_kF^i_j=F^i_kF^k_j$ in power series, and equating the coefficients, one reconstructs the coefficients of $\der_kF^i_j(\bm u_o)$. \newline
{\bf Step 2.} From the identities \eqref{midgm} and \eqref{midgm2}  for $F^i_j$, one can compute $\der_i F^i_j(\bm u_o)$ and $\der_j F^i_j(\bm u_o)$ provided that $u_{o,i}\neq u_{o,j}$.\newline
{\bf Step 3. }Assume that $u_{o,i}=u_{o,j}$. By taking the $\der_i$-derivative of both sides of \eqref{midgm} we obtain
\beq
\label{diDE4}
(\dl_j-\dl_i-2)\der_iF^i_j+(u^j-u^i)\der^2_iF^i_j=\sum_{k\neq i,j}(u^k-u^j)[\der_iF^i_kF^k_j+F^i_k\der_iF^k_j].
\eeq
By evaluating \eqref{diDE4} at $\bm u=\bm u_o$ we can compute all the numbers $\der_i F^i_{j}(\bm u_o)$, namely 
\[\der_i F^i_{j}(\bm u_o)=\frac{1}{\dl_j-\dl_i-2}\sum_{k\neq i,j}(u_{o}^k-u_{o}^j)\left[\der_i F^i_{k}(\bm u_o) F^k_{j}(\bm u_o)+ F^i_{k}(\bm u_o)F^k_{i}(\bm u_o) F^i_{j}(\bm u_o)\right].
\]
Notice that the only terms $\der_i F_{k}^i(\bm u_o)$ appearing in this sum are those computed in Step 2.\newline
{\bf Step 4.} If $u_{o,i}=u_{o,j}$, the numbers $\der_jF^i_j(\bm u_o)$ can be computed similarly as in Step 3, by invoking equation \eqref{midgm2}:
\[\der_j F^i_{j}(\bm u_o)=\frac{1}{\dl_j-\dl_i-2}\sum_{k\neq i,j}(u_{o}^k-u_{o}^i)\left[F^i_{j}(\bm u_o)F^j_{k}(\bm u_o) F^k_{j}(\bm u_o)+F^i_k(\bm u_o)\der_j F^k_{j}(\bm u_o)\right].
\]
\newline
This proves that all the first derivatives $\der_k F_{j}^i(\bm u_o)$ can be computed.
\vskip2mm
\noindent{\bf Inductive step. }Assume to know all the $N$-th derivatives $\der_{i_1}\dots \der_{i_N} F_{j}^i(\bm u_o)$. We show how to compute the number $\der_{h_1}\dots \der_{h_{N+1}} F_{j}^i(\bm u_o)$ for any $(N+1)$-tuple $(h_1,\dots, h_{N+1})$.\newline 
{\bf Step 1.} Assume that there exists $\ell\in\{1,\dots, N+1\}$ such that $h_\ell\neq i,j$. We have
\[\der_{h_1}\dots\der_{h_{N+1}}F_{j}^i=\der_{h_1}\dots\der_{h_{\ell-1}}\der_{h_{\ell+1}}\dots\der_{h_{N+1}}[\der_{h_\ell}F_{j}^i]=\der_{h_1}\dots\der_{h_{\ell-1}}\der_{h_{\ell+1}}\dots\der_{h_{N+1}}[F_{h_{\ell}}^iF_{j}^{h_{\ell}}].
\]By evaluation at $\bm u=\bm u_o$, we can compute all the numbers $\der_{h_1}\dots\der_{h_{N+1}}F_{j}^i(\bm u_o)$.
\vskip2mm
\noindent Now we need to compute the mixed derivatives $\der_i^p\der_j^{N+1-p}F_{j}^i(\bm u_o)$, with $0\leq p\leq N+1$. \newline
{\bf Step 2.} Assume $p>0$ and $u_{o,i}\neq u_{o,j}$. Take the $\der_i^{p-1}\der_j^{N+1-p}$-derivative of both sides of \eqref{midgm}: by evaluation at $\bm u=\bm u_o$ we can compute the numbers $\der_i^{p}\der_j^{N+1-p}F_{j}^i(\bm u_o)$.\newline
{\bf Step 3.} Assume $p>0$ and $u_{o,i}= u_{o,j}$. Take the $\der_i^{p-1}\der_j^{N+1-p}$-derivative of both sides of \eqref{diDE4}, to obtain
\begin{multline}\label{lisboa1}(\dl_j-\dl_i-p-1)\der_i^p\der_j^{N+1-p}F_{j}^i+(N+1-p)\der_i^{p+1}\der_j^{N-p}F_{j}^i+(u^j-u^i)\der_i^{p+1}\der_j^{N+1-p}F_{j}^i\\
=\der_i^{p-1}\der_j^{N+1-p}\sum_{k\neq i,j}(u^k-u^j)[\der_iF^i_kF^k_j+F^i_k\der_iF^k_j].
\end{multline}
Specialize \eqref{lisboa1} for $p=N+1$: by evaluation at $\bm u=\bm u_o$ of both sides, we can compute the derivative $\der_i^{N+1}F_{j}^i(\bm u_o)$.\newline
Specialize \eqref{lisboa1} for $p=N$: by evaluation at $\bm u=\bm u_o$ of both sides, we can compute the derivative $\der_i^{N}\der_jF_{j}^i(\bm u_o)$.\newline
Repeating this procedure, by decreasing $p\mapsto p-1$ at each step, we can compute all the mixed derivatives $\der_i^p\der_j^{N+1-p}F_{j}^i(\bm u_o)$.\newline
{\bf Step 4.} Assume $p=0$. The derivative $\der_j^{N+1}F_{j}^i(\bm u_o)$ can be computed, as in Steps 2 and 3, by invoking equation \eqref{midgm2}.
\vskip2mm
\noindent This proves that all the $(N+1)$-th derivatives $\der_{h_1}\dots\der_{h_{N+1}}F_{j}^i(\bm u_0)$ can be computed.
\endproof

\subsection{Convergence of semisimple admissible formal flat $F$-manifolds}\label{seconv}
We are now ready to prove the following result.

\begin{thm}\label{thconv}
Let $(H,\bm \Phi)$ be an admissible formal semisimple flat $F$-manifold over $\C$, with Euler field $E$. The oriented associativity potentials $\bm \Phi=(\Phi^1,\dots, \Phi^n)$ have a non-empty common domain of convergence.
\end{thm}
\proof
 Fix one ordering $\bm u_o\in\C^n$ of the eigenvalues of the operator $\mc U(\bm t)$ at $\bm t=0$. We have $n\times n$ matrix-valued (a priori) formal power series in $\bm u$
 \bea
 V(\bm u)= V_o+\sum_{k=1}^\infty\sum_{\ell_1,\dots,\ell_k=1}^n \frac{1}{k!}V^{(\bm \ell)}\prod_{j=1}^k\overline{u}_{\ell_j},\quad\quad
 V_i(\bm u)= V_{i,o}+\sum_{k=1}^\infty\sum_{\ell_1,\dots,\ell_k=1}^n \frac{1}{k!}V_i^{(\bm \ell)}\prod_{j=1}^k\overline{u}_{\ell_j},\\
 \Psi(\bm u)= \Psi_o+\sum_{k=1}^\infty\sum_{\ell_1,\dots,\ell_k=1}^n \frac{1}{k!}\Psi^{(\bm \ell)}\prod_{j=1}^k\overline{u}_{\ell_j},\quad\quad
 \Gm(\bm u)= \Gm_o+\sum_{k=1}^\infty\sum_{\ell_1,\dots,\ell_k=1}^n \frac{1}{k!}\Gm^{(\bm \ell)}\prod_{j=1}^k\overline{u}_{\ell_j},
 \eea
 where $\overline{u}_i:=u_i-u_{o,i}$ for $i=1,\dots,n$. These power series are well defined by the semisimplicity assumption, and they satisfy properties described in Propositions \ref{fundprop1} and \ref{fundprop2}.
 \vskip1mm
 Set $\bm H_o:=\Psi_o\widetilde\Psi_o^{-1}$, where $\widetilde\Psi_o:=\left(\left.\frac{\der u^i}{\der t^\al}\right|_{\bm t=0}\right)_{i,\al=1}^n$.
\vskip1mm
After fixing choices of normalizations of Section \ref{chnor}, we can introduce a system of monodromy data $(\la,\mu^\la,R,S_1,S_2,\La,C)$ for the formal flat $F$-structure, computed w.r.t.\,an admissible direction $\tau$ at $\bm u_o$, and the normalization $\bm H_o$ of Lam\'e coefficients at the origin.
Proposition \ref{propmonodatum} holds true, with the same proof. We can set the RHB problem $\mc P[\bm u_o,\tau,\frak M]$. 
This problem is solvable w.r.t.\,\,$\bm u$ on an open neighborhood $\mc V\setminus\Theta$ of $\bm u_o$, by Theorem \ref{teoimp}. The unique solution $G(z;\bm u)$ is holomorphic in $\bm u\in\mc V\setminus \Theta$, and with expansion
\bea
&G(z;\bm u)=I+\frac{1}{z}F_1^{\rm an}(\bm u)+O\left(\frac{1}{z^2}\right),&\quad z\to\infty,\quad z\in\Pi_{L/R},\\
&G(z;\bm u)=G_0(\bm u)+G_1(\bm u)z+G_2(\bm u)z^2+G_3(\bm u)z^3+O(z^4),&\quad z\to0.
\eea
Here the superscript ``an'' stands for {\it analytic}.
As output of Section \ref{RHBcon}, we also obtain a compatible joint system of differential equations (with analytic coefficients in $\bm u$, not just formal) of the form
\beq\label{ansyst}
\frac{\der}{\der u^i}X=\left(zE_i-V_i(\bm u)^T\right)X,\quad \frac{\der}{\der z}X=\left(U-\frac{1}{z}V(\bm u)^T\right)X,
\eeq
where 
\begin{align*}
&V^{\rm an}(\bm u):=[F_1^{\rm an}(\bm u)^T,U]-\La,\\
&V_i^{\rm an}(\bm u):=[F_1^{\rm an}(\bm u)^T,E_i]\equiv -\left(\frac{\der G_0}{\der u^i}\cdot G_0^{-1}\right)^T=\left(\frac{\der }{\der u^i}(G_0^T)^{-1}\right)\cdot G_0^T.
\end{align*}
We also have 
\[
V^{\rm an}(\bm u_o)=V_o,\quad G_0(\bm u_o)=(\Psi_o^{-1})^T,\quad \der_i G_0=V_i^{\rm an}G_0,\quad i=1,\dots,n.
\]
From the equality $[F_1^{\rm an}(\bm u_o)^T,E_i]=V_o=[\Gm_o,E_i]$ we deduce that $[F_1^{\rm an}(\bm u_o)'']^T=\Gm_o$. Moreover, by Lemma \ref{def1}, the off diagonal matrix $[F_1^{\rm an}(\bm u)'']^T$ solve equations \eqref{DE1},\eqref{DE2},\eqref{DE3},\eqref{midgm}, \eqref{midgm2}. By Lemma \ref{miracle}, we obtain $\Gm(\bm u)=[F_1^{\rm an}(\bm u)'']^T$. In particular, $\Gm(\bm u)$ is convergent. It follows that $\Psi(\bm u)=(G_0(\bm u)^{-1})^T$, $V_i(\bm u)=V_i^{\rm an}(\bm u)$, and $V(\bm u)=V^{\rm an}(\bm u)$ are convergent.

The oriented associativity potentials $\Phi^1,\dots,\Phi^n$ can be reconstructed via formulas \eqref{recprocedure}. The original formal structure $(H,\bm \Phi)$ turns out to be equivalent to the analytic flat $F$-manifold $\mc F[\bm u_o,\tau,\frak M]$. %
\endproof
{\bf Open question: }Does it exist a semisimple and strictly doubly resonant %
germ of flat $F$-structure which is purely formal?
\vskip1,5mm
A positive answer would imply the optimality of Theorem \ref{thconv}. The study of the strictly doubly resonant %
germs goes beyond the general theory developed in \cite{CDG1}.

\begin{rem}
Consider a trivial vector bundle $E$ on $\Pb^1$, equipped with a meromorphic connection $\nabla^o$ with connection matrix $\Om$ given by
\[\Om=-\left[U_o+\frac{1}{z}\left(\La+[(F_o'')^T,U_o]\right)\right]dz,\quad U_o={\rm diag}(u_o^1,\dots, u_o^n),\quad \La={\rm diag}(\la-\dl_1,\dots,\la-\dl_n).
\]Malgrange's Theorem \cite{Mal83a,Mal86} asserts that, if $\bm u_o\in\C^n\setminus\Dl$, the connection $\nabla^o$ has a germ of universal deformation. Its connection matrix is
\beq\label{undef}
-d(zU)-([F''(\bm u)^T,U]+\La)\frac{dz}{z}-[F''(\bm u)^T,dU],
\eeq
where $F''(\bm u)$ is the unique off-diagonal solution of the Darboux-Egoroff equations of Lemma \ref{miracle}.\\ 
If $\bm u_o\in\Dl$ and $\dl_i-\dl_j\notin\Z\setminus\{0\}$, C.\,Sabbah proved that $\nabla^o$ admits an integrable deformation of the form \eqref{undef}, see \cite{Sab18}. In this case, the deformation is not universal. Nevertheless, one can prove that there exists a class $\frak I$ of integrable deformations of $\nabla^o$ which are induced by pull-back (via a unique map) by the integrable deformation constructed by Sabbah. We refer the interested reader to \cite{Cot21}, where the class $\frak I$ and its generic elements are described. In both cases ($\bm u_o\notin\Dl$, or $\bm u_o\in\Dl$ with $\dl_i-\dl_j\notin\Z\setminus\{0\}$), the function $F''(\bm u)$ is analytic in a neighborhood of $\bm u_o$.\\
If $\bm u_o\in\Dl$ and $\dl_i-\dl_j\in\Z\setminus\{0\}$ for some $i\neq j$, Lemma \ref{miracle} does not hold, and the initial datum $(\bm u_o,F_o)$ does not identify a unique (formal) solution $F(\bm u)$ of the Darboux-Egoroff equations.
\end{rem}

\subsection{On the number of monodromy local moduli}Consider all $n$-dimensional germs of homogenous semisimple flat $F$-manifolds, modulo local isomorphisms. 

\begin{thm}\label{numbermoduli} 
The local isomorphism classes of $n$-dimensional germs of homogeneous semisimple flat $F$-manifolds generically depend on $n^2-n$ %
parameters.
The local isomorphism classes of $n$-dimensional germs of homogeneous semisimple Frobenius manifolds generically depend on $\frac{1}{2}(n^2-n)$ parameters. 
\end{thm}
{
\proof
We show that germs of flat $F$-manifolds are identifiable with points of a ``stratified space'' $X$, whose generic dimension is $n^2-n$. %
\vskip1mm
Let $(M,p_o)$ be a pointed semisimple homogeneous flat $F$-manifold. For generic $p_o$, the germ $(M,p_o)$ is admissible. After fixing 
\begin{enumerate}
\item a system of local flat coordinates $\bm t$ (such that $\frac{\der}{\der t^1}=e$), 
\item $\la_o\in\C$, 
\item an admissible direction $\tau$ at $p_o$, 
\item the initial value $\bm H_o$ of the Lam\'e coefficients,
\end{enumerate} we can introduce a tuple of monodromy data $\frak M=(\mu^{\la_o},R,S_1,S_2,\La_o,C)$. Conversely, from the knowledge of $(\bm u(p_o),\,\tau,\,\frak M)$ we can reconstruct the germ of the structure together with the above choices (1),(2), and (4): the original germ is indeed isomorphic to $\mc F[\bm u(p_o),\tau,\frak M]$. So, we can use the data $\frak M$ to parametrize the germs of flat $F$-structures on $(\C^n,\bm u(p_o))$.
\vskip1mm
For generic germs, the matrix $\mu^{\la_o}$ has simple spectrum, and with eigenvalues not differing by integers, so that $R=0$. So, generically, we have a tuple of matrices $(\mu^{\la_o}, \La_o,S_1,S_2,C)$ for a total of $2n^2+n$ entries. 
\vskip1mm
We need to impose several constraints, and to cut out redundancies in the counting.  
The matrices $(\mu^{\la_o}, \La_o,S_1,S_2,C)$ must satisfy equation \eqref{c2}, for a total of $n^2$ constraints. Moreover, given $h_1,h_2\in GL(n,\C)$ diagonal, the flat $F$-manifolds $\mc F[\bm u(p_o),\tau,\frak M],\,\,\mc F[\bm u(p_o),\tau,\frak M'],\,\,\mc F[\bm u(p_o),\tau,\frak M'']$ with
\[\frak M'=(\mu^{\la_o},\,\,\La_o, \,\,h_1^{-1}S_1h_1,\,\, h_1^{-1}S_2h_1,\,\, Ch_1),\quad \frak M''=(\mu^{\la_o},\,\,\La_o,\,\, S_1,\,\, S_2,\,\, h_2^{-1}C),
\]coincide up to isomorphisms, by Propositions \ref{propisorec1} and \ref{propisorec2}. In total, the number of free parameters equals $2n^2+n-n^2-2n=n^2-n$.
\vskip1mm
How to choose a set of $n^2-n$ independent parameters out of $(\mu^{\la_o}, \La_o,S_1,S_2,C)$? By equation \eqref{c2}, we have
\[
S_1^{-1}e^{2\pi\sqrt{-1}\La_o}S_2^{-1}\in \mc O(e^{-2\pi\sqrt{-1}\mu^{\la_o}}),
\]
where $\mc O(e^{-2\pi\sqrt{-1}\mu^{\la_o}})$ denotes the similarity orbit of $e^{-2\pi\sqrt{-1}\mu^{\la_o}}$. The codimension of the orbit $\mc O(e^{-2\pi\sqrt{-1}\mu^{\la_o}})$ in ${M}(n,\C)$ equals the dimension of the centralizer
\[\dim \{A\in M_n(\C)\colon\,\,\, [A,e^{-2\pi\sqrt{-1}\mu^{\la_o}}]=0\},
\]see \cite[\S 2.4]{Arn71}. Hence, for generic $\mu^{\la_o}$ we have 
\[\dim \mc O(e^{-2\pi\sqrt{-1}\mu^{\la_o}})=n^2-n.
\]Moreover, it is easy to see that if a matrix $A$ admits a LDU-decomposition\footnote{Here, as in Section \ref{secbrmon}, $L_n$ (resp.\,\,$U_n$) denotes the group of lower (resp.\,\,upper) triangular unipotent $n\times n$-matrices.}
\[A=G_1G_2G_3,\qquad  G_1\in L_n,\quad G_3\in U_n,\quad G_2={\rm diag}(g_1,\dots, g_n), 
\]then such a decomposition is unique, see e.g. \cite{HJ85}. From this, it follows that $(S_1,S_2)$ can be used as coordinates on $\mc O(e^{-2\pi\sqrt{-1}\mu^{\la_o}})$. The total number of parameters $(\mu^{\la_o},S_1,S_2)$ equals $n^2$.
\vskip1mm
Alternatively, one can choose $(\mu^{\la_o},C)$ as coordinates on $\mc O(e^{-2\pi\sqrt{-1}\mu^{\la_o}})$, provided that $C$ is defined up to left multiplication by diagonal invertible matrices. %
The point corresponding to $(\mu^{\la_o},C)$ is
\[
 C^{-1}e^{-2\pi\sqrt{-1}\mu^{\la_o}}C\in \mc O(e^{-2\pi\sqrt{-1}\mu^{\la_o}}).
\]In total, we get $n+(n^2-n)=n^2$ parameters. For parametrizing germs of flat $F$-manifolds:
\begin{itemize}
\item in the tuple $(\mu^{\la_o},S_1,S_2)$ the matrices $(S_1,S_2)$ should be taken up to conjugation $(S_1,S_2)\mapsto (h^{-1}S_1h,\,\, h^{-1}S_2h)$ with $h$ diagonal;
\item in the tuple $(\mu^{\la_o},C)$ should be taken up to right multiplication $C\mapsto Ch$ with $h$ diagonal and invertible.
\end{itemize}In total, we have $n^2-n$ free parameters.
At non-generic points, the space $X$ can get additional strata.
\vskip2mm
Let us also consider the subspace $X_{\rm Frob}\subseteq X$ of all $n$-dimensional pointed germs of homogeneous semisimple Frobenius manifolds. 
\vskip1mm
For Frobenius manifolds, we have standard choices $\la_o=\frac{d}{2}$ (with $d$ the {\it charge} of the manifold) and $ H_{o,i}=\eta(\der_i,\der_i)^\frac{1}{2}$ for $i=1,\dots,n$. The choice $\la_o=\frac{d}{2}$ further reduces the set of monodromy data, because we automatically have $\La_o=0$. Moreover, the pair $(S_1,S_2)$ must satisfy the relation $S_1^{-1}=S_2^T$, see Remark \ref{mondatfrob}. Hence, due to the constraint \eqref{c2}, the matrices $\mu^{\la_o}$ and $C$ can be reconstructed from the knowledge of $S_1$ only. 
\vskip1mm
The condition $S_1^{-1}=S_2^T$ imposes $\frac{1}{2}(n^2-n)$ constraints.  Notice that we do not need to quotient by the action $(S_1,h)\mapsto h^{-1}S_1h$, since the normalization of the Lam\'e coefficients is fixed.
We have a total amount of free $\frac{1}{2}(n^2-n)$ parameters out of $(\mu^{\la_o},S_1,S_2,C)$. This proves that generic stratum of $X_{\rm Frob}$ has dimension $\frac{1}{2}(n^2-n)$. \endproof
}

\begin{rem}
We underline that the tuple $(\mu^\la,R,S_1,S_2,\La,C)$ of monodromy data actually provides two equivalent systems of ``essential parameters'' classifying germs. For generic germs, one system is $(\mu^\la,[S_1],[S_2])$, the other is $(\mu^\la,[C])$. Here $[S_i]$, with $i=1,2$, denote the conjugacy classes $h^{-1}S_ih$ ($h\in GL(n,\C)$ diagonal), and $[C]$ denotes the equivalence class $h_1Ch_2$ ($h_1,h_2\in GL(n,\C)$ diagonal). In both cases have a total of $n^2-n$ essential parameters.
\end{rem}

\begin{rem}
There is another possible way to deduce that the local isomorphism class of a semisimple homogeneous flat $F$-manifold structure depends on $n^2-n$ moduli, namely by parametrizing the structure with the initial datum $\Tilde\Gm_o$ of the Darboux-Tsarev equations \eqref{ts1}, \eqref{ts2}, \eqref{hgm}. This joint system of differential equations, indeed, is {\it complete} in the sense of Darboux \cite[Livre III, Ch.\,I, pag.\,332]{Dar10}, that is of the form
\[\frac{\der \Tilde\Gm^i_j}{\der u^h}=f_{ijh}(\bm u,\Tilde\Gm),\qquad \Tilde\Gm=\left(\Tilde\Gm^i_j\right)_{i,j=1}^n,\quad i,j,h=1,\dots,n,\quad i\neq j,
\]whose coefficients $f_{ijh}(\bm u,\Tilde\Gm)$ satisfy the compatibility conditions
\[\frac{\der f_{ijh}}{\der u^k}+\sum_{\ell,m}\frac{\der f_{ijh}}{\der \Tilde\Gm^\ell_m}f_{\ell m k}=\frac{\der f_{ijk}}{\der u^h}+\sum_{\ell,m}\frac{\der f_{ijk}}{\der \Tilde\Gm^\ell_m}f_{\ell m h}.
\] A solution $\Tilde\Gm=\left(\Tilde\Gm^i_j\right)_{i,j=1}^n$ is uniquely determined by its initial datum at an arbitrary point $\bm u_o$, see \cite[Livre III, Ch.\,I, Th.\,II]{Dar10}. See also \cite[App.\,1, arXiv version]{AL15}. I thank P.\,Lorenzoni for many useful discussions on this point.
\end{rem}

\section{Applications to LM-CohFTs, $F$-CohFTs,  and open WDVV equations}\label{sec7}
\subsection{Losev-Manin moduli spaces and LM-CohFTs} A $n$-{\it pointed chain of projective lines} $(C; s_0,s_\infty; s_1,\dots, s_n)$ consists of the following data:
\begin{enumerate}
\item a nodal curve $C=C_1\cup\dots\cup C_m$ (over $\C$) whose irreducible components $C_j$ are projective lines;
\item each component $C_j$ is equipped with two marked points $p_j^\pm$, called {\it poles};
\item $C_i$ and $C_j$ intersect only if $|i-j|=1$;
\item $C_i$ and $C_{i+1}$ intersect transversally in $p_i^+=p_{i+1}^-$;
\item $s_0=p_1^-\in C_1$ and $s_\infty=p_m^+\in C_m$ are called {\it white points};
\item $s_1,\dots,s_n\in C\setminus\{p_1^\pm,\dots, p_m^\pm\}$ are called {\it black points}.
\end{enumerate}
A $n$-pointed chain of projective lines is {\it stable} if there is at least one black point on each irreducible components. Notice that black points are allowed to coincide.
\vskip1mm
Two $n$-pointed chains of projective lines $(C; s_0,s_\infty; s_1,\dots, s_n)$ and $(C'; s'_0,s'_\infty; s'_1,\dots, s'_n)$ are {\it isomorphic} if there exists an isomorphism $\phi\colon C\to C'$ such that $\phi(s_j)=s'_j$ for $j=0,1,\dots,n,\infty$.
\vskip1mm
\begin{figure}
\centering
\def\svgscale{.9}
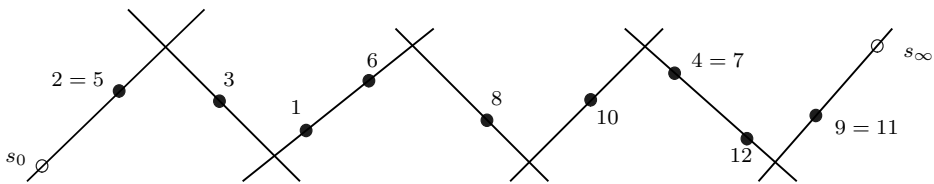
\caption{A stable 12-pointed chain of projective lines.}
\label{grLM}
\end{figure}
\noindent{\bf The spaces $\Ln$. }The {\it Losev-Manin moduli space} $\Ln$, with $n\geq1$, is defined as the fine-moduli space of stable $n$-pointed chains of projective lines \cite{LM00}.
\vskip1mm
The space $\Ln$ is a $(n-1)$-dimensional smooth toric variety (over $\C$): it contains the open dense torus $$L_n=\{(\Pb^1;0,\infty;s_1,\dots,s_n)\}/\text{iso}\cong (\C^*)^{n}/\C^*\cong (\C^*)^{n-1}.$$The space $\Ln$ is the toric variety associated with the convex polytope in $\C^n$ called {\it permutohedron}, defined as the convex hull of the $\frak S_n$-orbit of the point $(1,2,\dots, n)\in\C^n$, see \cite{LM00}. Such a toric variety can be constructed via iteration of blow-ups of $\Pb^{n-1}$. As a first step, blow-up $n$ points $p_1,\dots, p_n$ in general position in $\Pb^{n-1}$. Subsequently, blow-up the strict transforms of the $\frac{1}{2}n(n-1)$ lines passing through the pairs $(p_i,p_j)$ for all $i,j=1,\dots,n$. Continue this blowing-up procedure up to $(n-3)$-dimensional hyperplanes, see \cite[\S 4.3]{Kap93}.
\vskip1mm
Alternatively, $\Ln$ is the toric variety defined by the fan formed by the Weyl chambers of the roots system of type $A_{n-1}$, with $n\geq 2$, \cite{BB11}.
\vskip1mm
The group $\frak S_2\times \frak S_n$ naturally acts on $\Ln$ by permuting white and black points, respectively.
\vskip1mm
The cohomology ring $H^\bullet(\Ln,\Q)$ was studied in \cite{LM00, Man04}. It is algebraic: all odd cohomology groups vanish, and $H^\bullet(\Ln,\Q)$ is isomorphic to the Chow ring $A^\bullet(\Ln,\Q)$, \cite[Th. 2.7.1]{LM00}. See also, \cite{BM14} where the groups $H^\bullet(\Ln,\Q)$ are determined as representation of $\frak S_2\times \frak S_n$.
\vskip1mm
Given $n_1,n_2\geq 1$, we have a natural morphism $\LM{n_1}\times \LM{n_2}\mapsto \LM{n_1+n_2}$, defined by concatenation of white points. Furthermore, each boundary divisor of $\Ln$ is isomorphic to $\LM{n_1}\times \LM{n_2}$ with $n_1+n_2=n$.
\vskip1mm
Let $\DM{0}{n+2}$ be the moduli space of stable $(n+2)$-pointed trees of projective lines. We have a surjective birational morphism $p_n\colon \DM{0}{n+2}\to \Ln$ for any choice of two different labels $i, j$ in $(1, . . . , n + 2)$ (the chosen white points).
\vskip2mm
\noindent{\bf Losev-Manin cohomological field theories. } Let $V_1,V_2$ be two complex vector spaces of finite dimensions. A {\it LM-cohomological field theory} (for short, LM-CohFT), {\it on the pair $(V_1,V_2)$}, is the datum of polylinear maps $$\al_n\colon V_1^*\otimes V_1\otimes V_2^{\otimes n}\to H^\bullet (\Ln,\C),\quad\text{with }n\geq 1,$$such that, for any chosen bases\footnote{In the following paragraphs, if $(e_1,\dots, e_N)$ is a basis of a vector space $V$, then $(e_1^\vee,\dots, e_N^\vee)$ denotes the dual basis of $V^*$.}  $(v_1,\dots, v_{N_1})$ of $V_1$ and $(w_1,\dots, w_{N_2})$ of $V_2$, the following properties are satisfied:
\begin{enumerate}
\item $\al_n$ is $\frak S_n$-covariant w.r.t.\,\,the natural actions of $\frak S_n$ on both $V_2^{\otimes n}$ and $H^\bullet(\Ln,\C)$,
\item for any partition $I\coprod J=\{1,\dots,n\}$ with $|I|=n_1$ and $|J|=n_2$ we have\footnote{Einstein's summation rule over repeated Greek indices is used.}
$${\rm gl}^*\al_{n}(v_i^\vee\otimes v_h\otimes \bigotimes_{i=1}^{n} w_{\rho_i})=\al_{n_1}(v_i^\vee\otimes v_\mu\otimes \bigotimes_{i\in I} w_{\rho_i})\otimes \al_{n_2}(v_\mu^\vee\otimes v_h\otimes \bigotimes_{i\in J} w_{\rho_i}),$$
where $1\leq i,h\leq N_1$ and $1\leq \rho_1,\dots,\rho_n\leq N_2$, and ${\rm gl}\colon \LM{n_1}\times \LM{n_2}\to \LM{n_1+n_2}$ is the gluing map. 
\end{enumerate}

\begin{rem}
The spaces $\LM{n}$ and $\DM{0}{n}$, and their higher genus analogs, are two examples of moduli spaces of {\it weighted }stable pointed curves constructed in \cite{Has03}, corresponding to two different choices of weights. Losev-Manin CohFT's fit in a more general construction developed in \cite{BM09}, in the setting of moduli spaces of curves and maps with weighted stability conditions. We borrow the terminology ``Losev-Manin CohFT'' from \cite{SZ11}. 
\end{rem}

\noindent{\bf Commutativity equations. }Consider two complex vector spaces $V_1,V_2$ of dimension $N_1,N_2$ respectively. Fix a basis $(w_1,\dots, w_{N_2})$ of $V_2$ and let $\bm t:=(t^1,\dots, t^{N_2})$  be the dual coordinates. 
\vskip1mm
The {\it Losev-Manin commutativity equation} for $B\in\C[\![\bm t]\!]\otimes \End(V_1)$ is given by
\beq
\label{commeq}
dB\wedge dB=0.
\eeq In coordinates $\bm t$, equation \eqref{commeq} is equivalent to the commutation relations
\[\left[\frac{\der B}{\der t^i}, \frac{\der B}{\der t^j}\right]=0,\quad i,j=1,\dots,N_2.
\]
Fix a basis $(v_1,\dots, v_{N_1})$ of $V_1$, and let $(v_1^\vee, \dots, v_{N_1}^\vee)$ be the dual basis of $V_1^*$.
\vskip1mm
Given a LM-CohFT $(\al_n)_{n\geq 1}$, define the formal power series $\mc B^i_j\in\C[\![\bm t]\!]$, with $i,j=1,\dots, N_1$, by
\beq\label{Bser}
\mc B^i_j(\bm t):=\sum_{m=1}^\infty\sum_{\rho_1,\dots,\rho_m=1}^{N_2}\frac{t^{\rho_1}\dots t^{\rho_m}}{m!}\int_{\LM{m}}\al_m\left(v_i^\vee\otimes v_j\otimes \bigotimes_{\ell=1}^{m}w_{\rho_\ell}\right).
\eeq
The matrix $\mc B:=(\mc B^i_j)_{i,j=1}^{N_1}$ represents an element of $\C[\![\bm t]\!]\otimes \End(V_1)$ in the basis $v_1,\dots, v_{N_1}$ of $V_1$.

\begin{thm}\label{celm}
The matrix $\mc B$ is a solution of the commutativity equation \eqref{commeq}. Vice-versa, any solution $B$ of \eqref{commeq}, such that $B(0)=0$, has the form \eqref{Bser} for a unique LM-CohFT $(\al_n)_{n\geq 1}$. 
\end{thm}
\proof
This is an equivalent reformulation of \cite[Th.\,3.3.1, Prop.\,3.6.1]{LM00} and \cite[Th.\,5.1.1]{LM04}.
\endproof
\subsection{$F$-cohomological field theories} Let $V$ be a complex vector space of finite dimension $N$. Denote by $\DM{g}{n}$ the Deligne-Mumford moduli space of genus $g$ stable curves with $n$ marked points, defined for $g,n\geq 0$ in the stable regime $2g-2+n>0$.
\vskip1mm
As in \cite{BR18}, an $F$-{\it cohomological field theory} (for short $F$-CohFT) is the datum of 
\begin{itemize}
\item polylinear maps $c_{g,n+1}\colon V^*\otimes V^{\otimes n}\to H^{\rm ev}(\DM{g}{n+1},\C)$, for $2g-1+n>0$,
\item a distinguished vector $e_1\in V$,
\end{itemize}
such that, for any chosen basis $(e_1,\dots, e_N)$ of $V$, the following properties are satisfied:
\vskip2mm
\noindent (1) $c_{g,n+1}$ is $\frak S_n$-covariant w.r.t.\,\,the natural actions of $\frak S_n$ on both $V^*\otimes V^{\otimes n}$ (permutation of the $n$ copies of $V$) and on $H^{\rm ev}(\DM{g}{n+1},\C)$ (permutation of the last $n$ marked points);
\vskip1mm
\noindent (2) $\pi^*c_{g,n+1}(e_{\rho_0}^\vee\otimes \bigotimes_{i=1}^n e_{\rho_i})=c_{g,n+2}(e_{\rho_0}^\vee\otimes \bigotimes_{i=1}^n e_{\rho_i}\otimes e_1)$, for $1\leq \rho_1,\dots,\rho_n\leq N$, where $\pi\colon \DM{g}{n+2}\to \DM{g}{n+1}$ is the map forgetting the last marked point;
\vskip1mm
\noindent (3) $c_{0,3}(e_\al^\vee\otimes e_\bt\otimes e_1)=\dl^\al_\bt$, for $1\leq \al,\bt\leq N$;
\vskip1mm
\noindent (4) for any partition $I\coprod J=\{1,\dots, n\}$ with $|I|=n_1$ and $|J|=n_2$, we have\footnote{Einstein's summation rule over repeated Greek indices is used.} $${\rm gl}^*c_{g_1+g_2,n_1+n_2+1}(e_{\rho_0}^\vee\otimes \bigotimes_{i=1}^n e_{\rho_i})=c_{g_1,n_1+2}(e_{\rho_0}^\vee\otimes \bigotimes_{i\in I} e_{\rho_i}\otimes e_\mu)\otimes c_{g_2,n_2+1}(e_\mu^\vee\otimes\bigotimes_{j\in J}e_{\rho_j}),$$ for $1\leq \rho_1,\dots, \rho_n\leq N$, and ${\rm gl}\colon \DM{g_1}{n_1+2}\times \DM{g_2}{n_2+1}\to \DM{g_1+g_2}{n_1+n_2+1}$ is the corresponding gluing map.
\vskip2mm
The {\it genus 0 sector} (or {\it tree-level}) of a given $F$-CohFT is the datum of the maps $(c_{0,n})_{n\geq 2}$ and the distinguished vector $e_1\in V$ only.
\vskip1mm
Given a tree-level $F$-CohFT, fix a basis $(e_1,\dots, e_N)$ of $V$, and denote by $\bm t:=(t^1,\dots, t^N)$ the dual coordinates. Define the formal power series $F^\al\in\C[\![\bm t]\!]$, for $\al=1,\dots, N$, by
\beq
\label{potF}
F^\al(\bm t):=\sum_{n=2}^\infty\sum_{\rho_1,\dots,\rho_n=1}^N\frac{t^{\rho_1}\dots t^{\rho_n}}{n!}\int_{\DM{0}{n+1}}c_{0,n+1}\left(e_\al^\vee\otimes\bigotimes_{i=1}^n e_{\rho_i}\right).
\eeq

\begin{thm}\label{fmfco}
The functions $F^\al(\bm t)$ are solution of the oriented associativity equations 
\begin{align}
\label{oae1}
\frac{\der^2 F^\al}{\der t^1\der t^\bt}&=\dl^\al_\bt,\quad&\al,\bt=1,\dots,N,\\
\label{oae2}
\frac{\der^2 F^\al}{\der t^\mu\der t^\bt}\frac{\der^2 F^\mu}{\der t^\gm\der t^\dl}&=\frac{\der^2 F^\al}{\der t^\mu\der t^\gm}\frac{\der^2 F^\mu}{\der t^\bt\der t^\dl},\quad&\al,\bt,\gm,\dl=1,\dots,N,
\end{align}
and thus define a formal flat $F$-manifold structure on $V$ with unit $e_1$.
\vskip1mm
Vice-versa, any solution $(F^1,\dots, F^N)$ of \eqref{oae1}-\eqref{oae2}, with $F^\al(0)=0$ and $\der_\bt F^\al(0)=0$ for all $\al,\bt=1,\dots, N$, is of the form \eqref{potF} for a unique tree-level $F$-CohFT $(c_{0,n})_{n\geq 2}$.
\end{thm}
\proof
The first part of the statement follows from a simple computation, invoking properties (1)-(4) above. 
For a proof of the second part of the statement, see Appendix \ref{app}.
\endproof

\begin{rem}
The first part of Theorem \ref{fmfco} was already formulated in \cite[Sec.\,2.1]{BR18}. The second part of the statement, however, was not addressed in {\it loc.\,cit.}.
\end{rem}

\subsection{From tree-level $F$-CohFT to LM-CohFT, and vice-versa}Given a tree-level $F$-CohFT on $V$, a LM-CohFT is naturally defined on the pair $(V_1,V_2)=(V,V)$. For any $n\geq 1$ define
\[\al_n:=(p_{n})_*\circ c_{0,n+2}\colon V^*\otimes V^{\otimes (n+1)}\to H^\bullet (\Ln,\C),
\]where $p_n\colon \DM{0}{n+2}\to\Ln$ is the surjective birational morphism defined by a choice of two white points.

\begin{prop}\label{FtoLM}
The polylinear maps $(\al_n)_{n\geq 1}$ define a LM-CohFT on $(V,V)$.
\end{prop}
\proof
The $\frak S_n$-covariance of $\al_n$ follows from the $\frak S_{n+1}$-covariance of $c_{0,n+2}$. For $n_1+n_2=n$, we have the following commutative diagram
\[
\xymatrix{
\DM{0}{n_1+2}\times \DM{0}{n_2+2}\ar[r]^{\qquad\rm gl}\ar[d]_{p_{n_1}\times p_{n_2}}&\DM{0}{n+2}\ar[d]^{p_n}\\
\LM{n_1}\times \LM{n_2}\ar[r]^{\qquad\widetilde{\rm gl}}&\LM{n}
}
\]with proper vertical arrows and local complete intersections as horizontal arrows. The gluing property of $\al_n$ then follows from the gluing property of $c_{0,n+2}$ and the excess intersection formula \cite[Prop. 6.6 and Prop. 17.4.1]{Ful98}. Notice that the excess bundle $\mathbb E$ has rank 0 (both $\rm gl$ and $\widetilde{\rm gl}$ have codimension 1), hence
\[\widetilde{\rm gl}^*(p_n)_*\,x=(p_{n_1}\times p_{n_2})_*\,{\rm gl}^*x,
\]for all $x\in H^{\rm ev}(\DM{0}{n+2},\C)\cong A^\bullet(\DM{0}{n+2})_{\C}$.
\endproof

Vice-versa, given a LM-CohFT $(\al_n)_{n\geq 1}$ on $(V_1,V_2)$ we can reconstruct a formal flat $F$-manifold, provided that
\begin{itemize}
\item $\dim V_1=\dim V_2=N$,
\item we are given an extra amount of data, namely a {\it primitive vector}.
\end{itemize}
\begin{defn}
Let $B\in\C[\![\bm t]\!]\otimes \End(V_1)$ be a solution of commutativity equations \eqref{commeq}. A vector $h\in V_1$ is {\it primitive for $B$} if the vectors 
\[\left.\frac{\der B}{\der t^1}\right|_{\bm t=0}\cdot h,\quad\dots,\quad \left.\frac{\der B}{\der t^N}\right|_{\bm t=0}\cdot h
\]define a basis of $V_1$. Equivalently, $h$ is primitive if, for any chosen basis $(v_1,\dots, v_N)$ of $V_1$, we have 
\[
\det\left(\left.\frac{\der B^i_\mu}{\der t^k}\right|_0h^\mu\right)_{i,k=1}^N\neq 0,\quad \text{where }h=h^\mu v_\mu.
\]
\end{defn}

If $B$ admits a primitive vector $h$, we can identify $V_1$ and $V_2$ via the isomorphism $w_k\mapsto \left.\frac{\der B}{\der t^k}\right|_{0}h$, for $k=1,\dots, N$. Under such an identification, $\bm t$ can be thought as coordinates on $V_1\cong V_2$.

\begin{prop}[{\cite[Prop. 5.3.3]{LM04}}]\label{bpri}
If $B$ admits a primitive vector $h$, then there exists a formal flat $F$-manifold structure on $V_1\cong V_2$ with flat identity $h$. The oriented associativity potentials $\bm F=(F^1,\dots, F^N)$ satisfy $B^\al_\bt=\frac{\der F^\al}{\der t^\bt}$.
\end{prop}
\proof

Let $(v_1,\dots, v_N)$ be a basis of $V_1$ with $v_1=h$. By assumption, we have $\det\left(\left.\frac{\der B^i_1}{\der t^k}\right|_0\right)\neq 0$. Hence, up to change of basis $(w_1,\dots, w_N)$ of $V_2$, we can assume that in the coordinates $\bm t$ we have
\[\frac{\der B^i_1}{\der t^k}(\bm t)=\dl^i_k\quad \Rightarrow \quad B^i_1(\bm t)=t^i+c.
\]Consider the $\C[\![\bm t]\!]\otimes V_1$-valued differential form $dB\wedge d(Bh)$. In the bases $(w_i)_{i=1}^N$ and $(v_i)_{i=1}^N$ chosen as above, it has components 
\[[dB\wedge d(Bh)]^i=\left(\frac{\der B^i_\la}{\der t^\nu}dt^\nu\right)\wedge\left(\dl^\la_\mu dt^\mu\right)=
\frac{\der B^i_\mu}{\der t^\nu}\,dt^\nu\wedge dt^\mu.
\]On the other hand, $dB\wedge d(Bh)=(dB\wedge dB)h=0$,  since $h$ is $\bm t$-independent and $B$ a solution of \eqref{commeq}. Hence, 
\[\frac{\der B^i_\mu}{\der t^\nu}-\frac{\der B^i_\nu}{\der t^\mu}=0,\quad \mu,\nu=1,\dots, N.
\]This implies the existence of $F^i\in\C[\![\bm t]\!]$ such that $B^i_j=\der_jF^i$.
\endproof

\begin{thm}
The following notions are equivalent:
\begin{enumerate}
\item formal flat $F$-manifold,
\item tree-level $F$-cohomological field theory,
\item $LM$-cohomological field theory with primitive element.
\end{enumerate}
\end{thm}
\proof
It follows from Theorems \ref{celm}, \ref{fmfco} and Propositions \ref{FtoLM}, \ref{bpri}.
\endproof

\begin{rem}
Black marked points of stable $n$-pointed chains of projective lines are allowed to coincide. It would be tempting to compare these coincidences of black points with the coalescence phenomenon at irregular singularities of ordinary differential equations studied in \cite{CG18,CDG1}. Any contingent relation deserves further investigations. I thank Yu.I.\,Manin for pointing out such an analogy in a private communication.
\end{rem}

\subsection{Homogeneous $F$-CohFTs} A $F$-CohFT $(c_{g,n+1})_{g,n}$ is said to be {\it homogeneous} if
\begin{enumerate}
\item the vector spaces $V$ and $V^*$ are graded, with homogeneous bases $(e_1,\dots, e_N)$ and $(e_1^\vee,\dots, e_N^\vee)$,
\bea
&&\deg e_\al=-\deg e_\al^\vee=q_\al,\quad \al=1,\dots,N\qquad\deg e=0;
\eea
\item there exist $r^1,\dots, r^N,\gm\in\C$ such that
\begin{multline}\label{hg}
{\rm Deg}\,c_{g,n+1}\left(e_{\al_0}^\vee\otimes \bigotimes_{i=1}^n e_{\al_i}\right)+\pi_*\,c_{g,n+2}\left(e_{\al_0}^\vee\otimes \bigotimes_{i=1}^n e_{\al_i}\otimes r^\la e_\la\right)\\
=\left(\sum_{i=1}^nq_{a_i}-q_{\al_0}+\gm g\right)c_{g,n+1}\left(e_{\al_0}^\vee\otimes \bigotimes_{i=1}^n e_{\al_i}\right),
\end{multline}
\end{enumerate}
where ${\rm Deg}\colon H^\bullet(\DM{g}{n})\to H^\bullet(\DM{g}{n})$ rescale a $k$-th degree class by a factor $\frac{k}{2}$, and $\pi\colon \DM{g}{n+2}\to\DM{g}{n+1}$ is the morphism forgetting the last marked point.

\begin{prop}
If a tree-level $F$-CohFT is homogeneous then the associated formal flat $F$-manifold with potentials \eqref{potF} is homogeneous, with the Euler vector field 
\[E=\sum_{\al=1}^N((1-q_\al)t^\al+r^\al)\frac{\der}{\der t^\al}.
\]
\end{prop}
\proof
A simple computation shows that equations \eqref{hg}, specialized at $g=0$, imply equations \eqref{hEF} for the potentials \eqref{potF}.
\endproof

The following result follows from Theorem \ref{thconv}, and Abel Lemma.

\begin{thm}
Let $(H,\bm \Phi)$ be a formal semisimple flat $F$-manifold over $\C$, with $\dim_\C H=N$. Let $(c_{0,n+1})_{n\geq 2}$ and $(\al_n)_{n\geq 1}$ be the underlying tree-level $F$-CohFT and $LM$-CohFT, respectively. If $(H,\bm \Phi)$ is admissible, then there exist real positive constants $m,k_1,\dots, k_N\in\R_+$ such that 
\[\left|\int_{\DM{0}{|\bm n|+1}}c_{0,|\bm n|+1}\left(\Dl_\bt^\vee\otimes \bigotimes_{j=1}^N\Dl_j^{\otimes n_j}\right)\right|\leq m\,\bm n!\prod_{j=1}^Nk_j^{n_j},\quad \bm n\in\N^N,\quad \bt=1,\dots, N,
\]
\[\left|\int_{\LM{|\bm n|}}\al_{|\bm n|}\left(\Dl_\bt^\vee\otimes \Dl_\gm\otimes \bigotimes_{j=1}^N\Dl_j^{\otimes n_j}\right)\right|\leq m\,\bm n!\prod_{j=1}^Nk_j^{n_j},\quad \bm n\in\N^N,\quad \bt,\gm=1,\dots, N,
\]where we set $\bm n!:=\prod_{j=1}^N n_j$, and $|\bm n|:=\sum_{j=1}^Nn_j$.\qed
\end{thm}

\subsection{Open WDVV equations} Let $k$ be a commutative $\Q$-algebra. Consider a formal Frobenius manifold over $k$, with Euler vector field $E$, defined by the solution $F\in k[\![t^1,\dots, t^n]\!]$ of the WDVV equations: 
\begin{align}
\label{orwdvv1}
&\frac{\der^3F}{\der t^\al\der t^\bt\der t^\mu}\eta^{\mu\nu}\frac{\der^3F}{\der t^\nu\der t^\gm\der t^\dl}=\frac{\der^3F}{\der t^\dl\der t^\bt\der t^\mu}\eta^{\mu\nu}\frac{\der^3F}{\der t^\nu\der t^\gm\der t^\al},%
\\
\label{orwdvv2}
&\frac{\der^3F}{\der t^1\der t^\al\der t^\bt}=\eta_{\al\bt}={\rm const.,}\quad \eta=(\eta_{\al\bt})_{\al,\bt},\quad \eta^{-1}=(\eta^{\al\bt})_{\al,\bt},%
\\
\label{orwdvv3}
&E^\nu\frac{\der F}{\der t^\nu}=(3-d)F+Q(\bm t),\quad %
E^\nu=(1-q^\nu)t^\nu+r^\nu,
\end{align}where $\al,\bt,\gm\in\{1,\dots,n\}$, $d\in k$ is the conformal dimension (or charge) of the Frobenius manifold, $q^\nu,r^\nu\in k$, and $Q(\bm t)\in k[\bm t]$ is a quadratic polynomial in $\bm t$.
\vskip1mm
The {\it open WDVV equations} (OWDVV) are the following overdetermined system of PDEs for $F^o\in k[\![t^1,\dots, t^n,s]\!]$:
\begin{align}
\label{owdvv1}
&\frac{\der^3F}{\der t^\al\der t^\bt\der t^\mu}\eta^{\mu\nu}\frac{\der^2 F^o}{\der t^\nu\der t^\gm}+\frac{\der^2 F^o}{\der t^\al\der t^\bt}\frac{\der^2 F^o}{\der s\der t^\gm}=
\frac{\der^3F}{\der t^\gm\der t^\bt\der t^\mu}\eta^{\mu\nu}\frac{\der^2 F^o}{\der t^\nu\der t^\al}+\frac{\der^2 F^o}{\der t^\gm\der t^\bt}\frac{\der^2 F^o}{\der s\der t^\al},\\
\label{owdvv2}
&\frac{\der^3F}{\der t^\al\der t^\bt\der t^\mu}\eta^{\mu\nu}\frac{\der^2 F^o}{\der t^\nu\der s}+\frac{\der^2 F^o}{\der t^\al\der t^\bt}\frac{\der^2 F^o}{\der s^2}=\frac{\der^2 F^o}{\der s\der t^\bt}\frac{\der^2 F^o}{\der s\der t^\al},\\
\label{owdvv3}
&\frac{\der^2 F^o}{\der t^1\der t^\al}=0,\qquad \frac{\der^2 F^o}{\der t^1\der s}=1,\\
\label{owdvv4}
&E^\nu\frac{\der F^o}{\der t^\nu}+\left(\frac{1-d}{2}s+r^{n+1}\right)\frac{\der F^o}{\der s}=\frac{3-d}{2}F^o+L(\bm t),
\end{align}
where $\al,\bt,\gm\in\{1,\dots, n\}$, $r^{n+1}\in k$, and $L(\bm t)\in k[\bm t]$ is a linear polynomial in $\bm t$.
\vskip1mm
The OWDVV equations first appeared in \cite[Th. 2.7]{HS12}, in the context of open Gromov-Witten theory. These equations subsequently appeared in \cite{PST14,BCT18,BCT19}: although not explicitly mentioned in {\it loc.\,cit.}, the OWDVV equations follow from the open Topological Recursion Relations equations, see \cite[Th.\,1.5]{PST14}, \cite[Th.\,4.1]{BCT18}, \cite[Lem.\,3.6]{BCT19}, \cite[Sec.\,4]{Bur18}, \cite[Sec.\,1]{BB19}. The OWDVV equations play a central role in the general theory of relative quantum cohomology developed in \cite{ST19}.

\begin{prop}[P.\,Rossi, \cite{BB19}]\label{proprossi}The following conditions are equivalent:
\begin{enumerate}
\item $(F,F^o)$ is a solution of WDVV and OWDVV equations \eqref{orwdvv1}-\eqref{owdvv4},
\item $\left(\frac{\der F}{\der t^\mu}\eta^{\mu1},\dots,\frac{\der F}{\der t^\mu}\eta^{\mu n}, F^o\right)$ is a solution of the oriented associativity equations \eqref{wdvv1}-\eqref{wdvv2} in the coordinates $(\bm t,s)$, and the corresponding formal flat $F$-manifold is homogenous. 
\end{enumerate}
\end{prop}
\proof
The claim follows by a direct check.
\endproof
We will refer to the formal flat $F$-structure of point (2) of Proposition \ref{proprossi} as the formal flat $F$-manifold {\it underlying the pair} $(F,F^o)$.
As a corollary of Theorem \ref{thconv}, we deduce the following result.

\begin{thm}
Let $F\in\C[\![\bm t]\!]$, $F^o\in\C[\![\bm t, s]\!]$ be solutions of the WDVV and OWDVV equations. If the underlying formal flat $F$-manifold is semisimple, and it is not strictly doubly resonant, %
 then both $F$ and $F^o$ are convergent.\qed
\end{thm}

\begin{rem}
According to a conjecture of B.\,Dubrovin, the monodromy data of the quantum cohomology of a smooth projective variety encode information about the derived category $\mc D^b(X)$, see \cite{Dub98,CDG18,Cot20a,Cot20c}. It would be interesting to look for analog relations starting from the monodromy data, as defined here, of flat $F$-structures given by relative quantum cohomologies of \cite{ST19}. This will be addressed in a future project of the author.
\end{rem}

\appendix

\section{Proof of Theorem \ref{eulerunique}}\label{theuelrid}

Let $M$ be an analytic homogeneous semisimple flat $F$-manifold, and let $E_1,E_2\in\Gm(TM)$ be two Euler vector fields. 
\begin{lem}
We have $[E_1,E_2]=E_1-E_2$. 
\end{lem}
\proof By Proposition \ref{canuE}, we can choose canonical coordinates so that $E_1=\sum_ju^j\der_j$ and $E_2=\sum_j(u^j+c^j)\der_j$.
 We have
\[\pushQED{\qed} 
[E_1,E_2]=\sum_{j}(E_1^h\der_hE_2^j-E_2^h\der_hE_1^j)\der_j=-\sum_jc^j\der_j=E_1-E_2.\qedhere
\popQED
\]
\begin{lem}
We have $\nabla E_1=\nabla E_2$.
\end{lem}
\proof
Since $\nabla$ is torsionless, we have $\nabla_{E_1}E_2=\nabla_{E_2}E_1+[E_1,E_2]$. For arbitrary vector field $X$, we have
\[
\pushQED{\qed} 
\underbrace{\nabla_X\nabla_{E_1}E_2}_{0}=\underbrace{\nabla_X\nabla_{E_2}E_1}_{0}+\nabla_X[E_1,E_2].\qedhere
\popQED
\]
\proof[Proof of Theorem \ref{eulerunique}] Introduce the operators $\mc U_1(X):=E_1\circ X$ and $\mc U_2(X):=E_2\circ X$, and 
\[\mu(X):=X-\nabla_XE_1=X-\nabla_XE_2,\quad X\in\Gm(TM).
\]Choose canonical coordinates so that $E_1=\sum_ju^j\der_j$ and $E_2=\sum_j(u^j+c^j)\der_j$. Set $\widetilde\Psi^i_\al=\frac{\der u^i}{\der t^\al}$. The matrix $\widetilde\Psi$ diagonalizes both $\mc U_1$ and $\mc U_2$:
\[U_1=\widetilde\Psi\,\mc U_1\,\widetilde\Psi^{-1}={\rm diag}(u^1,\dots, u^n),
\]
\[U_2=\widetilde\Psi\,\mc U_2\,\widetilde\Psi^{-1}={\rm diag}(u^1+c^1,\dots, u^n+c^n).
\]
Set $\widetilde V:=\widetilde\Psi\,\mu\,\widetilde\Psi^{-1}$, $\widetilde V_i:=\der_i\widetilde\Psi\cdot\widetilde\Psi^{-1}$, and also introduce an off-diagonal matrix $\widetilde\Gm=(\widetilde\Gm^i_j)$ by
\[
\widetilde\Gm^i_j:=-(\widetilde V_i)^i_j,\quad i\neq j.
\] The matrix $\widetilde\Gm$ is a solution of the Darboux-Tsarev equations. Moreover, we have
\beq
\label{eq3}
\widetilde V=\widetilde V'+[\widetilde\Gm,U_1],\quad\text{and also}\quad \widetilde V=\widetilde V'+[\widetilde\Gm,U_2],
\eeq
where $V'$ denotes the diagonal part of $V$. This follows from Propositions \ref{pids} and \ref{DTeqs}. 
\vskip2mm
Let $j\neq k$, and take the $(j,k)$ entry of both equations \eqref{eq3}. We have
\[\widetilde V^j_k(\bm u)=\widetilde{\Gm}^j_k(\bm u)(u^k-u^j)=\widetilde{\Gm}^j_k(\bm u)(u^k-u^j)+\widetilde{\Gm}^j_k(\bm u)(c^k-c^j).
\]Hence, we have
\[\widetilde{\Gm}^j_k(\bm u)(c^k-c^j)=0,\quad \text{for any }j\neq k.
\]This means that 
\[
c^j\neq c^k\quad \Longrightarrow\quad \widetilde{\Gm}^j_k=\widetilde{\Gm}^k_j\equiv 0.
\]
Introduce the partition $\coprod_{r=1}^NI_r=\{1,\dots,n\}$ s.t.\,\,$c^i=c^j$ only if $i,j$ are in a same block, i.e. $i,j\in I_r$ for some $r$.
\vskip1mm
It follows that $\widetilde\Gm^i_j=0$ unless $i,j\in I_r$. Take $i,j\in I_r$ and $k\notin I_r$. We have
\[\der_k\widetilde\Gm^i_j=-\widetilde\Gm^i_j\widetilde\Gm^k_j+\widetilde\Gm^i_j\widetilde\Gm^j_k+\widetilde\Gm^i_k\widetilde\Gm^k_j=0.
\]We have proved that
\begin{itemize}
\item the function $\widetilde\Gm^i_j$ is not identically zero only if the indices $i,j$ are in the same block $I_r$;
\item the function $\widetilde\Gm^i_j$ only depends on coordinates $u^k$ with $i,j,k$ in the same block $I_r$.
\end{itemize}
It follows that all the matrices $\widetilde\Gm, \widetilde V, \widetilde V_i, \widetilde \Psi$ admits a direct sum decomposition
\[\widetilde\Gm=\bigoplus_{r=1}^N \widetilde\Gm^{(r)},\quad \widetilde V=\bigoplus_{r=1}^N\widetilde V^{(r)},\quad \widetilde V_i=\bigoplus_{r=1}^N\widetilde V_i^{(r)},\quad \widetilde \Psi=\bigoplus_{r=1}^N\widetilde \Psi^{(r)},
\]
and each summand $\widetilde\Gm^{(r)}, \widetilde V^{(r)}, \widetilde V_i^{(r)}, \widetilde \Psi^{(r)}$ only depends on canonical coordinates $u^k$ with $k\in I_r$. The original flat $F$-manifold $M$ locally decomposes into $N$ corresponding pieces:
\[M\stackrel{\text{loc.}}{\cong} \bigoplus_{j=1}^N M^{(j)}.
\]The flat $F$-manifold $M$ is irreducible if and only if $N=1$. This completes the proof.
\endproof 

\section{Proof of Theorem \ref{fmfco}}\label{app}
\begin{figure}
\centering
\def\svgscale{.9}
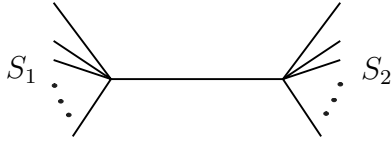
\caption{Isomorphism classes of one edges $n$-trees are parametrized by \sutp s $\si=\{S_1,S_2\}$ of $\{1,\dots,n\}$.}
\label{gr1}
\end{figure}

In order to complete the proof of Theorem \ref{fmfco}, we need to recall some preliminary known results on the (co)homology groups $H_\bullet(\DM{0}{n},\C)$ and $H^\bullet(\DM{0}{n},\C)$.
\vskip1mm
\noindent{\bf Graphs. }In what follows, a {\it graph} $\tau$ is an ordered family $(V_\tau, H_\tau, {\der}_\tau, j_\tau)$ where
\begin{itemize} 
\item $V_\tau$ is a finite set of {\it vertices},
\item $H_\tau$ is a finite set of {\it half-edges}, equipped with a vertex assignment function ${\der_\tau}\colon H_\tau\to V_\tau$, and an involution $j_\tau\colon H_\tau\to H_\tau$.
\end{itemize}
The set $E_\tau$ of 2-cycles of $j_\tau$ is the set of {\it edges} of $\tau$. The set $S_\tau$ of fixed points of $j_\tau$ is the set of {\it tails} of $\tau$. The datum of $j_\tau$ is thus equivalent to the datum of $E_\tau$ and $S_\tau$. For each vertex $v\in V_\tau$ define the set $H_\tau(v)$ of {\it half-edges attached at $v$} by $$H_\tau(v):=\der_\tau^{-1}(v),$$ the set $E_\tau(v)$ of {\it edges attached to }$v$ by $$E_\tau(v):=\{\{f_1,f_2\}\in E_\tau\colon \der_\tau(f_1)=v\text{ or }\der_\tau(f_2)=v\},$$
and the set $S_\tau(v)$ of {\it tails attached to $v$} by
\[S_\tau(v):=\{f\in S_\tau\colon \der_\tau f=v\}.
\] We clearly have a partition $H_\tau(v)\cong E_\tau(v)\coprod S_\tau(v)$.
\vskip1mm
An isomorphism $\tau_1\to \tau_2$ of graphs is the datum of two bijections $V_{\tau_1}\to V_{\tau_2}$ and $H_{\tau_1}\to H_{\tau_2}$ compatible with $\der$ and $j$.
\vskip1mm
A graph $\tau$ is conveniently identified with its associated topological space $\|\tau\|$. Vertices of $\tau$ are identified with $|V_\tau|$ distinct points $\{p(v)\}_{v\in V_\tau}$ on the curve $\mc C:=\{(t,t^2,t^3)\colon t\in\R\}\subseteq \R^3$, an edge $\{f_1,f_2\}\in E_\tau$ is identified with the segment\footnote{By Vandermode determinant any four points $p(v_1),\dots, p(v_4)$ are not coplanar, and they are vertices of a tetrahedron. This argument shows that segments representing edges intersect only at the appropriate vertices. } joining the points $p(\der_\tau(f_1))$ and $p(\der_\tau(f_2))$, tails at $v$ are identified with a star of $|S_\tau(v)|$ small segments originating from $p(v)$, intersecting neither edges nor other tails at other vertices. The space $\|\tau\|$ is the union of all these vertices and segments, equipped with the topology induced from $\R^3$. The graph $\tau$ is a {\it tree} if $\|\tau\|$ is connected and $H_1(\|\tau\|,\Z)=0$. A tree with $n$ tails, will be called an $n$-{\it tree}.
\vskip1mm
\noindent{\bf Dual stable graphs. }To each point $[C,(\bm x)]\in\DM{0}{n}$, we can attach a {\it dual stable graph} $\tau$ as follows:
\begin{enumerate}
\item the vertices of $\tau$ are in 1-1 correspondence with the irreducible components of $C$,
\item each node of $C$ is replaced by an edge connecting the vertices corresponding to the two sides of the node,
\item for each $i=1,\dots,n$ attach an \emph{tail} with label $i$ to the vertex corresponding to the irreducible component containing $x_i$.
\end{enumerate}
The resulting graph $\tau$ is always a $n$-tree satisfying the following \emph{stability condition}: each vertex has valence $|H_\tau(v)|$ at least 3. We say that $[C,(\bm x)]$ has \emph{combinatorial type} $\tau$.
\vskip1mm
Vice-versa, for any stable $n$-tree $\tau$, there exists a locally closed irreducible subscheme $D(\tau)\subseteq \DM{0}{n}$ parametrizing curves of combinatorial type $\tau$. The stratum $D(\tau)$ uniquely identifies the isomorphism class of $\tau$, and its codimension equals the number of edges $|E_\tau|$.\vskip1mm
For example, the $n$-tree with one vertex corresponds to the open stratum $\mc M_{0,n}$. The strata of codimension one are labelled by isomorphism classes of one-edge stable $n$-trees $\sigma$. Each such class can be identified with a {\it stable unordered 2-partition} of the set $\{1,\dots,n\}$. This consists of a set $\si=\{S_1,S_2\}$ such that $\{1,\dots, n\}=S_1\coprod S_2$, and $|S_i|\geq 2$ for $i=1,2$. See Figure \ref{gr1}.
\vskip1mm
Given $(i,j,k,l)\in\{1,\dots, n\}^4$ and a \sutp\, $\si$, we write $ij\si kl$ if $i,j$ and $k,l$ belong two the two different elements of $\si$.\vskip1mm
\noindent{\bf Keel's Theorem. }Introduce commuting indeterminates $D_\si$, indexed by stable unordered 2-partitions $\si$ of $\{1,\dots, n\}$. Consider the ideal $I_n\subseteq\C[(D_\si)_\si]$ generated as follows:
\begin{enumerate}
\item for each $(i,j,k,l)\in\{1,\dots, n\}^4$ set
\beq\label{krel}
R_{ijkl}:=\sum_{ij\si kl}D_\si-\sum_{kj\tau il}D_\tau\in I_n;
\eeq
\item if $\si$ and $\tau$ are such that $ij\si kl$ and $ij\tau kl$ for some $(i,j,k,l)\in\{1,\dots,n\}^4$, then set
\beq
D_\si D_\tau\in I_n.
\eeq
\end{enumerate}
\begin{thm}\label{thkeel}
We have an isomorphism of rings
\beq\label{kiso}
\C[(D_\si)_\si]/I_n\to H^\bullet(\DM{0}{n},\C)\cong A^\bullet(\DM{0}{n})_\C
\eeq
defined by
\[D_\si\mapsto \text{dual of the closed cycle }\overline{D(\si)}.
\]In particular, all odd cohomology groups vanish.
\end{thm}
\vskip1mm
\noindent{\bf Good monomials. }Consider a stable $n$-tree $\tau$. Any edge $e\in E_\tau$ defines a \sutp\, $\sigma(e)$ of $\{1,\dots,n\}$: by cutting $e$, we obtain two trees, whose tails (halves of $e$ excluded) form the two parts of $\si(e)$.
\vskip1mm
For each stable $n$-tree $\tau$, define the monomial $m(\tau):=\prod_{e\in E_\tau}D_{\si(e)}$. This is a monomial in $\C[(D_\si)_\si]$ of degree $|E_\tau|$. Monomials of this form are called {\it good monomials}.
\vskip1mm
Under Keel isomorphism \eqref{kiso}, $m(\tau)$ is the dual of the class $[\overline{D(\tau)}]\in H_\bullet(\DM{0}{n},\C)$. This follows from the fact that boundary components intersect transversally.
\begin{thm}[{\cite[Ch.III, \S3.6]{man}}]\label{mtspan}
Good monomials modulo $I_n$ span $\C[(D_\si)_\si]/I_n$. Equivalently, the classes $[\overline{D(\tau)}]$ span $H_\bullet(\DM{0}{n},\C)$.
\end{thm}
\vskip1mm
\noindent{\bf Manin's relations in higher codimensions. }We need information about linear relations among all good monomials of fixed degree, generalizing Keel's relations \eqref{krel}.
\vskip1mm
Let $n\geq 4$, and $\tau$ to be an $n$-tree. Given 
\begin{enumerate}
\item $v\in V_\tau$ with valence $|H_\tau(v)|\geq 4$,
\item $(i,j,k,l)\in H_\tau(v)^4$ pairwise distinct half-edges,
\end{enumerate}
set $T:=H_\tau(v)\setminus\{i,j,k,l\}$. For any ordered 2-partition $\al=(T_1,T_2)$ of $T$ (the case $T_i=\emptyset$ is allowed), we can define two trees $\tau'(\al)$ and $\tau''(\al)$, with $|E_\tau|+1$ edges each.
\vskip1mm
The tree $\tau'(\al)$ is defined by replacing the vertex $v$ with a new edge $e$, at whose vertices we have half-edges $\{i,j\}\cup T_1$ and $\{k,l\}\cup T_2$, respectively.

\begin{figure}
\centering
\def\svgscale{.9}
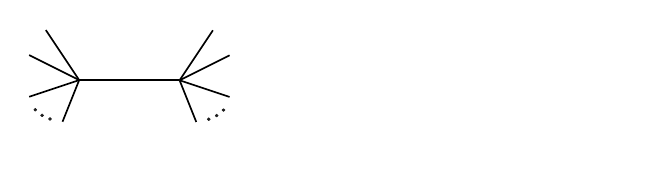
\caption{}
\label{gr2}
\end{figure}

The tree $\tau''(\al)$ is defined by replacing the vertex $v$ with a new edge $e$, at whose vertices we have half-edges $\{k,j\}\cup T_1$ and $\{i,l\}\cup T_2$, respectively. See Figure \ref{gr2}.
\vskip1mm
For each system $(\tau,v,i,j,k,l)$ as above, define the polynomial
\[R(\tau,v,i,j,k,l):=\sum_\al m(\tau'(\al))-m(\tau''(\al))\,\in \C[(D_\si)_\si].
\]
\begin{thm}\label{mrel}
We have $R(\tau,v,i,j,k,l)\in I_n$. Moreover, all linear relations modulo $I_n$ between good polynomials of degree $r+1$ are spanned by all the relations $R(\tau,v,i,j,k,l)$ with $|E_\tau|=r$.
\end{thm}
For a proof, see \cite[Ch.\,III, Prop.\,4.7.1, Th.\,4.8]{man}.
\vskip1mm
\noindent{\bf Proof of Theorem \ref{fmfco}. }We are now able to complete the proof. Given potentials 
\beq\label{potf1}
F^\al(\bm t)=\sum_{n\geq 2}\sum_{\rho_1,\dots, \rho_n=1}^N\frac{t^{\rho_1}\dots t^{\rho_n}}{n!}c^\al_{\rho_1\dots\rho_n},\quad c^\al_{\rho_1\dots\rho_n}\in\C,\quad \al=1,\dots,N,
\eeq
equations \eqref{oae1} and \eqref{oae2} are
\bean
\label{eqc1}
&c^\al_{1\bt}=\dl^\al_\bt,\quad c^\al_{1\bt\rho_1\dots\rho_n}=0\quad \text{for } n>0,\\
\label{eqc2}
&c^\al_{\mu\bt\rho_1\dots\rho_n}c^\mu_{\gm\dl\tau_1\dots\tau_k}= c^\al_{\mu\gm\rho_1\dots\rho_n}c^\mu_{\bt\dl\tau_1\dots\tau_k}.
\eean
Notice that the lower indices of the coefficients $c$'s can be arbitrarily permuted, i.e. $c^\al_{\rho_1\dots\rho_n}$ is uniquely identified by $\al$ and the set $\{\rho_1,\dots,\rho_n\}$. We need to prove that the potentials $F^\al$ are of the form \eqref{potF} for a unique existing tree-level $F$-CohFT $(c_{0,n+1})_{n\geq 2}$. We first prove the uniqueness, and hence the existence of such an $F$-CohFT.
\vskip1mm 
 {\bf Uniqueness. }Assume there exists a tree-level $F$-CohFT $(c_{0,n+1})_{n\geq 2}$ such that\footnote{For later notational convenience, we slightly changed the labelings: $\al\mapsto \rho_1$ and $\rho_i\mapsto \rho_{i+1}$, for $i=1,\dots,n$.}
\beq
\label{rel1}
\int_{\DM{0}{n+1}}c_{0,n+1}\left(e_{\rho_1}^\vee\otimes \bigotimes_{i=2}^{n+1} e_{\rho_i}\right)=c^{\rho_1}_{\rho_2\dots\rho_n},
\eeq
for any $1\leq \rho_1,\dots, \rho_{n+1}\leq N$ and any $n>1$. We claim that it is then possible to compute the numbers
\beq\label{intdtau}
\int_{\overline{D(\tau)}}c_{0,n+1}\left(e_{\rho_1}^\vee\otimes \bigotimes_{i=2}^{n+1} e_{\rho_i}\right),
\eeq
for all stable $(n+1)$-trees $\tau$. The homology classes $[\overline{D(\tau)}]$ span $H_\bullet(\DM{0}{n+1},\C)$, by Theorem \ref{mtspan}. 
Hence the datum of all possible numbers \eqref{intdtau}, for fixed indices $\rho_1,\dots,\rho_{n+1}\in\{1,\dots,N\}$, uniquely defines the cohomology class $c_{0,n+1}\left(e_{\rho_1}^\vee\otimes \bigotimes_{i=2}^{n+1} e_{\rho_i}\right)$ as linear functional on $H_\bullet(\DM{0}{n+1},\C)$. In other words, if we are able to compute all possible numbers \eqref{intdtau}, then any $F$-CohFT $(c_{0,n+1})_{n\geq 2}$ satisfying \eqref{rel1} is unique.
\vskip1mm
Given $\tau$, the number \eqref{intdtau} can be computed as follows, by iteration of the gluing property (4) of $F$-CohFT's. Denote by $v_o\in V_\tau$ the vertex of $\tau$ such that $1\in S_\tau(v_o)$, i.e. at which the tail 1 is attached. Orient the edges $e\in E_\tau$ in such a way that $v_o$ becomes an ``attractor''. In this way, at each vertex $v\in V_\tau\setminus\{v_o\}$ there is a single edge with outward orientation, all other edges being with inward orientation. At $v_o$ all edges are inward. Denote by $E_\tau^{\rm in}(v)$ the set of inward edges at $v$, and by $E_\tau^{\rm out}(v)$ the set of outward edges at $v$.
\vskip1mm
Consider now the following monomials attached to $\tau$. Each such monomial is product of coefficients $c$'s in \eqref{potf1}. In total we have $|V_\tau|$ factors $c$'s, one for each vertex $v\in V_\tau$. An index which is repeated inside the monomial (once upper and once lower) is said to be saturated. 
\begin{figure}
\centering
\def\svgscale{.9}
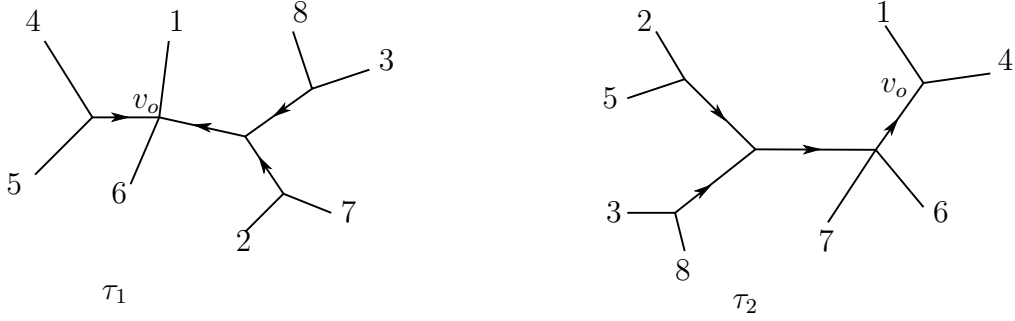
\caption{Two stable $8$-trees $\tau_1$ and $\tau_2$.}
\label{gr3}
\end{figure} 
The factor corresponding to $v\in V_\tau$ will have a total number of indices (upper and lower) equal to $|H_v|$, i.e. in bijection with half-edges. We will have 
\begin{enumerate}
\item a total number of $|E^{\rm out}_\tau(v)|\in\{0,1\}$ upper saturated indices, 
\item and a total number of $|E_\tau^{\rm in}(v)|$ lower saturated indices,
\item a total number of $|S_\tau(v)|$ lower indices selected from $(\rho_1,\dots, \rho_{n+1})$.
\end{enumerate}
The saturation of indices is dictated by edges: indices labelled by two halves of the same edge are saturated (one is up, the other is down). Non-saturated indices are dictated by the sets $S_\tau(v)$: the factor corresponding to $v\in V_\tau$ will have lower indices $\rho_{i}$ with $i\in S_\tau(v)$.
The vertex $v_o$ is the only vertex whose corresponding factor $c$ has upper index $\rho_1$.
\vskip1mm
The number \eqref{intdtau} equals the sum of all such monomials, over all possible values (ranging in $\{1,\dots, N\}$) of all saturated indices, according to Einstein's summation rule. For example, if $n=7$, and $\tau_1,\tau_2$ are the graphs of Figure \ref{gr3}, then
\[\int_{\overline{D(\tau_1)}}c_{0,8}\left(e_{\rho_1}^\vee\otimes\bigotimes_{j=2}^8e_{\rho_j}\right)=c^{\rho_1}_{\rho_6\al\bt}c^\al_{\rho_4\rho_5}c^\bt_{\gm\dl}c^\gm_{\rho_3\rho_8}c^\dl_{\rho_2\rho_7},
\]
\[\int_{\overline{D(\tau_2)}}c_{0,8}\left(e_{\rho_1}^\vee\otimes\bigotimes_{j=2}^8e_{\rho_j}\right)=c^{\rho_1}_{\rho_4 \al}c^\al_{\rho_7\rho_6\bt}c^\bt_{\gm\dl}c^\gm_{\rho_2\rho_5}c^\dl_{\rho_3\rho_8}.
\]
This is just an iteration of gluing rule (4) of an $F$-CohFT, which can be seen as a special instance of computation of \eqref{intdtau} for a one-edge $(n+1)$-tree. This proves uniqueness.
\vskip1mm
{\bf Existence. }In the previous part of the proof, we described an algorithm. For any fixed ${\bm \rho}=(\rho_1,\dots,\rho_{n+1})$, the algorithm associates with any stable $(n+1)$-tree $\tau$, a complex number $Y_{\bm \rho}(\tau)\in\C$, a polynomial expression in the coefficients $c$'s in \eqref{potf1}.  If we show that all linear relations between the homology classes $[\overline{D(\tau)}]$ are preserved by the map $\tau\mapsto Y_{\bm\rho}(\tau)$, then we would have a well-defined linear functional 
\[\tilde Y_{\bm \rho}\colon H_\bullet(\DM{0}{n+1},\C)\to \C,\qquad\overline{D(\tau)}\mapsto Y_{\bm \rho}(\tau),
\]i.e. a cohomology class. This would lead to a candidate as $F$-CohFT, 
\[
c_{0,n+1}\colon V^*\otimes V^{\otimes n}\to H^\bullet(\DM{0}{n+1},\C),\quad \left(e_{\rho_1}^\vee\otimes \bigotimes_{i=2}^{n+1}e_{\rho_i}\right)\mapsto \tilde Y_{\bm \rho}.
\]
\begin{figure}
\centering
\def\svgscale{.9}
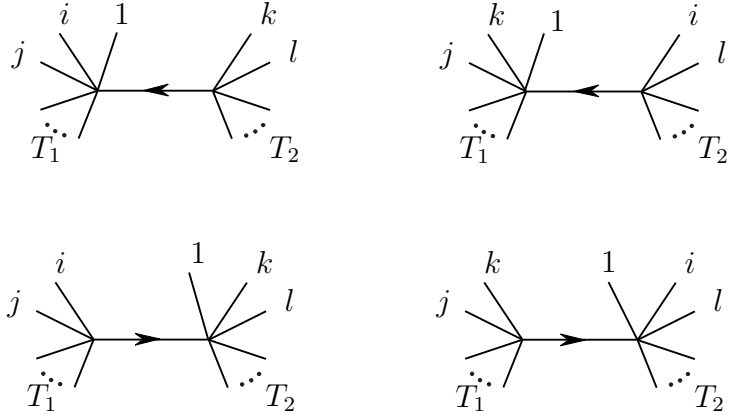
\caption{On the left (resp. right), we draw the two graphs $\tau'(\alpha)$ (resp. $\tau''(\alpha)$) which contribute to the l.h.s.\,\,(resp. r.h.s.) of equation \eqref{relY} in the case $v=v_o$.}
\label{gr4}
\end{figure} Indeed, the properties (1)-(3) of $F$-CohFT for $c_{0,n+1}$ would follow from the symmetry of $c^{\rho_1}_{\rho_2\dots\rho_{n+1}}$ in the lower indices, and equations \eqref{eqc1}, \eqref{eqc2}. Also, the gluing property (4) would follow from the definition of the numbers $Y_{\bm \rho}(\tau)$. This would complete the proof.
\vskip1mm
By Theorem \ref{mrel}, it is then sufficient to prove that, for any fixed system $(\tau,v,i,j,k,l)$, all the relations $R(\tau,v,i,j,k,l)$ are preserved by $Y_{\bm \rho}$, i.e.
\beq\label{relY}
\sum_\al Y_{\bm \rho}(\tau'(\al))=\sum_\bt Y_{\bm \rho}(\tau''(\bt)).
\eeq
The trees $\tau',\tau''$ are obtained from $\tau$ by replacing the vertex $v\in E_\tau$ with an edge. There are many ways to do this, labelled by 2-partitions of $H_\tau(v)$. They induced a 2-partition of ${i,j,k,l}$. We put on the l.h.s.\,\,of \eqref{relY} those which split $\{i,j,k,l\}$ in two pieces $\{i,j\}\coprod\{k,l\}$, and we put on the r.h.s.\,\,of \eqref{relY} those which split $\{i,j,k,l\}$ in two pieces $\{k,j\}\coprod\{i,l\}$. The remaining partitions do not contribute.
\vskip1mm
We have in total $2^5$ possible cases to consider, according wether $v$ coincides with the marked vertex $v_o$ or not, and wether each of $i,j,k,l$ is an edge or a tail.

Consider for example the case in which $v=v_o$ and each of $i,j,k,l$ is a tail. In the l.h.s.\,\,of \eqref{relY} we have the contributions coming from two possible graphs, according to the resulting position of the distinguished tail labeled by 1. Analogously, in the r.h.s.\,\,we have the contributions coming from two graphs. See Figure \ref{gr4}.

Equations \eqref{relY} thus reduces to an identity of the form
\[\sum_\al\left(\clr{ c^{\rho_1}_{\rho_i\rho_j\la\dots}c^\la_{\rho_k\rho_l\dots}}+c^{\rho_1}_{\rho_k\rho_l\la\dots}c^\la_{\rho_i\rho_j\dots}\right)=
\sum_\bt\left(\clr{ c^{\rho_1}_{\rho_j\rho_k\la\dots}c^\la_{\rho_i\rho_l\dots}}+c^{\rho_1}_{\rho_i\rho_l\la\dots}c^\la_{\rho_k\rho_j\dots}\right),
\]where dots stand for all possible partitions of indices, induced by $\al$ and $\bt$. Both red terms and black terms in this equation cancel, due to equations \eqref{eqc2}.
\vskip1mm
The reader can check that all other 31 possible cases can be handled similarly. One can always recognize in equation \eqref{relY} a linear combination of identities \eqref{eqc2}, whose left and right sides correspond to the inserted edge in $\tau'$ and $\tau''$, respectively.
\vskip1mm
This completes the proof.

\end{document}